\documentclass{amsart}

\usepackage{amssymb}
\usepackage{graphicx}
\usepackage{enumerate}
\usepackage{amscd}

\title{On genus-$1$ simplified broken Lefschetz fibrations}
\author{Kenta Hayano}
\address{Department of Mathematics, Graduate School of Science, Osaka University, Toyonaka, Osaka 560-0043, Japan}
\email{k-hayano@cr.math.sci.osaka-u.ac.jp}

\theoremstyle{plain}
\newtheorem{thm}{Theorem}[section]
\newtheorem{cor}[thm]{Corollary}
\newtheorem{lem}[thm]{Lemma}
\newtheorem{prop}[thm]{Proposition}

\theoremstyle{definition}
\newtheorem{defn}[thm]{Definition}

\newtheorem{conj}[thm]{Conjecture}
\newtheorem{rem}[thm]{Remark}

\makeatletter

\@addtoreset{figure}{section}
\makeatother

\makeatletter

\@addtoreset{equation}{section}
\makeatother

\begin{document}

\begin{abstract}

Auroux, Donaldson and Katzarkov introduced broken Lefschetz fibrations as a generalization of Lefschetz fibrations in order to describe near-symplectic $4$-manifolds. 
We first study monodromy representations of higher sides of genus-$1$ simplified broken Lefschetz fibrations. 
We then completely classify diffeomorphism types of such fibrations with connected fibers and with less than six Lefschetz singularities. 
In these studies, we obtain several families of genus-$1$ simplified broken Lefschetz fibrations, which we conjecture contain all such fibrations, 
and determine the diffeomorphism types of the total spaces of these fibrations. 
Our results are generalizations of Kas' classification theorem of genus-$1$ Lefschetz fibrations, 
which states that the total space of a non-trivial genus-$1$ Lefschetz fibration over $S^2$ is diffeomorphic to an elliptic surface $E(n)$, for some $n\geq 1$. 

\end{abstract}

\maketitle

\section{Introduction}

Broken Lefschetz fibrations were introduced by Auroux, Donaldson and Katzarkov \cite{ADK} as a generalization of Lefschetz fibrations. 
Donaldson \cite{Don} proved that symplectic $4$-manifolds admit Lefschetz fibrations after blow-ups and Gompf \cite{Gompf} proved the converse, 
i.e., the total space of every non-trivial Lefschetz fibration admits a symplectic structure. 
In \cite{ADK}, Auroux, Donaldson and Katzarkov extended these results to $4$-manifolds with near-symplectic structure, 
which is a closed $2$-form symplectic outside a union of circles where it vanishes transversely. 
In fact, they proved that a $4$-manifold $M$ is near-symplectic if and only if $M$ admits a broken Lefschetz pencil structure with certain conditions. 
So it is quite natural to try to generalize various results about Lefschetz fibrations to those for broken Lefschetz fibrations. 

\par

The classification of diffeomorphism types of total spaces $M$ of non-trivial genus-$1$ Lefschetz fibrations over $S^2$ has already been done by Kas \cite{Kas} and Moishezon \cite{Moish}, independently. 
He proved that $M$ is diffeomorphic to an elliptic surface $E(n)$, for some $n\geq 1$. 
Especially, the diffeomorphism type of $M$ is determined only by the Euler characteristic of $M$. 
Our main results generalize Kas and Moishezon's theorem to genus-$1$ simplified broken Lefschetz fibrations. 

\par

Williams \cite{Wil} proved that every closed connected oriented smooth $4$-manifold admits a simplified purely wrinkled fibration by using singularity theory 
(we will not mention the definition of simplified purely wrinkled fibrations in this paper. For precise definition of these fibrations, see \cite{Wil}). 
Simplified purely wrinkled fibrations can be changed into simplified broken Lefschetz fibrations 
by using Lemma 2.5 of \cite{Ba2} and Lekili's move in \cite{Lek} which changes a single cusped indefinite folds into indefinite folds and a Lefschetz singularity. 
As a result, it is known that every closed connected oriented smooth $4$-manifold admits a simplified broken Lefschetz fibration structure. 
So it is natural to examine the minimal genus of all simplified broken Lefschetz fibrations $M\rightarrow S^2$ on a fixed $4$-manifold $M$. 
As far as the author knows, there are few results about the above question so far. 
This problem is another motivation of our study and our main theorems give a partial answer to the question. 
Moreover, we will show that the minimal genus of all simplified broken Lefschetz fibrations on $\sharp k\mathbb{CP}^2$ is greater than $1$ for all $k>0$ (c.f. Corollary \ref{definite}). 

\par

Our main results is as follows. 

\noindent
{\bf Main Theorem A.} {\it 
The following $4$-manifolds admit genus-$1$ simplified broken Lefschetz fibration structures with non-empty round singular locus. 

\begin{itemize}

\item $\sharp k\mathbb{CP}^2\sharp \ell\overline{\mathbb{CP}^2}$, where $\ell>0$ and $k\geq 0$. 

\item $\sharp kS^2\times S^2$, where $k\geq 0$. 

\item $L\sharp k\overline{\mathbb{CP}^2}$, where $k\geq 0$ and $L$ is either of the manifolds $L_n$ and $L_n^{\prime}$ defined by Pao \cite{Pao}. 

\end{itemize}}
\noindent
\\[-10pt]
{\bf Main Theorem B.} {\it
Let $M$ be a total space of a genus-$1$ simplified broken Lefschetz fibration with non-empty round singular locus and with $r$ Lefschetz singularities, where $r$ is a non-negative integer. 
If $r\leq 5$, then $M$ is diffeomorphic to one of the following $4$-manifolds. 

\begin{itemize}

\item $\sharp k\mathbb{CP}^2\sharp(r-k)\overline{\mathbb{CP}^2}$, where $0\leq k \leq r-1$. 

\item $\sharp kS^2\times S^2$, where $k=\dfrac{r}{2}$ for even $r$. 

\item $S^1\times S^3\sharp S\sharp r\overline{\mathbb{CP}^2}$, where $S$ is either of the manifolds $S^2\times S^2$ and $S^2\tilde{\times}S^2$. 

\item $L\sharp r\overline{\mathbb{CP}^2}$, where $L$ is either of the manifolds $L_n$ and $L_n^{\prime}$. 

\end{itemize}}

\begin{rem}

In private talks with R. $\dot{\text{I}}$. Baykur, he told the author about an alternative proof of a part of Main Theorem A. 
His proof depends on singularity theory, while ours on Kirby calculus. 
The latter includes essential tools used in the proof of Main Theorem B. 

\end{rem}

While the author was writing this paper, Baykur and Kamada \cite{BK} posted a paper about genus-$1$ simplified broken Lefschetz fibrations on arXiv. 
They obtained some results about genus-$1$ simplified broken Lefschetz fibrations independent of the author. 
They classified genus-$1$ simplified broken Lefschetz fibrations up to blow-ups.  

\par

In Section 3, we examine monodromy representations of higher sides of genus-$1$ simplified broken Lefschetz fibrations by using a graphical method, called the chart description, 
which is used in \cite{KMMW} to determine monodromy representations of genus-$1$ Lefschetz fibrations. 
Although our method is similar to that of \cite{KMMW}, the definition of our chart is slightly different from that in \cite{KMMW}. 
We first give the definition of chart and review some basic techniques used in \cite{KMMW}. 
We then prove that the monodromy of the higher side of a genus-$1$ simplified broken Lefschetz fibration can be represented by a certain normal form after successive application of chart moves. 

\par

In Section \ref{examples}, we construct some families of genus-$1$ simplified broken Lefschetz fibrations by giving monodromy representations of higher sides of such fibrations. 
We then determine diffeomorphism types of total spaces of these fibrations by using Kirby calculus. 
Main Theorem A will be proved in the last of this section. 
We will also prove that $4$-manifolds with positive definite intersection form cannot admit genus-$1$ simplified broken Lefschetz fibration structures by using Kirby calculus. 
Some examples obtained in this section are also used to prove Main Theorem B in the last section. 

\par

We devote the last section to proving Main Theorem B. 
We prove that all genus-$1$ simplified broken Lefschetz fibrations with less than six Lefschetz singularities are contained in the families obtained in Section 4 
by moving certain word sequences representing monodromy representations. 
We end this section with explaining why the number of Lefschetz singularities is limited in Main Theorem B. 
\\[7pt]
{\bf Acknowledgments.} The author wishes to express his gratitude to Hisaaki Endo for his encouragement and many useful suggestions. 
He would also like to thank Refik $\dot{\text{I}}$nan\c{c} Baykur, Seiichi Kamada and Osamu Saeki for commenting on the draft of this paper and for helpful discussions. 
The author is supported by Yoshida Scholarship 'Master21' and he is deeply grateful to Yoshida Scholarship for their support.

\section{preliminaries}

\subsection{Broken Lefschetz fibrations}

Let $M$ and $B$ be compact connected oriented smooth manifolds of dimension $4$ and $2$, respectively. 

\begin{defn}

A smooth map $f:M\rightarrow B$ is called a {\it broken Lefschetz fibration} if it satisfies the following conditions: 

\begin{enumerate}[(1)]

\item $\partial M=f^{-1}(\partial B)$; 

\item there exist a finite set $\mathcal{C}=\{p_1,\ldots,p_n\}$ and a one dimensional submanifold $Z$ in int$M$ such that 
$df_p$ is surjective for all $p\in M\setminus(\mathcal{C}\cup Z)$; 

\item for each $p_i$, there exist a complex local coordinate $(z_1,z_2)$ around $p_i$ and a complex local coordinate $\xi$ around $f(p_i)$ such that 
$z_1(p_i)=z_2(p_i)=0$ and $f$ is locally written as $\xi=f(z_1,z_2)=z_1z_2$; 

\item for each $p\in Z$, there exist a real local coordinate $(t,x_1,x_2,x_3)$ around $p$ and a real local coordinate $(y_1,y_2)$ around $f(p)$ such that 
$Z=\{(t,0,0,0)|t\in\mathbb{R}\}$ and $f$ is locally written as $(y_1,y_2)=f(t,x_1,x_2,x_3)=(t,{x_1}^2+{x_2}^2-{x_3}^2)$; 

\item the restriction map of $f$ to $Z\cup\mathcal{C}$ is injective; 

\item for any $q\in B\setminus f(Z)$, $f^{-1}(q)$ has no $(-1)$-spheres. 

\end{enumerate}

We call a smooth map $f$ an {\it achiral broken Lefschetz fibration} if $f$ satisfies the conditions (1), (2), (4), (5), (6) and the following condition: 

\begin{enumerate}[(3)$^\prime$]

\item for each $p_i$, there exist a complex local coordinate $(z_1,z_2)$ around $p_i$ and a complex local coordinate $\xi$ around $f(p_i)$ such that 
$z_1(p_i)=z_2(p_i)=0$ and $f$ is locally written as $\xi=f(z_1,z_2)=z_1z_2$ or $z_1\overline{z_2}$. 

\end{enumerate}

\end{defn}

\begin{rem}

The above definition is rather special one from the previous definitions of broken Lefschetz fibrations. 
Indeed, all the definitions of broken Lefschetz fibrations in previous papers do not contain the condition (6). 
However, in this paper, we adopt the above definition. 

\end{rem}

If $Z=\phi$, we call $f$ simply a {\it (achiral) Lefschetz fibration}. 
A regular fiber of such a fibration is a closed oriented connected surface. 
We call its genus the {\it genus} of $f$. 
For simplicity, we will refer (achiral) broken Lefschetz fibrations and (achiral) Lefschetz fibrations by (A)BLF and (A)LF, respectively. 

\par

We call a singularity locally written as $f(z_1,z_2)=z_1z_2$ (resp. $f(z_1,z_2)=z_1\overline{z_2}$) a {\it Lefschetz singularity} (resp. {\it achiral Lefschetz singularity}). 
Also, we call $Z$ and $f(Z)$ in the definition above the {\it round singular locus} and the {\it round singular image}, respectively. 
An inverse image of $\nu f(Z)$ is called a {\it round cobordism} of $f$, where $\nu f(Z)$ is a tubelar neighborhood of $f(Z)$. 

\par

Let $f:M\rightarrow S^2$ be a BLF. 
Suppose that $f$ has a connected round singular locus. 
Then the round singular image of $f$ is an embedded circle and the circle divides $S^2$ into two $2$-disks $D_1$ and $D_2$. 
Moreover, one of the following occurs: 

\begin{itemize}

\item a regular fiber of the fibration $f^{-1}(D_i)\rightarrow D_i$ is connected, while a regular fiber of the fibration $f^{-1}(D_j)\rightarrow D_j$ is disconnected; 

\item any regular fibers are connected and a genus of a regular fiber of the fibration $f^{-1}(D_i)\rightarrow D_i$ is just one higher than a genus of a regular fiber of the fibration $f^{-1}(D_j)\rightarrow D_j$, 

\end{itemize}
where $\{i,j\}=\{1,2\}$. 
We call an inverse image of $D_i$ the {\it higher side} of $f$ and an inverse image of $D_j$ the {\it lower side} of $f$. 

\begin{defn}

A BLF $f:M\rightarrow S^2$ is {\it simplified} if $f$ satisfies the following conditions: 

\begin{enumerate}[(1)]

\item $f$ has a connected round singular locus; 

\item every Lefschetz singularity of $f$ is contained in the higher side of $f$; 

\item all fibers of $f$ are connected. 

\end{enumerate}

The genus of the higher side of $f$ as an LF is called the {\it genus} of $f$. 

\end{defn}

We sometimes refer to simplified BLF as SBLF for short.

\subsection{Monodromy representations}\label{monodromy}
We review the definition of monodromy representations.
For more details, see \cite{GS}.
Let $f:M\rightarrow B$ be a genus-$g$ LF, $\mathcal{C}\subset M$ the set of critical points of $f$ and 
$\psi_0:f^{-1}(y_0)\rightarrow \Sigma_g$ an orientation preserving diffeomorphism for a point $y_0\in B\setminus f(\mathcal{C})$. 
For a loop $\gamma:(I,\partial I)\rightarrow (B\setminus f(\mathcal{C}),y_0)$, the pull-back $\gamma^{\ast}f=\{(t,x)\in I\times M\mid \gamma(t)=f(x)\}$ of $f$ by $\gamma$ is a trivial $\Sigma_g$-bundle over $I$.
Let $\Psi:\gamma^{\ast}f\rightarrow I\times \Sigma_g$ be a trivialization of this bundle which is equal to $\psi_0$ on $\{0\}\times f^{-1}(y_0)$ and $\Psi(t,x)=(t,\psi_t(x))$.
Then $[\psi_1\circ \psi_0^{-1}]\in \mathcal{M}_g$ is independent of the choice of a trivialization $\Psi$, where $\mathcal{M}_g$ is the mapping class group of $\Sigma_g$. 
We define $\rho_{f}:\pi_1(B\setminus f(\mathcal{C}),y_0)\rightarrow \mathcal{M}_g$ as follows:
\[
\rho_{f}([\gamma])=[\psi_1\circ \psi_0^{-1}]. 
\]
This map is well-defined and called a {\it monodromy representation} of $f$. 
If we define the group multiplication of $\mathcal{M}_g$ as $[f]\cdot[g]=[g\circ f]$ for $[f],[g]\in \mathcal{M}_g$, then $\rho_{f}$ is a homomorphism. 
In this paper, the group multiplication of $\mathcal{M}_g$ is always defined as above. 
We remark that if we change the choice of a diffeomorphism $\psi_0$, 
we obtain a new monodromy representation by composing an inner automorphism by a fixed element $\psi\in\mathcal{M}_g$. 
Two monodromy representations $\rho_f:\pi_1(B\setminus f(\mathcal{C}),y_0)\rightarrow\mathcal{M}_g$ and $\rho_{f^\prime}:\pi_1(B^\prime\setminus f^\prime(\mathcal{C}^\prime),y_0^\prime)\rightarrow\mathcal{M}_g$ 
are {\it equivalent} if there exist an element $g\in\mathcal{M}_g$ and an orientation preserving diffeomorphism $h:(B,f(\mathcal{C}),y_0)\rightarrow(B^\prime,f^\prime(\mathcal{C}^\prime),y_0^\prime)$ 
such that 
\[
\text{conj}(g)\circ\rho_f=\rho_{f^\prime}\circ h_\ast, 
\]
where $\text{conj}(g)$ is the inner automorphism of $\mathcal{M}_g$ by the element $g$ and 
$h_\ast:\pi_1(B\setminus f(\mathcal{C}),y_0)\rightarrow\pi_1(B^\prime\setminus f^\prime(\mathcal{C}^\prime),y_0^\prime)$ is the isomorphism induced by $h$. 
\par

A monodromy representation of an LF over $D^2$ is written by a sequence of elements of the mapping class group as follows. 
We put $f(\mathcal{C})=\{y_1,\ldots,y_n\}$.
Let $A_1,\ldots,A_n$ be embedded paths in $D^2$, beginning at $y_0$ and otherwise disjoint, 
connecting $y_0$ to the respective critical values $y_1,\ldots,y_n$. 
We choose indices so that the paths $A_1,\ldots,A_n$ are cyclically ordered by traveling counterclockwise around $y_0$. 
We obtain an ordered basis $a_1,\ldots,a_n$ of $\pi_1(D^2\setminus f(\mathcal{C}),y_0)$ by connecting a counterclockwise circle around each $y_i$ to the base point $y_0$ by using $A_i$. 
Then the element $a_1\cdot\cdots\cdot a_n$ represents $[\partial D^2]$ in $\pi_1(D^2\setminus f(\mathcal{C}),y_0)$. 
We put $W_f=(\rho_f(a_1),\ldots,\rho_f(a_n))$. 
It is the sequence of elements of the mapping class group of $\Sigma_g$. 
We call this a {\it Hurwitz system of $f$} or a {\it monodromy factorization}. 
We put $w(W_f)=\rho_f(a_1)\cdot\cdots\cdot\rho_f(a_n)\in \mathcal{M}_g$. 

\par

\begin{figure}[htbp]
\begin{center}
\includegraphics[width=90mm]{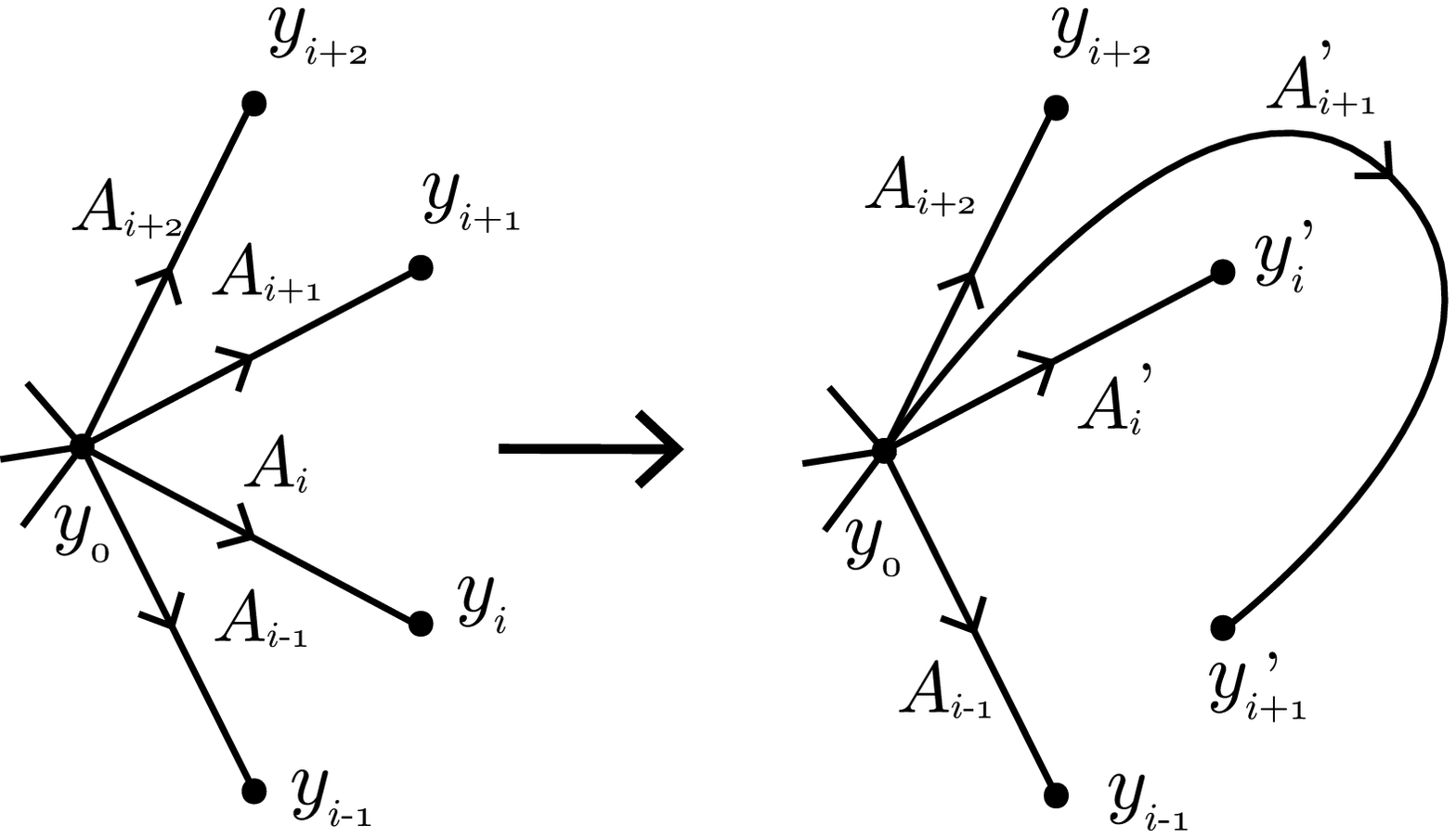}
\end{center}
\caption{}
\label{elem_trans}
\end{figure}

We obtain a sequence of elements of $\mathcal{M}_g$ from an LF $f$ over $D^2$ 
by choosing an identification $\psi_0$ and paths $A_1,\ldots,A_n$ as above. 
However, this sequence depends on these choices. 
We now review the effects of changing an identification and paths, see also \cite{GS}, \cite{Matsumoto2}. 
First, if we change an identification $\psi_0$, then all the elements of the sequence are conjugated by a fixed element $\psi$ of $\mathcal{M}_g$, 
i.e.$(\rho_f(a_1),\ldots,\rho_f(a_n))\rightarrow (\psi^{-1}\rho_f(a_1)\psi,\ldots,\psi^{-1}\rho_f(a_n)\psi)$. 
We call this transformation a {\it simultaneous conjugation} by $\psi$. 
Second, if we change paths in $D^2$ as Figure \ref{elem_trans}, 
then the pair of the elements $(\varphi_i,\varphi_{i+1})=(\rho_f(a_i),\rho_f(a_{i+1}))$ in $W_f$ is replaced by the pair $(\varphi_{i+1},\varphi_{i+1}^{-1}\cdot\varphi_i\cdot\varphi_{i+1})$. 
These transformations and their inverse are called {\it elementary transformations}. 
Lastly, if we change paths to other paths in $D^2$ satisfying the conditions mentioned above, 
then the sequence is changed by successive application of elementary transformations 
since any two paths satisfying the conditions are connected by successive application of isotopies in $D^2\setminus f(\mathcal{C})$ and the moves shown in Figure \ref{elem_trans}. 

From the above arguments, we obtain:

\begin{lem}

Let $f_i:M_i\rightarrow D^2$($i=1,2$) be a genus-$g$ LF and $W_{f_i}$ a Hurwitz system ($i=1,2$).
If $W_{f_1}$ is changed into $W_{f_2}$ by successive application of simultaneous conjugations and elementary transformations, 
then $f_1$ and $f_2$ are isomorphic.(i.e.there exist orientation preserving diffeomorphisms $\Phi:M_1\rightarrow M_2$ and $\varphi:D^2\rightarrow D^2$ such that $\varphi\circ f_1=f_2\circ\Phi$)

\end{lem}

Let $[f]$ be an element of $\mathcal{M}_1$. 
$[f]$ induces an isomorphism $f_\ast:H_1(T^2;\mathbb{Z})\rightarrow H_1(T^2;\mathbb{Z})$. 
We fix generators $\mu$ and $\lambda$ of $H_1(T^2;\mathbb{Z})$. 
Then $f$ is represented by a matrix $A\in SL(2,\mathbb{Z})$ as follows: 
\[
(f_\ast(\mu),f_\ast(\lambda))=(\mu,\lambda){}^tA. 
\]
The correspondence $[f]\mapsto A$ induces an isomorphism between the groups $\mathcal{M}_1$ and $SL(2,\mathbb{Z})$. 
We assume that $\mu\cdot\lambda=1$, where $\mu\cdot\lambda$ represents the intersection number of $\mu$ and $\lambda$. 
In this paper, we always identify the group $\mathcal{M}_1$ with the group $SL(2,\mathbb{Z})$ via the above isomorphism. 

Let $\gamma_{p,q}$ be a simple closed curve on $T^2$ representing the element $p\mu+q\lambda$ in $H_1(T^2;\mathbb{Z})$ for relatively prime integers $p,q$ 
and $T_{p,q}\in\mathcal{M}_1$ the right-handed Dehn twist along $\gamma_{p,q}$. 
Then, by the Picard-Lefschetz formula, 
\[
T_{1,0}=
\begin{pmatrix}
1 & 0 \\
1 & 1 
\end{pmatrix}, 
\hspace{2em}
T_{0,1}=
\begin{pmatrix}
1 & -1 \\
0 & 1
\end{pmatrix}. 
\]
If we put $X_1=
\begin{pmatrix}
1 & 0 \\
1 & 1 
\end{pmatrix}
,X_2=
\begin{pmatrix}
1 & -1 \\
0 & 1
\end{pmatrix}
$, then $SL(2,\mathbb{Z})$ has the following finite presentation \cite{MKS}:
\[
SL(2,\mathbb{Z})=<X_1,X_2|(X_1X_2)^6,X_1X_2X_1X_2^{-1}X_1^{-1}X_2^{-1}>. 
\]

Let $f:M\rightarrow S^2$ be a genus-$1$ SBLF, then the higher side of $f$ is a genus-$1$ LF over $D^2$.
So we can take a Hurwitz system of this fibration. 
We call this a {\it Hurwitz system} of $f$.

\subsection{Kirby diagrams of broken Lefschetz fibrations}\label{Kirby diagrams of BLF}
The technique of handle decomposition of BLFs was studied in \cite{Ba2}. 
We review this technique in this subsection. 
The reader should turn to \cite{Ba2} for the details of this technique and the examples of Kirby diagrams of BLFs.
We begin with the definition of ($4$-dimensional) round handles, which arises in the discussion of this subsection. 

\begin{defn}
Let $M$ be a smooth $4$-manifold and we put 
\[
R_{i}^{\pm}=I\times D^i\times D^{3-i}/((1,x_1,x_2,x_3)\sim(0,\pm x_1,x_2,\pm x_3)). 
\] 
where $i=1,2$.
Let $\psi:I\times \partial D^i\times D^{3-i}/\sim\rightarrow \partial M$ be an embedding. 
We call $M\bigcup_{\psi}R_{i}^{\pm}$ a $4$-manifold obtained by attaching a {\it round $i$-handle} and $R_{i}^{+}$(resp. $R_{i}^{-}$) {\it (4-dimensional) untwisted} (resp. {\it twisted}) {\it round $i$-handle}. 
\end{defn}

\begin{rem}
In \cite{Ba2}, round handles of arbitrary dimension are defined, but in this paper, only $4$-dimensional ones appear. 
Throughout this paper, we assume that round handles are always $4$-dimensional ones. 
Both untwisted and twisted round handles are diffeomorphic to $S^1\times D^3$, but these two round handles have different attaching regions. 
The attaching region of an untwisted round $i$-handle is the trivial $S^{i-1}\times D^{3-i}$-bundle over $S^1$, 
while the attaching region of a twisted round $i$-handle is a non-trivial $S^{i-1}\times D^{3-i}$-bundle over $S^1$. 
\end{rem}

The following lemma shows that a round handle attachment is described by attaching two handles with consecutive indices: 

\begin{lem}[\cite{Ba2}]\label{decomposition of round handles}
For $i\in\{1,2\}$, round $i$-handle attachment is given by $i$-handle attachment followed by $(i+1)$-handle attachment 
whose attaching sphere goes over the belt sphere of the $i$-handle geometrically twice, algebraically zero times if the round handle is untwisted and twice if the round handle is twisted.
\end{lem}

The next lemma shows relation between round handle attachment and handle decomposition of broken Lefschetz fibrations:

\begin{lem}[\cite{Ba2}]\label{round cobordisms}
Let $W$ be a round cobordism of a BLF with connected round singular locus $Z\subset W$ and 
$\partial_{-}W\subset \partial W$ the union of connected components of the boundary of $W$ facing to the lower side. 
Then $W$ is obtained by attaching an untwisted or twisted round $1$-handle to $\partial_{-}W\times I$. 
Moreover, the restriction of the BLF to $\psi(I\times\partial D^1\times\{0\}/\sim)\subset\partial_{-}W\times\{1\}$ is a double cover of the image of this map, 
where $\psi$ is the attaching map of the round handle. 
\end{lem}

\begin{rem}\label{dual decomp.}
We also obtain a round cobordism $W$ by a round $2$-handle attachment via dual handle decomposition. 
More precisely, let $\partial_{+}W\subset\partial W$ be the union of connected components of the boundary of $W$ facing to the higher side, 
then $W$ is obtained by attaching an untwisted or twisted round $2$-handle to $\partial_{+}W\times I$. 
Moreover, as above lemma, when the BLF restricts to $\psi(I\times\partial D^2\times\{0\}/\sim)\subset\partial_{+}W\times I$, 
this map is $S^1$-bundle over the image of the map, which is trivial if the round handle is untwisted and non-trivial if the round handle is twisted, 
where $\psi$ is the attaching map of the round handle. 
\end{rem}

This lemma says that a round cobordism is described by a round handle attachment. 
Conversely, let $f:M\rightarrow B$ be a BLF and $\tilde{M}$ the manifold obtained by attaching a round handle whose attaching map satisfies the condition mentioned above, 
then we can extend the map $f$ to a BLF $\tilde{f}:\tilde{M}\rightarrow B$ 
which has connected round singular locus in the round handle at the center of the core of the fiberwise attached $1$-handle (resp. $2$-handle) if the index of the round handle is $1$ (resp. $2$). 
\par
\begin{figure}[htbp]
\begin{center}
\includegraphics[width=60mm]{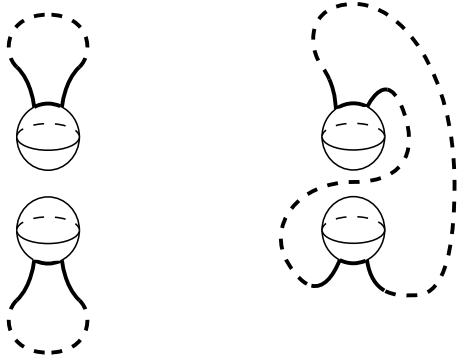}
\end{center}
\caption{the left diagram is the untwisted round handle, while the right diagram is the twisted one.}
\label{round1handle}
\end{figure}

We are ready to discuss Kirby diagrams of BLF. 
We begin with discussion of diagrams of round $1$-handles. 
By Lemma \ref{decomposition of round handles}, we can describe a diagram of a round $1$-handle as shown in Figure \ref{round1handle}. 
Moreover, by Lemma \ref{round cobordisms}, the bold curves in Figure \ref{round1handle} are double covers if a round $1$-handle describes a round cobordism of BLF. 

A round cobordism can be also described by a round $2$-handle attachment, which is given by $2$-handle attachment followed by $3$-handle attachment. 
By the condition mentioned in Remark \ref{dual decomp.}, 
round $2$-handle attachment to a BLF is realized as a fiberwise $2$-handle attachment 
if the round handle describes a round cobordism of an extended BLF. 
Let $H_2$ (resp. $H_3$) be a $2$-handle (resp. $3$-handle) of such a round $2$-handle. 
Then the attaching circle of $H_2$ is in a regular fiber of the surface bundle in the boundary of the BLF 
and its framing is along the regular fiber. 
Moreover, the attaching circle of $H_2$ is preserved under the monodromy of the surface bundle up to isotopy. 
As usual, we do not draw the $3$-handle $H_3$, which is forced to be attached in a way that it completes fiberwise $2$-handle attachment. 
Thus we only draw the attaching circle of the $2$-handle when we draw the round $2$-handle and the difference between the untwisted and twisted cases is somewhat implicit. 
However, we can distinguish these two cases by the action of the monodromy representation of the boundary surface bundle; 
the round $2$-handle is untwisted if this action preserves the orientation of the attaching circle and twisted if this action reverses the orientation of the attaching circle.

\section{Chart descriptions}
In \cite{KMMW}, chart descriptions were introduced to describe a monodromy representation of a genus-$1$ ALF. 
We modify the definition of a chart description to describe a monodromy representation of a higher side of a genus-$1$ SBLF 
and show that such a monodromy representation satisfies a kind of condition in this section. 
We begin with the modified definition of a chart. 

\begin{defn}
Let $\Gamma$ be a finite graph in $D^2$ (possibly being empty or having {\it hoops} that are closed edges without vertices). 
Then $\Gamma$ is called a {\it chart} if $\Gamma$ satisfies the following conditions: 

\begin{enumerate}[(1)]
\item the degree of each vertex is equal to $1$, $6$ or $12$; 

\item each vertex in $\partial D^2$ has degree-$1$; 

\item each edge in $\Gamma$ is labeled $1$ or $2$ and oriented; 

\item for a degree-$1$ vertex which is in $\text{int}D^2$, the incident edge is oriented inward (see Figure \ref{vertex}); 

\item for a degree-$6$ vertex, the six incident edges are labeled alternately with $1$ and $2$, 
and three consecutive edges are oriented inward and the other three edges are oriented outward (see Figure \ref{vertex}, where $\{i,j\}=\{1,2\}$); 

\item for a degree-$12$ vertex, the twelve incident edges are labeled alternately with $1$ and $2$, 
and all edges are oriented inward or all edges are oriented outward (see Figure \ref{vertex}, where $\{i,j\}=\{1,2\}$); 

\item an interior of each edge is in $\text{int}D^2$; 

\item let $\{v_1,\ldots,v_n\}$ be the set of vertices in $\partial D^2$. 
We assume that the indeces are chosen so that $v_1,\ldots,v_n$ appear in this order when we travel counterclockwise on $\partial D^2$. 
We define a pair $(i_k,\varepsilon_k)$ by the following rules: 

\begin{enumerate}[(i)]
\setlength{\itemindent}{8pt}
\item $i_k$ is a label of the edge $e$ whose end is $v_k$; 

\item $\varepsilon_k$ is $+1$ if $e$ is oriented away from $v$ and $-1$ if $e$ is oriented toward $v$. 
\end{enumerate}
\noindent
Then, for some $k$, the sequence $((i_k,\varepsilon_k),\ldots,(i_n,\varepsilon_n),(i_1,\varepsilon_1),\ldots,(i_{k-1},\varepsilon_{k-1}))$ consists of following subsequences: 

\begin{enumerate}[(a)]
\setlength{\itemindent}{14pt}
\item $((1,\varepsilon))$; 

\item $((i,\varepsilon),(j,\varepsilon),(i,\varepsilon),(j,\varepsilon),(i,\varepsilon),(j,\varepsilon))$  $(\{i,j\}=\{1,2\})$, 

\end{enumerate}
where $\varepsilon$ is equal to $\pm 1$. 

\end{enumerate}

\begin{figure}[htbp]
\begin{center}
\includegraphics[width=140mm]{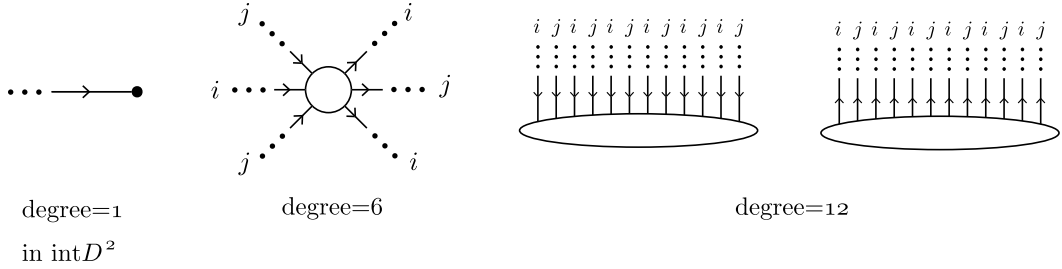}
\end{center}
\caption{verteces of a chart.}
\label{vertex}
\end{figure}

\end{defn}

\begin{figure}[htbp]
\begin{center}
\includegraphics[width=65mm]{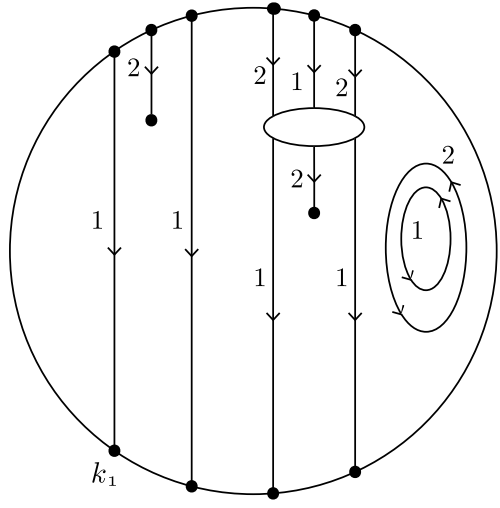}
\end{center}
\caption{An example of a chart}
\label{charts}
\end{figure}

\begin{rem}

In this paper, we will use ovals as vertices to make the diagrams easier to understand as shown in Figure \ref{vertex} and Figure \ref{charts}, for example. 
The reader may confuse these ovals with hoops. 
However, we can distinguish these two diagrams easily since hoops have orientation and labels while ovals representing vertices do not. 
For example, the chart in Figure \ref{charts} has two hoops and one oval representing degree-$6$ vertex. 

\end{rem}

An example of a chart is illustrated in Figure \ref{charts}. 
For this chart, the sequence mentioned in the condition (8) of the definition of charts is as follows: 
\[
((1,-1),(1,-1),(1,-1),(1,-1),(2,+1),(1,+1),(2,+1),(1,+1),(2,+1),(1,+1)), 
\]
which satisfies the condition mentioned in the definition. 

For a chart $\Gamma$, we denote by $V(\Gamma)$ the set of all the vertices of $\Gamma$, 
and by $S_\Gamma$ the subset of $V(\Gamma)$ consisting of the degree-$1$ vertices in $\text{int}D^2$. 
Let $v$ be a vertex of $\Gamma$. 
An edge $e$ incident to $v$ is called an {\it incoming edge} of $v$ if $e$ is oriented toward $v$ and an {\it outgoing edge} of $v$ if $e$ is oriented away from $v$. 

\par

A degree-$1$ or $12$ vertex of a chart is {\it positive} (resp.{\it negative}) if all the edges incident to the vertex is outgoing edge (resp. incoming edge) of the vertex. 
We remark that each degree-$1$ vertex in $\text{int}D^2$ is negative by the definition of charts. 

\par

Among the six edges incident to a degree-$6$ vertex $v$ of a chart, three consecutive edges are incoming edges of $v$ and the other three edges are outgoing edges of $v$. 
We call the middle edge of the three incoming edges or the three outgoing edges a {\it middle edge} and another edge a {\it non-middle edge}. 

\par

An edge in a chart is called a {\it $(d_1,d_2)$-edge} if its end points are a degree-$d_1$ vertex and a degree-$d_2$ vertex, where $d_1,d_2\in\{1,6,12\}$ and $d_1\leq d_2$. 
An edge in a chart is called a {\it $(\partial,d)$-edge} if one of its end points is in $\partial D^2$ and the other is degree-$d$ vertex, where $d\in\{1,6,12\}$, 
and we call an edge whose two end points are in $\partial D^2$ a {\it $(\partial,\partial)$-edge}. 
A $(\partial,\ast)$-edge is called a {\it boundary edge}, where $\ast\in\{1,6,12,\partial\}$. 

\par

For a chart $\Gamma$, let $((i_1,\varepsilon_1),\ldots,(i_n,\varepsilon_n))$ be a sequence determined by the rule mentioned in the condition (8) of the definition of charts. 
We assume that indices are chosen so that this sequence consists of the two subsequences; 
$((1,\varepsilon))$ and $((i,\varepsilon),(j,\varepsilon),(i,\varepsilon),(j,\varepsilon),(i,\varepsilon),(j,\varepsilon))$, where $\{i,j\}=\{1,2\}$ and $\varepsilon=\pm 1$. 
We call such a sequence a {\it boundary sequence} of $\Gamma$ and two subsequences above the {\it unit subsequences}. 
For a fixed decomposition of the sequence into the unit subsequences, 
the union of six vertices which corresponds to a subsequence $((i,\varepsilon),(j,\varepsilon),(i,\varepsilon),(j,\varepsilon),(i,\varepsilon),(j,\varepsilon))$ and the six edges incident to the vertices is called a {\it boundary comb} of $\Gamma$ with respect to the fixed decomposition.

\par

Let $\Gamma$ be a chart in $D^2$. 
A path $\eta:[0,1]\rightarrow D^2$ is said to be {\it in general position} with respect to $\Gamma$ 
if $\eta([0,1])\cap\Gamma$ is empty or consists of finite points in $\Gamma\setminus V(\Gamma)$ and $\eta$ intersects edges of $\Gamma$ transversely. 
Let $\eta$ be such a path. 
We put $\eta([0,1])\cap\Gamma=\{p_1,\ldots,p_n\}$. 
We assume that $p_1,\ldots,p_n$ appear in this order when we go along $\eta$ from $\eta(0)$ to $\eta(1)$. 
For each $p_i$, we define a letter $w_i=X_k^\varepsilon$ by following rules: 

\begin{enumerate}[(i)]

\item $k$ is the label of the edge of $\Gamma$ containing the point $p_i$; 

\item $\varepsilon$ is equal to $+1$ if the intersection number $\eta\cdot(\text{the edge containing }p_i)$ is $+1$ and $-1$ if the intersection number $\eta\cdot(\text{the edge containing }p_i)$ is $-1$.

\end{enumerate}

\noindent
We put $w_\Gamma(\eta)=w_1,\ldots,w_n$ and call this word the {\it intersection word} of $\eta$ with respect to $\Gamma$. 
We assume that this word represents an element of $SL(2,\mathbb{Z})$ by regarding the letters $X_1$, $X_2$ as the matrices defined in subsection \ref{monodromy}. 

\begin{defn}
Let $\Gamma$ be a chart in $D^2$. 
We fix a point $y_0\in D^2\setminus V(\Gamma)$. 
We define a homomorphism
\[
\rho_\Gamma:\pi_1(D^2\setminus S_\Gamma,y_0)\rightarrow SL(2,\mathbb{Z})
\]
as follows: For an element $x\in\pi_1(D^2\setminus S_\Gamma,y_0)$, we choose a representative path
\[
\eta:[0,1]\rightarrow D^2\setminus S_\Gamma
\]
of $x$ so that $\eta$ is in general position with respect to $\Gamma$. 
Then we put $\rho_\Gamma(x)=w_\Gamma(\eta)$. 
We call the homomorphism $\rho_\Gamma$ the {\it monodromy representation associated with $\Gamma$}. 
\end{defn}

We can prove the following lemma by an argument similar to that given in the proof of Lemma 12 of \cite{KMMW}. 

\begin{lem}
The homomorphism $\rho_\Gamma:\pi_1(D^2\setminus S_\Gamma,y_0)\rightarrow SL(2,\mathbb{Z})$ is well-defined. 
\end{lem}

Since the monodromy of the boundary of the higher side of a genus-$1$ SBLF fixes a simple closed curve in a regular fiber, 
the monodromy is represented by an element $\pm X_1^m$, where $m$ is an integer. 
So we can also prove the following lemma by the same argument as that in the proof of Theorem 15 of \cite{KMMW}. 

\begin{lem}\label{relation between BLF and chart}
Let $f:M\rightarrow S^2$ be a genus-$1$ SBLF. 
Then there exists a chart $\Gamma$ in $D^2$ such that 
the monodromy representation of the higher side of $f$ is equal to the monodromy representation associated with $\Gamma$ up to inner automorphisms of $SL(2,\mathbb{Z})$. 
\end{lem}

By Lemma \ref{relation between BLF and chart}, the monodromy representation of the higher side of a genus-$1$ SBLF is represented by a chart. 
However, such a chart is not unique. 
We next introduce moves of charts which do not change the associated monodromy representations. 

\begin{lem}\label{CI}
Let $\Gamma_1$ and $\Gamma_2$ be charts and $E\subset D^2$ a $2$-disk. 
We assume the following conditions:

\begin{enumerate}[(a)]
\setlength{\itemindent}{5pt}
\item $E\cap S_{\Gamma_i}=\emptyset$ ($i=1,2$); 

\item $\Gamma_1$ and $\Gamma_2$ are identical outside of $E$; 

\item $D^2\setminus E$ is path connected. 

\end{enumerate}

\noindent
Then the monodromy representation associated with $\Gamma_1$ is equal to the one associated with $\Gamma_2$. 

\end{lem}

We can prove Lemma \ref{CI} similarly to the proof of Lemma 16 in \cite{KMMW}. 

\begin{defn}
When two charts $\Gamma_1$ and $\Gamma_2$ are in the situation of Lemma \ref{CI}, we say that $\Gamma_1$ is obtained from $\Gamma_2$ by a {\it CI-move} in $E$. 
In particular, a CI-move described in Figure \ref{channel change} is called a {\it channel change}. 

\begin{figure}[htbp]
\begin{center}
\includegraphics[width=70mm]{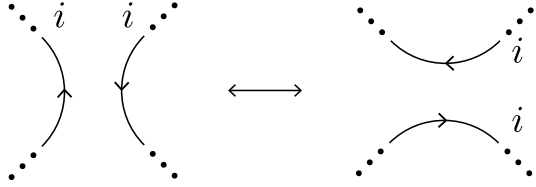}
\end{center}
\caption{a channel change}
\label{channel change}
\end{figure}

\end{defn}

\begin{figure}[htbp]
\begin{center}
\includegraphics[width=130mm]{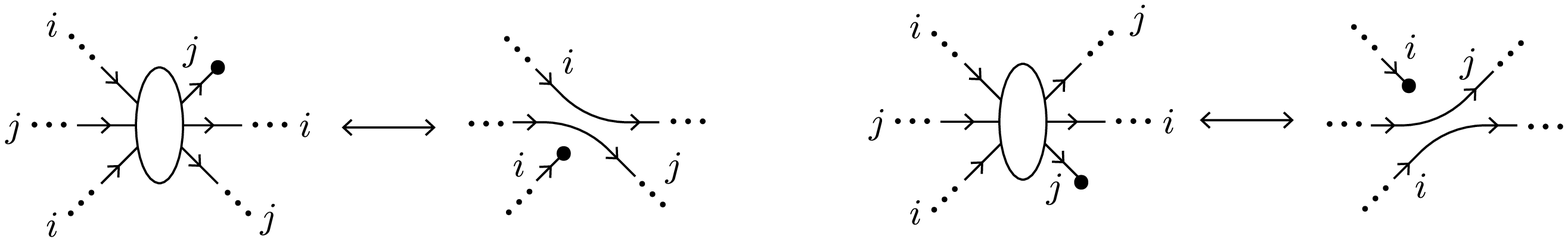}
\end{center}
\caption{CII-moves}
\label{CII-move}
\end{figure}

\begin{lem}\label{CII}

Let $\Gamma_1$ and $\Gamma_2$ be charts and $E\subset D^2$ a $2$-disk. 
We assume that the following conditions:

\begin{enumerate}[(a)]
\setlength{\itemindent}{5pt}
\item $\Gamma_1$ and $\Gamma_2$ are differ by one of Figure \ref{CII-move} in $E$; 

\item $\Gamma_1$ and $\Gamma_2$ are identical outside of $E$; 

\item $D^2\setminus E$ is path connected. 

\end{enumerate}
\noindent
Then the monodromy representation associated with $\Gamma_1$ is equal to the one associated with $\Gamma_2$ up to an equivalence. 

\end{lem}

The proof of Lemma \ref{CII} is quite similar to that of Lemma 18 of \cite{KMMW}. 
So we omit it. 

\begin{defn}

When two charts $\Gamma_1$ and $\Gamma_2$ are in the situation of Lemma \ref{CII}, we say that $\Gamma_1$ is obtained from $\Gamma_2$ by a {\it CII-move} in $E$. 

\end{defn}

\begin{rem}

CII-moves defined in \cite{KMMW} have four types of substitution (see Figure 6 of \cite{KMMW}), 
while CII-moves we defined have only two types of substitution illustrated in Figure \ref{CII-move}. 
This is because charts defined in \cite{KMMW} may have positive degree-$1$ vertices in $\text{int}D^2$ but our charts never have by the definition. 

\end{rem}

By a {\it C-move}, we mean a CI-move, CII-move or an isotopic deformation in $D^2$. 
Two charts are {\it C-move equivalent} if they are related by a finite sequence of C-moves. 
The monodromy representations associated with such charts is equivalent by Lemma \ref{CI} and Lemma \ref{CII}. 

\begin{thm}\label{main in chart sec}
Let $f:M\rightarrow S^2$ be a genus-$1$ SBLF and $W_f$ a Hurwitz system of $f$. 
Then by successive application of simultaneous conjugations and elementary transformations, we can change $W_f$ into a sequence $W$ which satisfies the following conditions: 
\begin{enumerate}[(a)]
\setlength{\itemindent}{10pt}
\item $w(W)=\pm X_1^m$; 

\item $W=(X_1,\ldots,X_1,X_1^{-n_1}X_2X_1^{n_1},\ldots,X_1^{-n_s}X_2X_1^{n_s})$, 

\end{enumerate}
\noindent
where $m,n_1,\ldots,n_s$ are integers. 

\end{thm}

\begin{lem}\label{vanish12}

Let $\Gamma$ be a chart. Then by successive application of C-moves, we can change $\Gamma$ into a chart which has no degree-$12$ vertices. 

\end{lem}

\begin{figure}[htbp]
\begin{center}
\includegraphics[width=110mm]{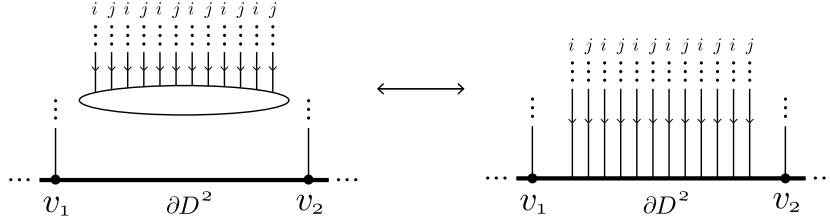}
\end{center}
\caption{CI-move used to prove Lemma \ref{vanish12}, 
where the bold lines represent $\partial D^2$. }
\label{moves12}
\end{figure}

{\it Proof}. We first remark that a chart move illustrated in Figure\ref{moves12} is a CI-move, 
where two vertices $v_1$ and $v_2$ satisfy one of the following conditions: 

\begin{enumerate}[(a)]

\item $v_1$ and $v_2$ are not contained in same boundary comb; 

\item $v_1$ and $v_2$ are end points of a boundary comb. 

\end{enumerate}

We choose a decomposition of the boundary sequence of $\Gamma$ into the unit subsequences. 
Let $v_1$ and $v_2$ be consecutive vertices in $\partial D^2$ which satisfy one of the conditions (a) and (b) as above. 
We denote by $S$ the connected component of $\partial D^2\setminus(\partial D^2\cap\Gamma)$ between $v_1$ and $v_2$. 
We can move all the degree-$12$ vertices in $\Gamma$ into a region of $\partial D^2\setminus \Gamma$ containing $S$ 
by using CI-moves illustrated in Figure 12 of \cite{KMMW}. 
Then, by the CI-move illustrated in Figure \ref{moves12}, we can eliminate all the degree-$12$ vertices in $\Gamma$. 
This completes the proof of Lemma \ref{vanish12}. \hfill $\square$

\begin{lem}\label{about(1,6)}

Let $\Gamma$ be a chart. 
Then, by successive application of C-moves, 
we can change $\Gamma$ into a chart $\Gamma^\prime$ such that each $(1,6)$-edge $e$ in $\Gamma^\prime$ satisfies the following conditions: 

\begin{enumerate}[(i)]
\item $e$ is a middle edge; 

\item the label of $e$ is $2$; 

\item let $K$ be the connected component of $D^2\setminus\Gamma^\prime$ whose closure contains $e$. 
Then $K\cap\partial D^2$ is not empty. 

\end{enumerate}

\end{lem}

The idea of the proof of Lemma \ref{about(1,6)} is similar to the proof of Lemma 22 in \cite{KMMW}. 
But the two proofs are slightly different because of the difference of the definition of charts. 
So we give the full proof below. 
\\
\par

{\it Proof}. 
Let $n(\Gamma)$ be the sum of the number of degree-$6$ vertices and the number of $(1,6)$-edges in $\Gamma$. 
The proof proceeds by induction on $n(\Gamma)$. 

\par

If $n(\Gamma)=0$, the conclusion of Lemma \ref{about(1,6)} holds since $\Gamma$ has no $(1,6)$-edges. 
We assume that $n(\Gamma)>0$ and there exists a $(1,6)$-edge which does not satisfy at least one of the conditions {\it (i)}, {\it (ii)} or {\it (iii)} of Lemma \ref{about(1,6)}. 

\par

{\it Case.1}: Suppose that $\Gamma$ has a non-middle $(1,6)$-edge. 
Let $v$ be a degree-$6$ vertex which is an end point of a $(1,6)$-edge. 
We can apply a CII-move around $v$ and eliminate this vertex. 
Then one of the following occurs: 

\begin{itemize}

\item both the number of degree-$6$ vertices and the number of $(1,6)$-edges decrease; 

\item the number of $(1,6)$-edges is unchanged, but the number of degree-$6$ vertices decreases; 

\end{itemize}
\noindent
In each case, $n(\Gamma)$ decreases and the conclusion holds by the induction hypothesis. 

\par

{\it Case.2}: Suppose that $\Gamma$ has a middle $(1,6)$-edge whose label is $1$. 
Let $e$ be the edge and $v_0$ and $v_1$ the end points of $e$ whose degrees are $1$ and $6$, respectively. 
We denote by $K$ the connected component of $D^2\setminus\Gamma$ whose closure contains $v_0$, $v_1$ and $e$. 
We take a sequence $f_1,\ldots,f_m$ of edges of $\Gamma$ with signs as in the proof of Lemma 22 of \cite{KMMW}. 
For each $f_i$, we take a letter $w(f_i)=X_k^{\varepsilon}$, where $k$ is equal to the label of the edge $f_i$ 
and $\varepsilon$ is equal to the sign of $f_i$. 
We remark that both $f_1$ and $f_m$ are equal to $e$ and the sign of $f_1$ is negative, while the sign of $f_m$ is positive, since the vertex $v_0$ is negative. 

\par

{\it Case.2.1}: There exists a consecutive pair $f_i$ and $f_{i+1}$ such that two edges are incident to a common vertex and
\[
(w(f_i),w(f_{i+1}))=(X_1^{-1},X_2^{-1}). 
\]

{\it Case.2.2}: There exists a consecutive pair $f_i$ and $f_{i+1}$ such that two edges are incident to a common vertex and
\[
(w(f_i),w(f_{i+1}))=(X_2,X_1). 
\]

{\it Case.2.3}: $K\cap\partial D^2=\emptyset$. 

If one of the above three cases occurs, then the conclusion holds by the same argument as that in Lemma 22 of \cite{KMMW}. 

\begin{figure}[htbp]
\begin{center}
\includegraphics[width=145mm]{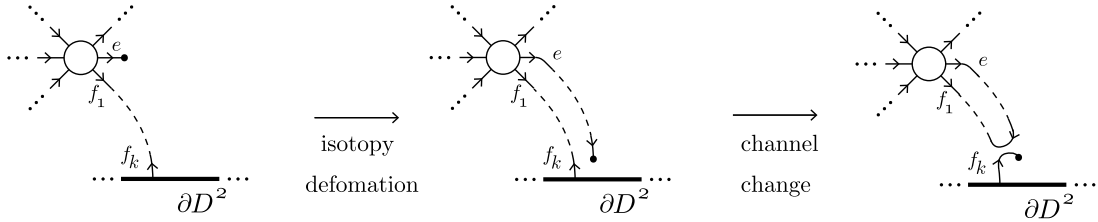}
\end{center}
\caption{The bold line in the figure describes $\partial D^2$. }
\label{move_case24.1}
\end{figure}

{\it Case.2.4}: Suppose that $K\cap\partial D^2\neq\emptyset$. 
Then one of the edges $f_1,\ldots,f_m$ is a boundary edge. 
By Cases 2.1 and 2.2, we can assume that $(w(f_i),w(f_{i+1}))$ is not equal to the subsequences $(X_1^{-1},X_2^{-1})$ and $(X_2,X_1)$ if $f_i$ and $f_{i+1}$ are incident to a common vertex. 
Let $f_k$ be a boundary edge with the smallest index. 
By the assumption above, $w(f_k)$ is equal to either $X_1^{-1}$ or $X_2$. 
If $w(f_k)=X_1^{-1}$, we can decrease the number of $(1,6)$-edges by applying C-moves illustrated in Figure \ref{move_case24.1}. 
So the conclusion holds by the induction hypothesis. 

Suppose that $w(f_k)=X_2$. 
One of $f_{k+1},\ldots,f_m$ is a boundary edge but not a $(\partial,1)$-edge. 
Let $f_l$ be such an edge with the smallest index. 

\par

\begin{figure}[htbp]
\begin{center}
\includegraphics[width=145mm]{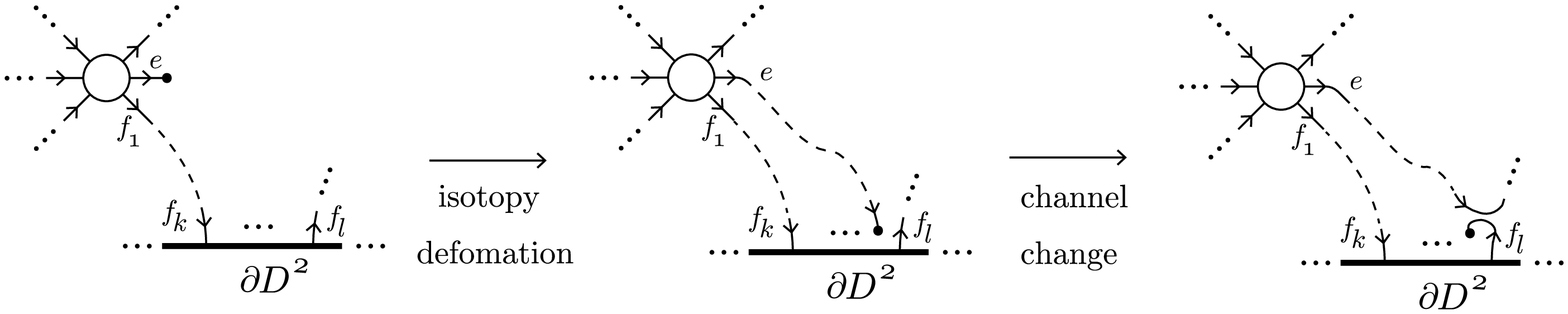}
\end{center}
\caption{}
\label{move_case24.2}
\end{figure}

{\it Case.2.4.1}: Suppose that $w(f_l)=X_1$. 
Then we can decrease the number of $(1,6)$-edges by applying C-moves illustrated in Figure\ref{move_case24.2} 
and the conclusion holds by the induction hypothesis. 

\begin{figure}[htbp]
\begin{center}
\includegraphics[width=145mm]{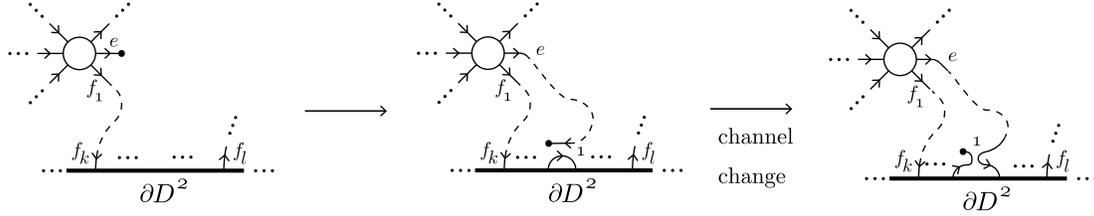}
\end{center}
\caption{We first apply CI-move between the two boundary comb which contain $f_k$ and $f_l$, respectively, and we obtain a new $1$-labeled $(\partial,\partial)$-edge. 
Then we move $v_0$ near this edge by isotopy deformation and apply a channel change. }
\label{move_case24.3}
\end{figure}

{\it Case.2.4.2}: Suppose that $w(f_l)=X_2$. 
When we fix a decomposition of the boundary sequence of $\Gamma$ into unit subsequences, 
$f_k$ is contained in a boundary comb distinct from that of $f_l$. 
So we can apply C-moves as shown in Figure \ref{move_case24.3} and conclusion holds by induction. 

{\it Case.2.4.3}: Suppose that $w(f_l)=X_2^{-1}$. 
If both $f_k$ and $f_l$ were contained in a same boundary comb, there would be at least one $(\partial,1)$-edge between $f_k$ and $f_l$. 
But all the degree-$1$ edges are negative. 
This contradiction says that a boundary comb that contains $f_k$ is different to a boundary comb that contains $f_l$. 
So we can apply C-moves similar to the C-moves we use in Case.2.4.2 and the conclusion holds by induction hypothesis. 

{\it Case.2.4.4}: Suppose that $w(f_l)=X_1^{-1}$. 
If each $f_{l+1},\ldots,f_{m}$ were not a boundary edge, then, 
\[
(w(f_{l}),\ldots,w(f_{m}))=(X_1^{-1},X_2,X_1^{-1},X_2,\ldots)
\]
This contradicts $w(f_m)=X_1$. 
So at least one of $f_{l+1},\ldots,f_{m}$ is a boundary edge. 
Let $f_{k^\prime}$ be such an edge with smallest index. 
Then $w(f_{k^\prime})$ is equal to either $X_1^{-1}$ or $X_2$. 
If $w(f_{k^\prime})=X_1^{-1}$, the conclusion holds by the above argument. 
If $w(f_{k^\prime})=X_2$, one of four cases above occurs for $f_{k^\prime}$. 
When one of the former three cases occurs, the conclusion holds by the same argument. 
When Case.2.4.4 occurs for $f_{k^\prime}$, we can take $f_{k^{\prime\prime}}$ again as we take $f_{k^\prime}$. 
We can repeat the above argument and thus the conclusion holds since $m$ is finite. 

\par

{\it Case.3}: Suppose that $\Gamma$ has a middle $(1,6)$-edge whose label is $2$. 
We define $K$ as we defined in Case.2. 
If $K\cap\partial D^2=\emptyset$, we can prove the conclusion by the same argument as that in Cases 2.1, 2.2 and 2.3. 
So we conclude that $K\cap\partial D^2\neq\emptyset$. 

\par

Combining the conclusions of Cases.1, 2 and 3, we complete the proof of Lemma \ref{about(1,6)}.  \hfill $\square$
\\
\par

\begin{figure}[htbp]
\begin{center}
\includegraphics[width=65mm]{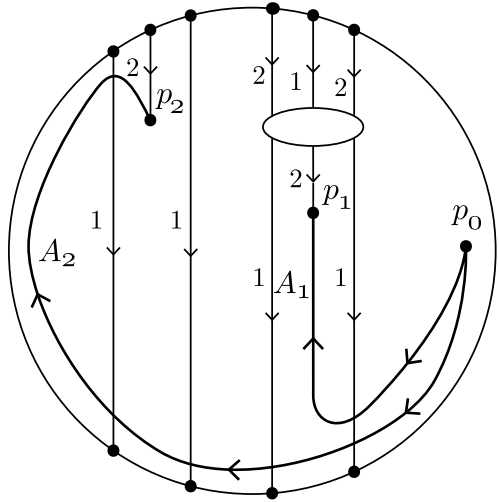}
\end{center}
\caption{Examples of paths $A_1,\ldots,A_n$ determined by the condition (a) and the constructions (b) and (c). }
\label{ex_path}
\end{figure}

{\it Proof of Theorem \ref{main in chart sec}}: 
By Lemma \ref{relation between BLF and chart}, 
we can take a chart $\Gamma$ that $\rho_\Gamma$ is equal to the monodromy representation of the higher side of $f$ up to inner automorphisms of $SL(2,\mathbb{Z})$. 
We first eliminate degree-$12$ vertices in $\Gamma$ by applying Lemma \ref{vanish12}. 
Then, by applying Lemma \ref{about(1,6)}, we change the chart $\Gamma$ into a chart such that 
all the $(1,6)$-edges satisfy conditions {\it (i)}, {\it (ii)} and {\it (iii)} in Lemma \ref{about(1,6)}. 
In the process to prove Lemma \ref{about(1,6)}, no new degree-$12$ vertices appear. 
So the chart obtained by the above process has no degree-$12$ vertices. 
Let $\{v_1,\ldots,v_m\}$ be the set of degree-$1$ vertices of $\Gamma$ in $\partial D^2$. 
We choose the indices of $v_i$ so that $v_1,\ldots,v_m$ appear in this order when we travel counterclockwise on $\partial D^2$ 
and that $v_1$ and $v_m$ satisfy one of the following conditions: 

\begin{itemize}

\item $v_1$ and $v_m$ are not contained in same boundary comb; 

\item $v_1$ and $v_m$ are end points of a boundary comb. 

\end{itemize}
\noindent
We denote by $e_i$ a boundary edge incident to $v_i$.
We put $S_{\Gamma}=\{p_1,\ldots,p_n\}$. 
Let $K_i$ be a connected component of $D^2\setminus\Gamma$ whose closure contains $p_i$. 
By the assumption about $\Gamma$, each $p_i$ is an end point of either $(1,6)$-edge or $(\partial,1)$-edge. 
For each $p_i$ which is an end point of $(1,6)$-edge, we choose a connected component $E_i$ of $K_i\cap\partial D^2$. 
We denote the two points of $\partial E_i$ by $v_{k_i}$ and $v_{k_{i}+1}$, where $k_i\in\{1,\ldots,m\}$ and $v_{m+1}=v_1$. 
Let $V$ be a sufficiently small collar neighborhood of $\partial D^2$ in $D^2$ and $p_0$ a point in $V\cap K$, 
where $K$ is a connected component of $D^2\setminus \Gamma$ whose closure contains a connected component of $\partial D^2\setminus\{v_1,\ldots,v_m\}$ between $v_m$ and $v_1$. 
We take embedded paths $A_i$ ($i=1,\ldots,n$) in $D^2$ starting from $p_0$ as follows: 

\begin{enumerate}[(a)]

\item if $i\neq j$, then $A_i\cap A_j=\{p_0\}$; 

\item if $p_i$ is incident to a $(\partial,1)$-edge $e_j$, then $A_i$ starts from $p_0$, 
travels in $V$ across the edges $e_1,\ldots,e_{j-1}$, goes into $K_i$ and ends at $p_i$; 

\item if $p_i$ is incident to a $(1,6)$-edge, then $A_i$ starts from $p_0$, 
travels in $V$ across the edges $e_1,\ldots,e_{k_{i}-1}$, goes into $K_i$ and ends at $p_i$. 

\end{enumerate}
\noindent
For example, the paths $A_1,\ldots,A_n$ are as shown in Figure\ref{ex_path} for the charts described in Figure\ref{charts}. 

Let $a_i$ be an element of $\pi_1(D^2\setminus S_\Gamma,p_0)$ 
which is represented by a curve obtained by connecting counterclockwise circle around $p_i$ to the base point $p_0$ by using $A_i$. 
It is sufficient to prove that each $\rho_\Gamma(a_i)$ is equal to either $X_1$ or $X_1^{-n}X_2X_1^n$, where $n$ is an integer. 

\par

{\it Case.1}: Suppose that $p_i$ is an end point of $(1,6)$-edge and $e_{k_{i}}$ is not contained in a boundary comb. 
Then the intersection word of $A_i$ is equal to $X_1^{n}$. 
So $\rho_\Gamma(a_i)$ is equal to $X_1$ if the label of the $(1,6)$-edge is $1$ and $X_1^nX_2X_1^{-n}$ if the label of the $(1,6)$-edge is $2$. 

\par

{\it Case.2}: Suppose that $p_i$ is an end point of $(\partial,1)$-edge and the edge is not contained in a boundary comb. 
Then the intersection word of $A_i$ is equal to $X_1^n$ and the conclusion holds. 

\par

\begin{figure}[htbp]
\begin{center}
\includegraphics[width=140mm]{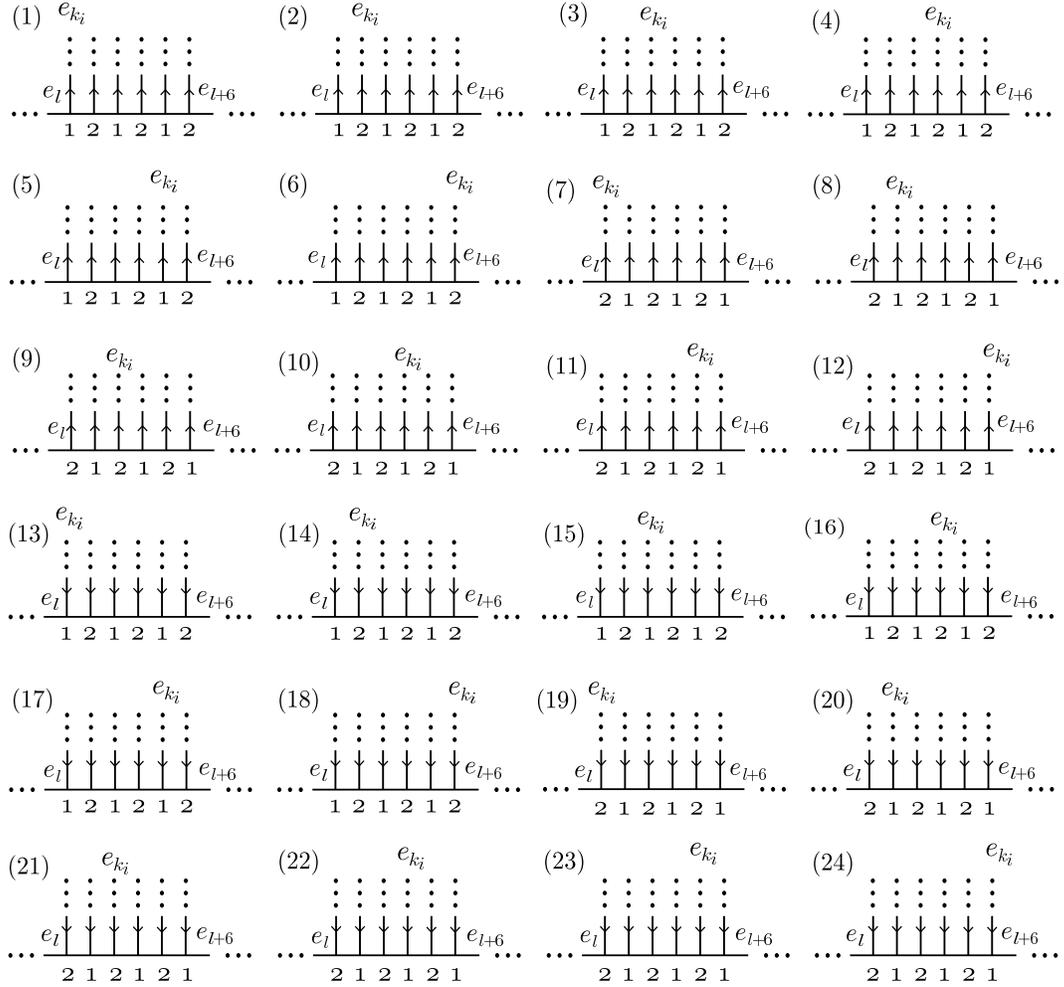}
\end{center}
\caption{24 cases about $e_{k_i}$ and the boundary comb containing $e_{k_i}$. }
\label{24cases}
\end{figure}

{\it Case.3}: Suppose that $p_i$ is an end point of $(1,6)$-edge and $e_{k_i}$ is contained in a boundary comb. 
Let $e_l$ and $e_{l+6}$ be two edges at the end of the boundary comb which contains $e_{k_i}$. 
Then one of 24 cases illustrated in Figure\ref{24cases} occurs. 

The intersection word of a path which starts from $p_0$, travels in $V$ across the edges $e_1,\ldots,e_{l-1}$, ends near the boundary comb is equal to $X_1^n$, where $n$ is an integer. 
Since the label of the $(1,6)$-edge incident to $p_i$ is $2$, $\rho_\Gamma(a_i)$ is calculated as follows: 
{\allowdisplaybreaks
\begin{align*}
\rho_\Gamma(a_i) & =\begin{cases}
X_1^{n+1}X_2X_1^{-n-1} & \text{(If the case (1) or (2) occurs.)}, \\
X_1^{n+1}X_2X_1X_2X_1^{-1}X_2^{-1}X_1^{-n-1} & \text{(If the case (3) or (4) occurs.)}, \\
X_1^{n+1}X_2X_1X_2X_1X_2X_1^{-1}X_2^{-1}X_1^{-1}X_2^{-1}X_1^{-n-1} & \text{(If the case (5) or (6) occurs.)}, \\
X_1^{n}X_2X_1^{-n} & \text{(If the case (7) occurs.)}, \\
X_1^{n}X_2X_1X_2X_1^{-1}X_2^{-1}X_1^{-n} & \text{(If the case (8) or (9) occurs.)}, \\
X_1^{n}X_2X_1X_2X_1X_2X_1^{-1}X_2^{-1}X_1^{-1}X_2^{-1}X_1^{-n} & \text{(If the case (10) or (11) occurs.)}, \\
X_1^{n}X_2X_1X_2X_1X_2X_1X_2X_1^{-1}X_2^{-1}X_1^{-1}X_2^{-1}X_1^{-1}X_2^{-1}X_1^{-n} & \text{(If the case (12) occurs.)}, \\
X_1^{n-1}X_2X_1^{-n+1} & \text{(If the case (13) or (14) occurs.)}, \\
X_1^{n-1}X_2^{-1}X_1^{-1}X_2X_1X_2X_1^{-n+1} & \text{(If the case (15) or (16) occurs.)}, \\
X_1^{n-1}X_2^{-1}X_1^{-1}X_2^{-1}X_1^{-1}X_2X_1X_2X_1X_2X_1^{-n+1} & \text{(If the case (17) or (18) occurs.)}, \\
X_1^{n}X_2X_1^{-n} & \text{(If the case (19) occurs.)}, \\
X_1^{n}X_2^{-1}X_1^{-1}X_2X_1X_2X_1^{-n} & \text{(If the case (20) or (21) occurs.)}, \\
X_1^{n}X_2^{-1}X_1^{-1}X_2^{-1}X_1^{-1}X_2X_1X_2X_1X_2X_1^{-n} & \text{(If the case (22) or (23) occurs.)}, \\
X_1^{n}X_2^{-1}X_1^{-1}X_2^{-1}X_1^{-1}X_2^{-1}X_1^{-1}X_2X_1X_2X_1X_2X_1X_2X_1^{-n} & \text{(If the case (24) occurs.)}. 
\end{cases}
\end{align*}
}
Since $X_1X_2X_1X_2^{-1}X_1^{-1}X_2^{-1}=(X_1X_2)^6=E$, we obtain: 
\begin{align*}
\rho_\Gamma(a_i) & =\begin{cases}
X_1^{n+1}X_2X_1^{-n-1} & \text{(If the case (1), (2), (22) or (23) occurs.)}, \\
X_1 & \text{(If the case (3), (4), (8), (9), (15), (16), (20) or (21) occurs.)}, \\ 
X_1^{n-1}X_2X_1^{-n+1} & \text{(If the case (5), (6), (10), (11), (13) or (14) occurs.)}, \\
X_1^{n}X_2X_1^{n} & \text{(If the case (7), (12), (17), (18), (19) or (24) occurs.)}. 
\end{cases}
\end{align*}
For each case, the conclusion holds. 

\begin{figure}[h!]
\begin{center}
\includegraphics[width=150mm]{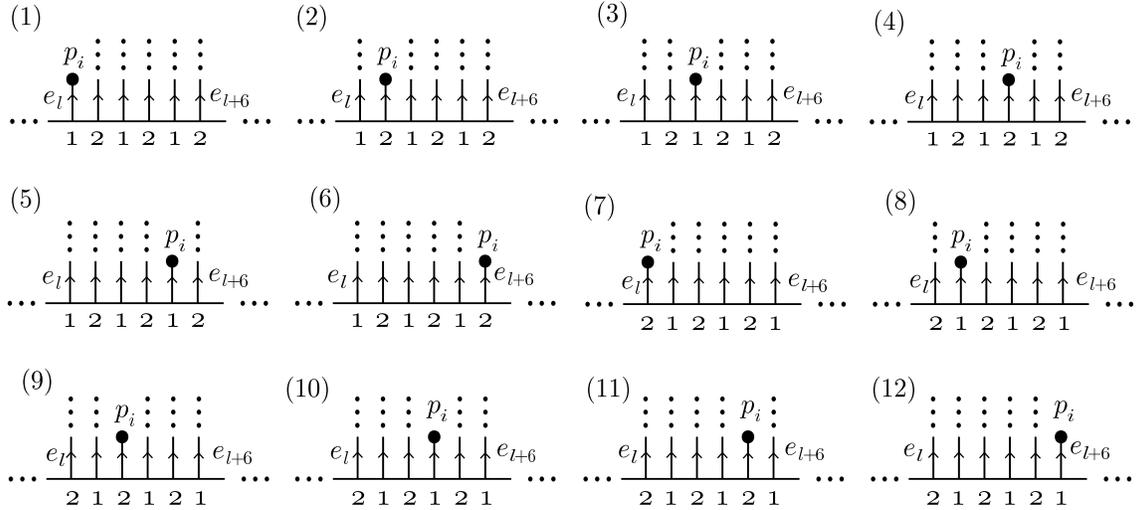}
\end{center}
\caption{12 cases about $e_{j}$ and the boundary comb containing $e_{j}$. }
\label{12cases}
\end{figure}

{\it Case.4}: Suppose that $p_i$ is an end point of $(\partial,1)$-edge $e_j$ and $e_j$ is contained in a boundary comb. 
Let $e_l$ and $e_{l+6}$ be two edges at the end of the boundary comb which contains $e_{j}$. 
Since the degree-$1$ vertex $p_i$ is negative, one of 12 cases illustrated in Figure\ref{12cases} occurs. 
We assume that the intersection word of a path which starts from $p_0$, travels in $V$ across the edges $e_1,\ldots,e_{l-1}$, ends near the boundary comb is equal to $X_1^n$, where $n$ is an integer. 
Then $\rho_\Gamma(a_i)$ is calculated as follows: 

\begin{align*}
\rho_\Gamma(a_i) & =\begin{cases}
X_1 & \text{(If the case (1) occurs.)}, \\
X_1^{n+1}X_2X_1^{-n-1} & \text{(If the case (2) occurs.)}, \\
X_1^{n+1}X_2X_1X_2^{-1}X_1^{-n-1} & \text{(If the case (3) occurs.)}, \\
X_1^{n+1}X_2X_1X_2X_1^{-1}X_2^{-1}X_1^{-n-1} & \text{(If the case (4) occurs.)}, \\
X_1^{n+1}X_2X_1X_2X_1X_2^{-1}X_1^{-1}X_2^{-1}X_1^{-n-1} & \text{(If the case (5) occurs.)}, \\
X_1^{n+1}X_2X_1X_2X_1X_2X_1^{-1}X_2^{-1}X_1^{-1}X_2^{-1}X_1^{-n-1} & \text{(If the case (6) occurs.)}, \\
X_1^{n}X_2X_1^{-n} & \text{(If the case (7) occurs.)}, \\
X_1^{n}X_2X_1X_2^{-1}X_1^{-n} & \text{(If the case (8) occurs.)}, \\
X_1^{n}X_2X_1X_2X_1^{-1}X_2^{-1}X_1^{-n} & \text{(If the case (9) occurs.)}, \\
X_1^{n}X_2X_1X_2X_1X_2^{-1}X_1^{-1}X_2^{-1}X_1^{-n} & \text{(If the case (10) occurs.)}, \\
X_1^{n}X_2X_1X_2X_1X_2X_1^{-1}X_2^{-1}X_1^{-1}X_2^{-1}X_1^{-n} & \text{(If the case (11) occurs.)}, \\
X_1^{n}X_2X_1X_2X_1X_2X_1X_2^{-1}X_1^{-1}X_2^{-1}X_1^{-1}X_2^{-1}X_1^{-n} & \text{(If the case (12) occurs.)}. \\
\end{cases}
\end{align*}

By the relations of $SL(2,\mathbb{Z})$, we obtain: 

\begin{align*}
\rho_\Gamma(a_i) & =\begin{cases}
X_1 & \text{(If the case (1), (4), (9) or (12) occurs.)}, \\
X_1^{n+1}X_2X_1^{-n-1} & \text{(If the case (2) or (5) occurs.)}, \\
X_1^nX_2X_1^{-n} & \text{(If the case (3), (6), (7) or (10) occurs.)}, \\
X_1^{n-1}X_2X_1^{-n+1} & \text{(If the case (8) or (11) occurs.)}. \\
\end{cases}
\end{align*}

For each cases, the conclusion holds. 

Combining the conclusions we obtain in Cases.1, 2, 3 and 4, we complete the proof of Theorem \ref{main in chart sec}.  \hfill $\square$

\section{Examples}\label{examples}
In this section, we give examples of genus-$1$ SBLFs. 
We have already known that the monodromy representations of the higher sides of such fibrations satisfy a certain condition by Theorem \ref{main in chart sec}. 
We first give some examples of sequences of elements of $SL(2,\mathbb{Z})$ and examine what $4$-manifolds these sequences represent. 
We denote by $T(n_1,\ldots,n_s)$ the following sequence: 
\[
(X_1^{-n_1}X_2X_1^{n_1},\ldots,X_1^{-n_s}X_2X_1^{n_s}). 
\]

\begin{prop}\label{ex_seq}

We define sequences $S_r$ and $T_s$ as follows: 
\begin{align*}
S_r & =(X_1,\ldots,X_1) \hspace{1em}\text{($r$ $X_1$'s stand in a line.)}. \\
T_s & =\begin{cases}
(X_1^{-1}X_2X_1,X_1X_2X_1^{-1}) & (s=2), \\
(X_1^{-3}X_2X_1^{3},X_2,X_1^{3}X_2X_1^{-3}) & (s=3), \\
(X_1^{-2s+3}X_2X_1^{2s-3},X_1^{-2s+6}X_2X_1^{2s-6},X_1^{-2s+10}X_2X_1^{2s-10},X_1^{-2s+14}X_2X_1^{2s-14},\ldots \\
 \hspace{2em}\ldots ,X_1^{2s-10}X_2X_1^{-2s+10},X_1^{2s-6}X_2X_1^{-2s+6},X_1^{2s-3}X_2X_1^{-2s+3}) & (s\geq 4). 
\end{cases}
\end{align*}
that is, $T_s=T(n_1,\ldots,n_s)$, where $n_1=2s-3$, $n_s=-2s+3$, $n_i=2s-6+4(i-1)$ ($i=2,\ldots,s-1$). 
Then, $w(S_r)={X_1}^r$ and $w(T_s)=(-1)^{s+1}X_1^{-5s+6}$. 
So these sequences are Hurwitz systems of some genus-$1$ SBLF. 

\end{prop}

{\it Proof}: It is obvious that $w(S_r)$ is equal to ${X_1}^r$. 
We prove $w(T_s)=(-1)^{s+1}X_1^{-5s+6}$ by induction on $s$. 
Since $(X_1X_2)^3=-E$ and $X_1X_2X_1=X_2X_1X_2$, we obtain: 
{\allowdisplaybreaks
\begin{align*}
X_2{X_1}^2X_2{X_1}^2 & =X_2X_1(X_2X_1X_2)X_1 \\
& =-E. 
\end{align*}
}
So $w(T_2)$ and $w(T_3)$ is computed as follows: 
{\allowdisplaybreaks
\begin{align*}
w(T_2) & =(X_1^{-1}X_2X_1)(X_1X_2X_1^{-1}) \\
& =X_1^{-1}(X_2{X_1}^2X_2)X_1^{-1} \\
& =X_1^{-1}(-X_1^{-2})X_1^{-1} \\
& =-X_1^{-4}. \\
w(T_3) & =(X_1^{-3}X_2{X_1}^3)X_2({X_1}^3X_2X_1^{-3}) \\
& =X_1^{-3}X_2X_1({X_1}^2X_2{X_1}^2)X_1X_2X_1^{-3} \\
& =X_1^{-3}X_2X_1(-X_2^{-1})X_1X_2X_1^{-3} \\
& =-X_1^{-3}(X_2X_1X_2^{-1})X_1X_2X_1^{-3} \\
& =-X_1^{-3}(X_1^{-1}X_2X_1)X_1X_2X_1^{-3} \\
& =-X_1^{-4}(X_2{X_1}^2X_2)X_1^{-3} \\
& =X_1^{-9}. 
\end{align*}
}

By the definition of $T_s$, $w(T_s)$ is represented by $w(T_{s-2})$ as follows: 

{\allowdisplaybreaks
\begin{align*}
w(T_s) & =(X_1^{-2s+3}X_2X_1^{2s-3})(X_1^{-2s+6}X_2X_1^{2s-6})(X_1^{-2s+7}X_2^{-1}X_1^{2s-7})w(T_{s-2}) \\
& \hspace{1.4em} (X_1^{2s-7}X_2^{-1}X_1^{-2s+7})(X_1^{2s-6}X_2X_1^{-2s+6})(X_1^{2s-3}X_2X_1^{-2s+3}) \\
& =X_1^{-2s+3}X_2{X_1}^2(X_1X_2X_1)X_2^{-1}X_1^{2s-7}w(T_{s-2})X_1^{2s-7}X_2^{-1}(X_1X_2X_1){X_1}^2X_2X_1^{-2s+3} \\
& =X_1^{-2s+3}(X_2{X_1}^2X_2)X_1^{2s-6}w(T_{s-2})X_1^{2s-6}(X_2{X_1}^2X_2)X_1^{-2s+3} \\
& =X_1^{-2s+3}(-X_1^{-2})X_1^{2s-6}w(T_{s-2})X_1^{2s-6}(-X_1^{-2})X_1^{-2s+3} \\
& =X_1^{-5}w(T_{s-2})X_1^{-5}. 
\end{align*}
}
\noindent
Thus the conclusion holds by the induction hypothesis. 
This completes the proof of Proposition \ref{ex_seq}.  \hfill $\square$

\begin{thm}\label{aboutS_r}
Let $f:M\rightarrow S^2$ be a genus-$1$ SBLF. 
Suppose that $W_f=S_r$. 
Then $M$ is diffeomorphic to one of the following $4$-manifolds: 
\begin{enumerate}[(1)]

\item $\sharp r\overline{\mathbb{CP}^2}$; 

\item $L\sharp r\overline{\mathbb{CP}^2}$; 

\item $S^1\times S^3\sharp S \sharp r\overline{\mathbb{CP}^2}$, 

\end{enumerate}
\noindent
where $S$ is either of the manifolds $S^2\times S^2$ and $S^2\tilde{\times}S^2$ 
and $L$ is either of the manifolds $L_n$ and $L_n^\prime$. 

\end{thm}

Before proving Theorem \ref{aboutS_r}, we review the definition and some properties of $L_n$ and $L_n^\prime$. 
For more details, see \cite{Pao}. 
Let $N_0$ and $N_1$ be $4$-manifolds diffeomorphic to $D^2\times T^2$. 
The boundaries of $N_0$ and $N_1$ are $\partial D^2\times T^2$. 
Let $(t,x,y)$ be a coordinate of $\mathbb{R}^3$. 
We identify $\partial D^2\times T^2$ with $\mathbb{R}^3/\mathbb{Z}^3$. 
The group $GL(3,\mathbb{Z})$ naturally acts on $\mathbb{R}^3$ and this action descends to an action on the lattice $\mathbb{Z}^3$. 
So $GL(3,\mathbb{Z})$ acts on $\partial D^2\times T^2$. 
For an element $A$ of $GL(3,\mathbb{Z})$, we denote by $f_A$ a self-diffeomorphism of $\partial D^2\times T^2$ defined as follows: 
\[
f_A([t,x,y])=[(t,x,y){}^tA]. 
\]
We define elements $A_n$ and $A_n^\prime$ of $GL(3,\mathbb{Z})$ as follows: 
\[
A_n=
\begin{pmatrix}
0 & 1 & 1 \\
0 & n & n-1 \\
1 & n & 0 
\end{pmatrix}
\hspace{0.2em},\hspace{0.2em}
A_n^\prime=
\begin{pmatrix}
0 & 1 & 1 \\
0 & n & n-1 \\
1 & n-1 & 0 
\end{pmatrix}.
\]
\noindent
Let $B^3\subset D^2\times S_y^1$ be an embedded $3$-ball, 
where $S_x^1$ and $S_y^1$ represent circles with coordinates $x$ and $y$, respectively. 
Then, 
\[
B^3\times {S_x}^1\subset(D^2\times {S_y}^1)\times {S_x}^1=D^2\times T^2=N_0. 
\]
We take a diffeomorphism $\iota:S^2\rightarrow\partial B^3$. 
We obtain a diffeomorphism as follows: 
\[
h=\iota\times id:S^2\times S^1\rightarrow \partial B^3\times S^1. 
\]
We define $L_n$ and $L_n^\prime$ as follows: 
\begin{align*}
L_n & =D^2\times S^2\cup_h(N_0\setminus(\text{int}B^3\times S^1))\cup_{f_{A_n}}N_1, \\
L_n^\prime & =D^2\times S^2\cup_h(N_0\setminus(\text{int}B^3\times S^1))\cup_{f_{A_n^\prime}}N_1. 
\end{align*}

\begin{rem}

The original definitions of $L_n$ and $L_n^\prime$ are different from the above definitions. 
However, both two definitions are equivalent (c.f. Lemma V.7 in \cite{Pao}). 
We also remark that these manifolds were constructed in Example 1 of section 8.2 of \cite{ADK}, 
although they did not state that their examples were actually the manifolds $L_n$ and $L_n^\prime$. 
Indeed, $N_1$ (resp. $N_0\setminus(\text{int}B^3\times S^1)$, $D^2\times S^2$) in our paper corresponds to $X_-$ (resp. $W$, $X_+$) in \cite{ADK}. 
Moreover, in the construction in \cite{ADK}, we obtain $L_n$ (resp. $L_n^\prime$) 
if we glue $X_-$ by using the element $(k,l)\in\mathbb{Z}^2\cong \pi_1(\text{Diff}(T^2))$ which satisfies gcd$(k,l)=n$ and $X_+$ by using the trivial (resp. non-trivial) element of $\pi_1(\text{Diff}(S^2))\cong \mathbb{Z}/2$. 

\end{rem}

\begin{figure}[htbp]
\begin{center}
\includegraphics[width=60mm]{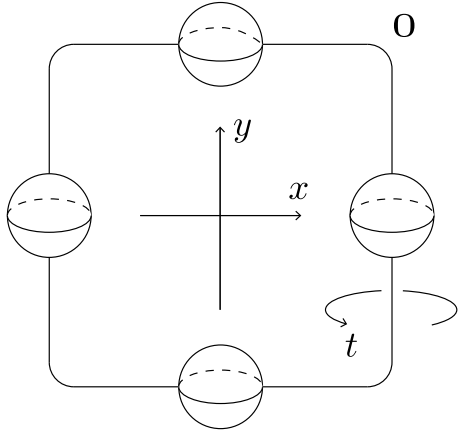}
\end{center}
\caption{A Kirby diagram of $N_1$. 
$t$ represents the coordinate of $\partial D^2$ and $x$ and $y$ represent the coordinate of $T^2$. }
\label{N_1}
\end{figure}

We next take handle decompositions of $L_n$ and $L_n^\prime$. 
Since $N_1=D^2\times T^2$, a Kirby diagram of $N_1$ is as shown in Figure \ref{N_1}. 
The coordinate $(t,x,y)$ is also described as in Figure \ref{N_1}. 

\begin{figure}[htbp]
\begin{center}
\includegraphics[width=65mm]{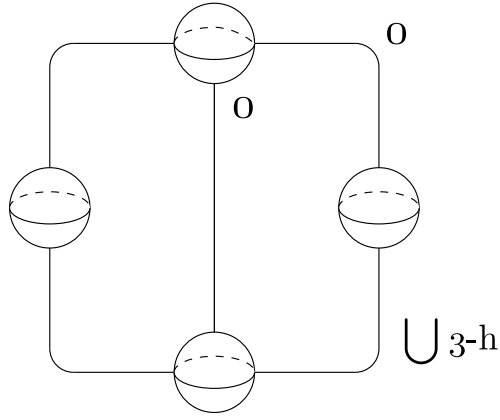}
\end{center}
\caption{A diagram of $(N_0\setminus(\text{int}B^3\times S^1))\cup_{f_{A_n}}N_1$ and $(N_0\setminus(\text{int}B^3\times S^1))\cup_{f_{A_n^\prime}}N_1$. }
\label{N_0N_1}
\end{figure}

$B^3\times S^1$ has a handle decomposition consisting of a $0$-handle and a $1$-handle. 
So we can decompose $N_0\setminus(\text{int}B^3\times S^1)$ as follows: 
\[
N_0\setminus(\text{int}B^3\times S^1)=\partial N_0\times I\cup(2\text{-handle})\cup(3\text{-handle}). 
\]
Let $C_1\subset\partial N_0$ be an attaching circle of the $2$-handle. 
By the construction of the decomposition, we obtain: 
\[
C_1=\{[t,0,0]\in\partial N_0|t\in[0,1]\}. 
\]
Since $f_{A_n}([t,0,0])=f_{A_n^\prime}([t,0,0])=[0,0,t]$, an attaching circle of the $2$-handle is in a regular fiber and along $y$-axis in the diagram of $N_1$. 
Since $f_{A_n}([t,0,\varepsilon])=f_{A_n^\prime}([t,0,\varepsilon])=[\varepsilon,(n-1)\varepsilon,t]$ for sufficiently small $\varepsilon>0$, 
the framing of the $2$-handle is along a regular fiber. 
Hence the diagram of $(N_0\setminus(\text{int}B^3\times S^1))\cup_{f_{A_n}}N_1$ and $(N_0\setminus(\text{int}B^3\times S^1))\cup_{f_{A_n^\prime}}N_1$ is as shown in Figure \ref{N_0N_1}. 

We can decompose $D^2\times S^2$ as follows: 
\[
D^2\times S^2=\partial D^2\times S^2\times I\cup(2\text{-handle})\cup(4\text{-handle}). 
\]

Let $C_2\subset\partial N_0$ be the image under $h$ of the attaching circle of the $2$-handle of $D^2\times S^2$. 
After moving $C_2$ by isotopy in $N_0$, we obtain: 
\[
C_2=\{[0,t,\delta]\in\partial N_0|t\in[0,1]\}, 
\]
where $\delta>0$ is sufficiently small.

\begin{figure}[htbp]
\begin{center}
\includegraphics[width=135mm]{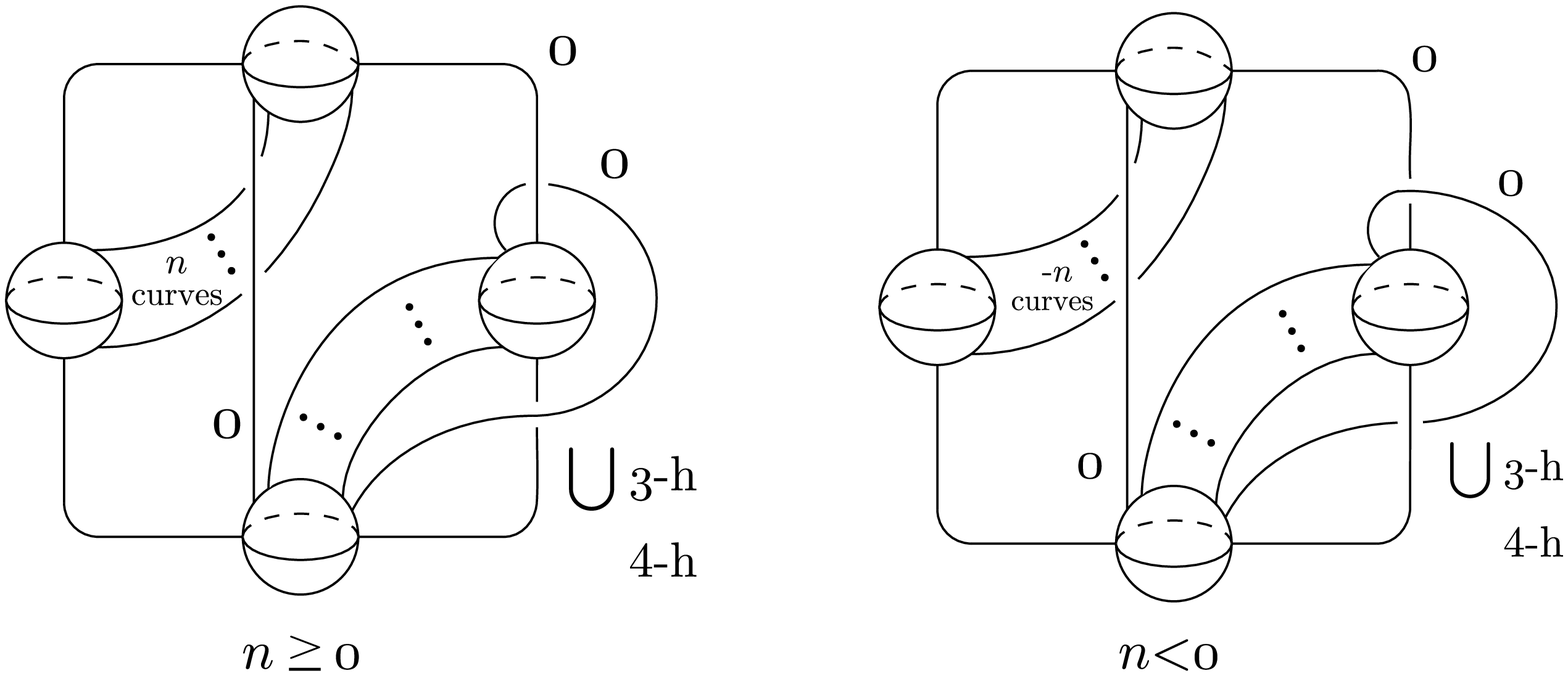}
\end{center}
\caption{Left: A Kirby diagram of $L_n$ for $n\geq 0$. Right: A Kirby diagram of $L_n$ for $n<0$. }
\label{L_n}
\end{figure}

The framing of the $2$-handle is $\{[0,t,\delta^\prime]\in\partial N_0|t\in[0,1]\}$, where $\delta^\prime>\delta$ is sufficiently small. 
Since $f_{A_n}([0,t,\delta])=[t+\delta,nt+(n-1)\delta,nt]$, we can describe the attaching circle of the $2$-handle of $D^2\times S^2$ contained in $L_n$ and the knot representing the framing of the $2$-handle in the diagram described in Figure \ref{N_0N_1}. 
Eventually, we can describe a diagram of $L_n$ as shown in Figure \ref{L_n}. 

\begin{figure}[htbp]
\begin{center}
\includegraphics[width=135mm]{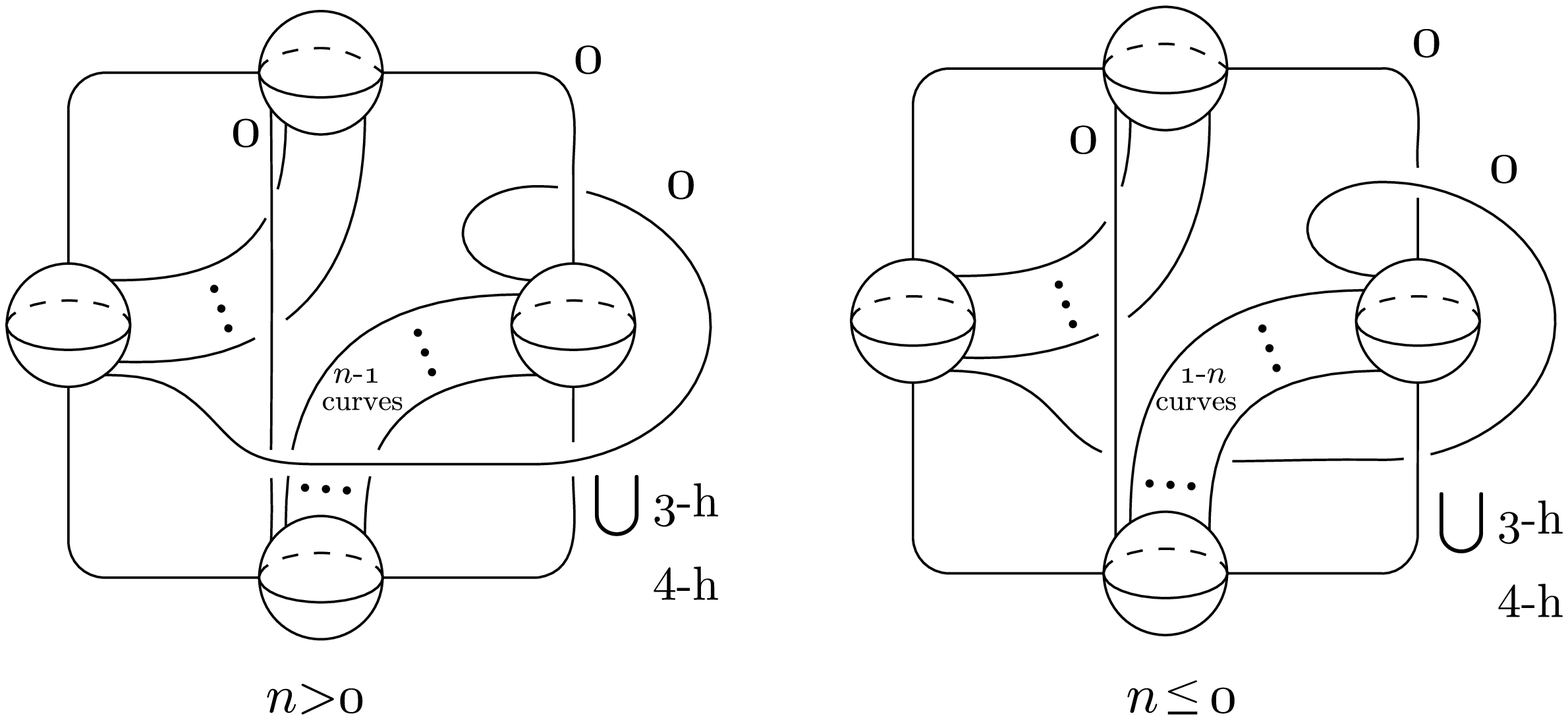}
\end{center}
\caption{Left: A Kirby diagram of $L_n^\prime$ for $n\geq 0$. Right: A Kirby diagram of $L_n^\prime$ for $n<0$. }
\label{L_n2}
\end{figure}

Similarly, a diagram of $L_n^\prime$ is described as shown in Figure \ref{L_n2}. 

By the diagrams of $L_n$ and $L_n^\prime$ described in Figure \ref{L_n} and Figure \ref{L_n2}, both $L_n$ and $L_n^\prime$ admit genus-$1$ SBLFs without Lefschetz singularities. 
We can easily prove by Kirby calculus that $L_{-n}$ (resp. $L_{-n}^\prime$) is diffeomorphic to $L_n$ (resp. $L_n^\prime$). 

\par

\begin{figure}[htbp]
\begin{center}
\includegraphics[width=130mm]{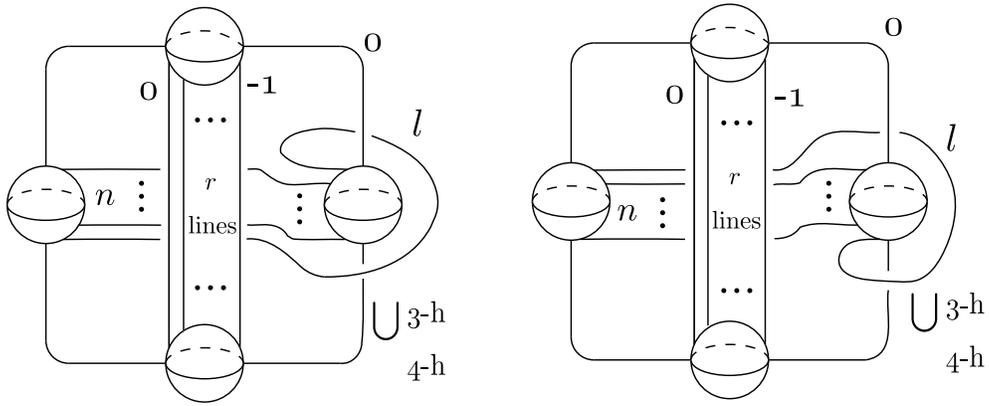}
\end{center}
\caption{A Kirby diagram of a genus-$1$ SBLF whose Hurwitz system is $S_r$. 
Framings of $r$ $2$-handles parallel to the $2$-handle of the round $2$-handle are all $-1$. }
\label{BLF-S_r}
\end{figure}

{\it Proof of Theorem \ref{aboutS_r}}: The higher side of $f$ is obtained by attaching $r$ $2$-handles to a trivial $T^2$ bundle over $D^2$. 
Each attaching circle of the $2$-handle is in a regular fiber and isotopic to a simple closed curve $\gamma_{1,0}$. 
Since $w(W_f)={X_1}^r$, a $2$-handle of a round $2$-handle is attached along $\gamma_{1,0}$ in a regular fiber of the boundary of the higher side. 
We obtain $M$ by attaching a $2$-handle and a $4$-handle to the $4$-manifold obtained by successive handle attachment to $D^2\times T^2$. 
If the attaching circle of the $2$-handle of $D^2\times S^2$ goes through the $1$-handle that the $2$-handle of the round handle goes through, 
we can slide the $2$-handle of $D^2\times S^2$ to the $2$-handle of the round handle so that the $2$-handle of $D^2\times S^2$ does not go through the $1$-handle. 
Thus a Kirby diagram of $M$ is one of the diagrams in Figure \ref{BLF-S_r}, where $n$ and $l$ are integers. 
It is obvious that both two $4$-manifolds illustrated in Figure \ref{BLF-S_r} are diffeomorphic to each other. 
We denote by $M_{n,l}$ the $4$-manifold illustrated in Figure \ref{BLF-S_r}. 
We remark that $l$ framed knot in Figure \ref{BLF-S_r} represents a $2$-handle of $D^2\times S^2$ and the attachment of the lower side depends only on the parity of $l$. 
So $M_{n,l}$ and $M_{n,l^\prime}$ are diffeomorphic to each other if $l\equiv l^\prime$ (mod $2$).

We will determine what $4$-manifold $M_{n,l}$ is by Kirby calculus. 

\par

\begin{figure}[htb]
\begin{center}
\includegraphics[width=135mm]{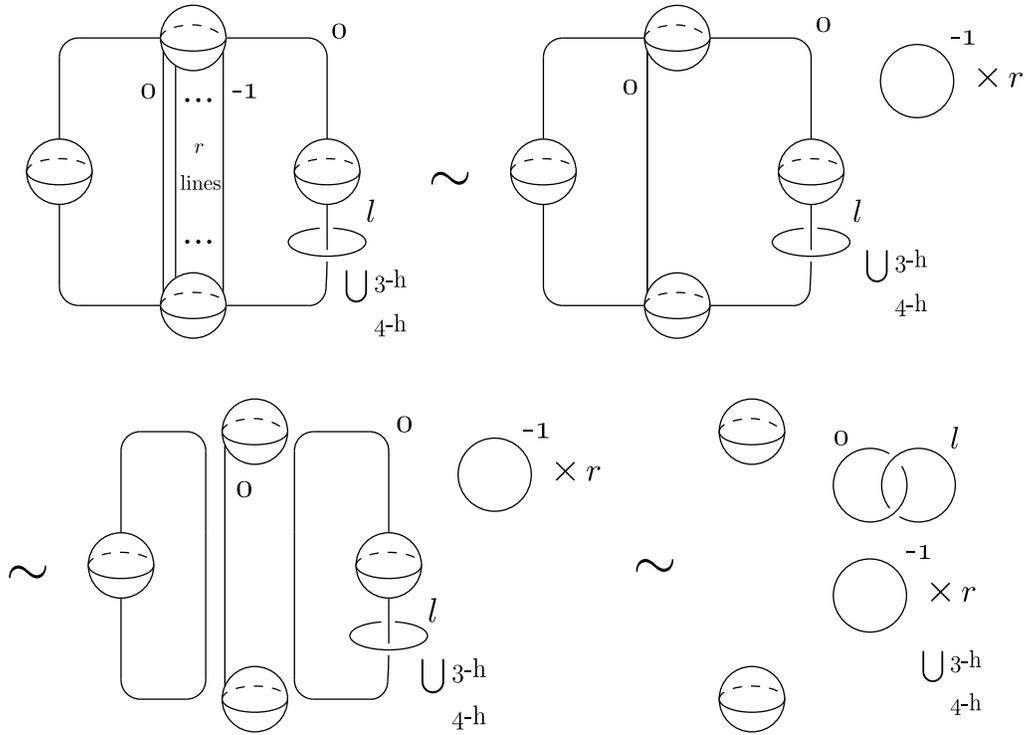}
\end{center}
\caption{A Kirby diagram of $M_{0,l}$. }
\label{S_r-n=0}
\end{figure}

We first examine the case $n=0$. 
A Kirby diagram of $M_{0,l}$ is as shown in Figure \ref{S_r-n=0}. 
We first slide $r$ $2$-handles representing Lefschetz singularities to the $2$-handle of the round $2$-handle. 
Then we slide the $2$-handle of $D^2\times T^2$ to the $2$-handle of the round $2$-handle 
and move this $2$-handle so that the attaching circle of the $2$-handle does not go through $1$-handles. 
We obtain the last diagram of Figure \ref{S_r-n=0} by eliminating the obvious canceling pair. 

Thus we obtain: 
\begin{align*}
M_{0,l}=S^1\times S^3\sharp S\sharp r\overline{\mathbb{CP}^2}, 
\end{align*}
where $S$ is equal to $S^2\times S^2$ if $l$ is even and $S^2\tilde{\times} S^2$ if $l$ is odd. 

\par

We next examine the case $n=1$. 
We first slide $r$ $(-1)$-framed $2$-handles and the $2$-handle of $D^2\times T^2$ as we slide in the case $n=0$. 
Then we can eliminate the obvious canceling pair. 
We next move the $2$-handle of $D^2\times T^2$ so that the attaching circle of the $2$-handle does not go through the $1$-handle. 
Then we can eliminate two pairs of handles. 
Thus we can prove: 
\[
M_{1,l}=\sharp r\overline{\mathbb{CP}^2}. 
\]

\par

Lastly, we examine the case $n\geq 2$. 
A Kirby diagram of $M_{n,l}$ is illustrated in Figure \ref{S_r-other} and 
we change the diagram as shown in Figure \ref{S_r-other} by performing similar calculus to the calculus illustrated in Figure \ref{S_r-n=0}. 

\begin{figure}[htbp]
\begin{center}
\includegraphics[width=135mm]{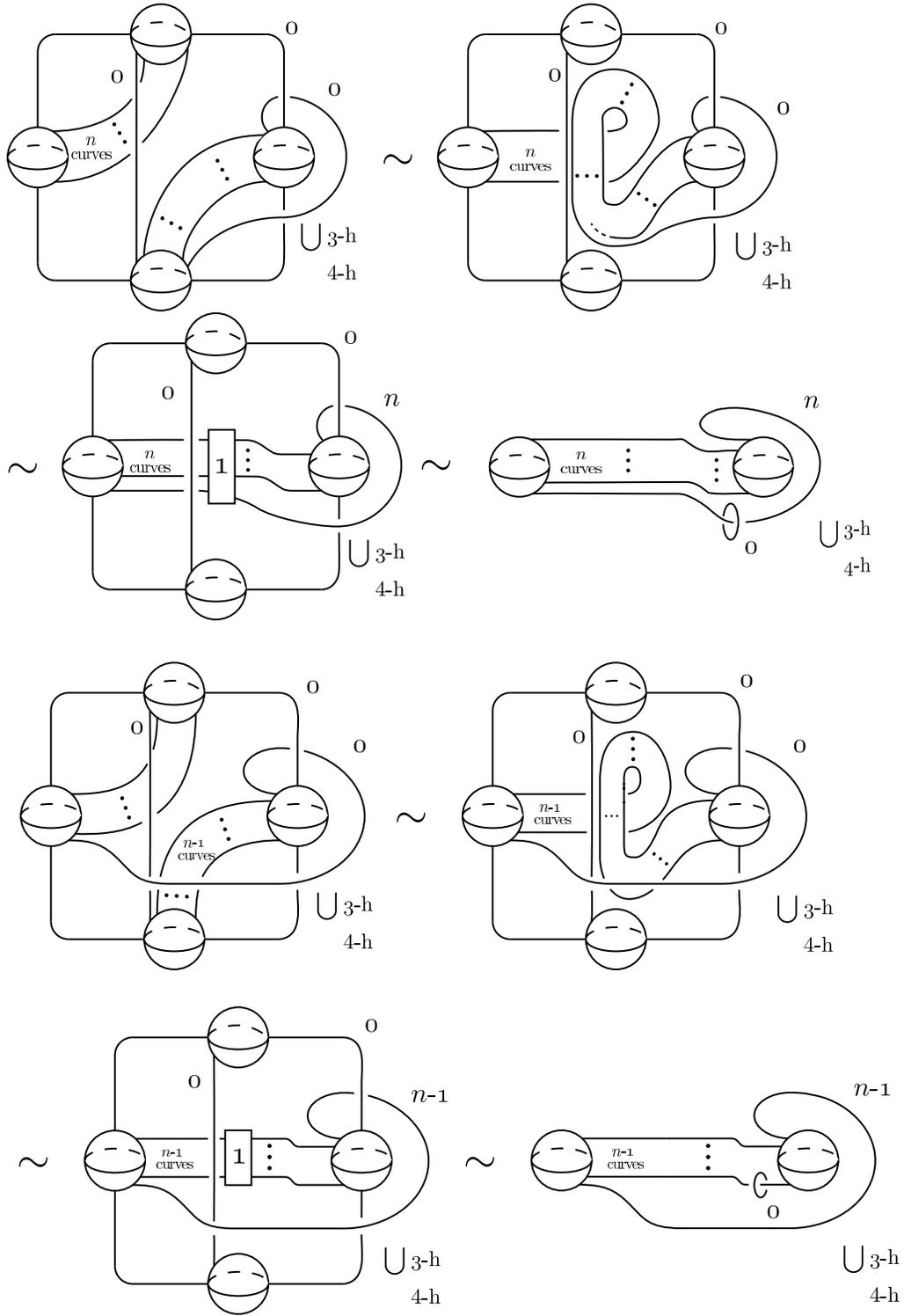}
\end{center}
\caption{The upper four diagrams describe $L_n$, while the lower four diagrams describe $L_n^\prime$. }
\label{movesL_n}
\end{figure}

\begin{figure}[t!]
\begin{center}
\includegraphics[width=145mm]{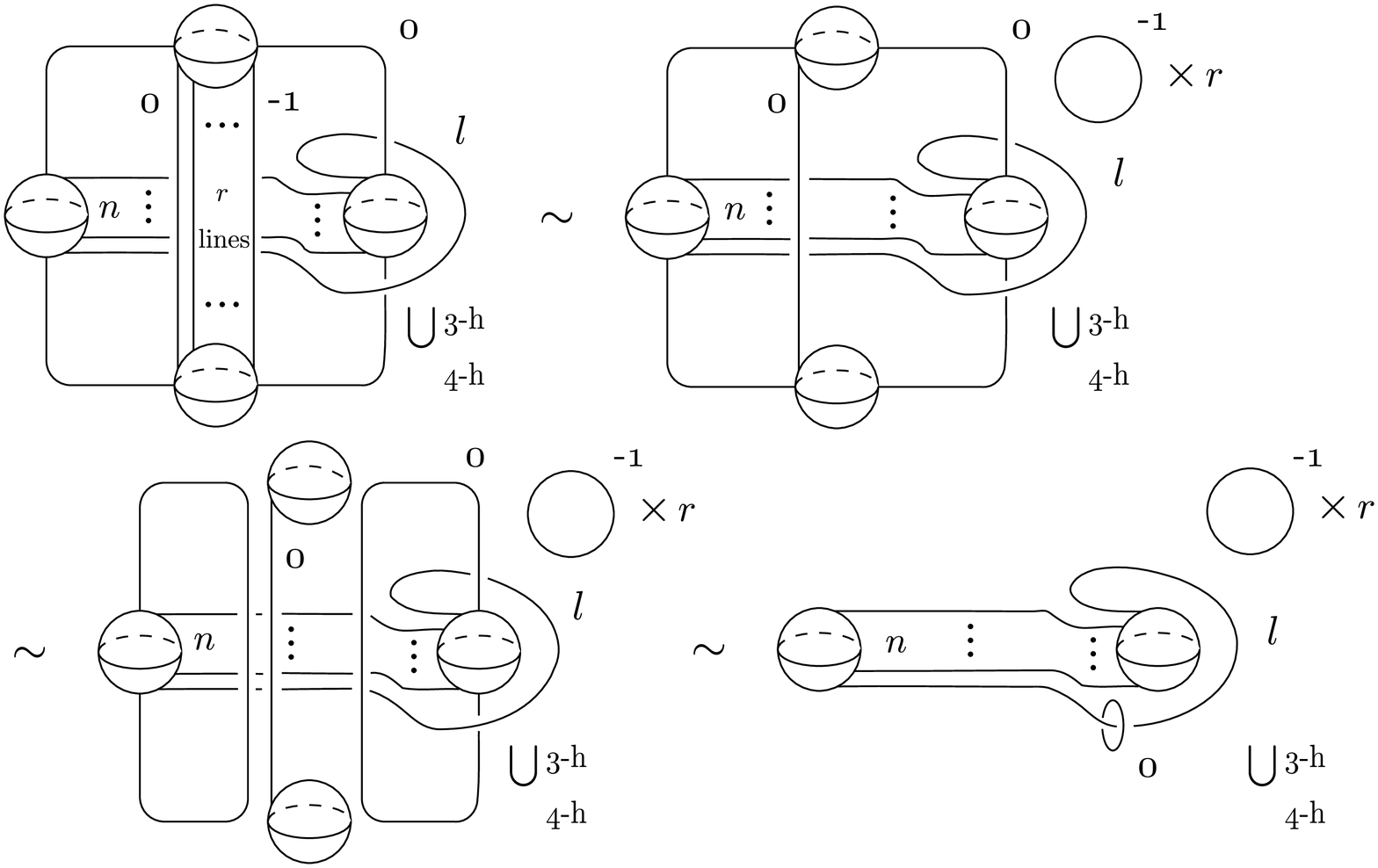}
\end{center}
\caption{A Kirby diagram of $M_{n,l}$, where $n\geq 2$. }
\label{S_r-other}
\end{figure}

We can change the diagrams of $L_n$ and $L_n^\prime$ illustrated in Figure \ref{L_n} and \ref{L_n2} as shown in Figure \ref{movesL_n}. 
The upper four diagrams in Figure \ref{movesL_n} describe $L_n$, where $n\geq 0$. 
We first slide $2$-handle of $D^2\times S^2$ to the $2$-handle of $N_0$. 
Then we move the $2$-handle of $D^2\times S^2$ by isotopy and slide the $2$-handle of $D^2\times T^2$ to the $2$-handle of $N_0$. 
To obtain the last diagram of the four, we first untwist the $2$-handle of $D^2\times S^2$ by using $0$-framed meridian of this. 
The framing of this handle is unchanged in these moves. 
Then we eliminate a canceling pair. 
The lower four diagrams in Figure \ref{movesL_n} describe $L_n^\prime$, where $n\geq 0$. 
By the same method as above, we obtain the last diagram of the four. 

Thus we obtain: 
\[
M_{n,l}=L\sharp r\overline{\mathbb{CP}^2}, 
\]
where $L$ is either $L_n$ or $L_n^\prime$. 
This completes the proof of Theorem \ref{aboutS_r}.  \hfill $\square$

\begin{thm}\label{aboutT_s}

Let $f:M\rightarrow S^2$ be a genus-$1$ SBLF. 
Suppose that $W_f=T_s$. 
Then $M$ is diffeomorphic to $S\sharp (s-2)\mathbb{CP}^2$, where $S$ is either of the manifolds $S^2\times S^2$ and $S^2\tilde{\times}S^2$. 

\end{thm}

Before proving Theorem \ref{aboutT_s}, we change $T_s$ by applying elementary transformations and simultaneous conjugations 
so that we can easily determine what $4$-manifold $T_s$ represents. 

\begin{lem}\label{transT_s}

Suppose that $s\geq 3$. 
Then, by successive application of elementary transformations and simultaneous conjugations, we can change $T_s$ into a following sequence: 
\[
(T_{2,1},T_{2,3},\ldots,T_{2,2s-5},T_{1,s-1},T_{1,-1}). 
\]

\end{lem}

{\it Proof}: By applying simultaneous conjugation, we can transform $T_s$ as follows: 
{\allowdisplaybreaks
\begin{align*}
T_s & \rightarrow\begin{cases}
(X_1^{-5}X_2X_1^{5},X_1^{-2}X_2X_1^2,X_1X_2X_1^{-1}) & (s=3) \\
(X_1^{-4s+7}X_2X_1^{4s-7},X_1^{-4s+10}X_2X_1^{4s-10},X_1^{-4s+14}X_2X_1^{4s-14},X_1^{-4s+18}X_2X_1^{4s-18},\ldots \\
\hspace{.5em}\ldots ,X_1^{-6}X_2X_1^6,X_1^{-2}X_2X_1^2,X_1X_2X_1^{-1}) & (s\geq 4) 
\end{cases} \\
& = \hspace{.2em}\begin{cases}
(T_{5,1},T_{2,1},T_{-1,1}) & (s=3), \\
(T_{4s-7,1},T_{4s-10,1},T_{4s-14,1},T_{4s-18,1},\ldots,T_{6,1},T_{2,1},T_{-1,1}) & (s\geq 4). 
\end{cases}
\end{align*}
}

We first prove that the following formula holds up to orientation of simple closed curves: 
\begin{equation}\label{scc1}
T_{2,2l-1}\circ T_{2,2l-3}\circ\cdots\circ T_{2,1}(\gamma_{4k-10,1})=\gamma_{4k-10-4l,4lk-4l^2-10l+1}, 
\end{equation}
where $k\geq 4$ and $l\leq k-3$. 

\par

We prove the formula (\ref{scc1}) by induction on $l$. 
By the Picard-Lefschetz formula, we obtain: 
{\allowdisplaybreaks
\begin{align*}
T_{2,1}(\gamma_{4k-10,1}) &= \gamma_{4k-10,1}-((4k-10)-2\cdot 1)\gamma_{2,1} \\
&= \gamma_{4k-10,1}-(4k-12)\gamma_{2,1} \\
&= \gamma_{-4k+14,-4k+13} 
\end{align*}
}
Thus the formula (\ref{scc1}) holds for $l=1$. 

\par

For general $l$, by the induction hypothesis, we obtain: 
{\allowdisplaybreaks
\begin{align*}
& T_{2,2l-1}\circ T_{2,2l-3}\circ\cdots\circ T_{2,1}(\gamma_{4k-10,1}) \\
= & T_{2,2l-1}(\gamma_{4k-6-4l,4lk-4k-4l^2-2l+7}) \\
= & \gamma_{4k-6-4l,4lk-4k-4l^2-2l+7}-((2l-1)(4k-6-4l)-2(4lk-4k-4l^2-2l+7))\gamma_{2,2l-1} \\
= & \gamma_{4k-6-4l,4lk-4k-4l^2-2l+7}-(4k-4l-8)\gamma_{2,2l-1} \\
= & \gamma_{-4k+10+4l,-4lk+4l^2+10l-1}. 
\end{align*}
}
So the formula (\ref{scc1}) holds for general $l$ and we obtain: 
\[
\text{conj}(T_{2,2k-7}\circ T_{2,2k-9}\circ\cdots\circ T_{2,1})(T_{4k-10,1})=T_{2,2k-5}. 
\]
Thus, for $s\geq 4$, we can transform $T_s$ as follows: 
{\allowdisplaybreaks
\begin{align*}
T_s & \rightarrow (T_{4s-7,1},T_{4s-10,1},T_{4s-14,1},T_{4s-18,1},\ldots,T_{6,1},T_{2,1},T_{-1,1}) \\
& \rightarrow (T_{4s-7,1},T_{4s-10,1},T_{4s-14,1},T_{4s-18,1},\ldots,T_{10,1},T_{2,1},T_{2,3},T_{-1,1}) \\
& \rightarrow (T_{4s-7,1},T_{4s-10,1},T_{4s-14,1},T_{4s-18,1},\ldots,T_{14,1},T_{2,1},T_{2,3},T_{2,5},T_{-1,1}) \\
& \rightarrow \cdots \\
& \rightarrow (T_{4s-7,1},T_{2,1},T_{2,3},\ldots,T_{2,2s-5},T_{-1,1}). 
\end{align*}
}

We next prove that the following formula holds up to orientation: 
\begin{equation}\label{scc2}
T_{2,2l-1}\circ T_{2,2l-3}\circ\cdots T_{2,1}(\gamma_{4k-7,1})=\gamma_{4k-7-4l,4lk-4l^2-7l+1}, 
\end{equation}
where $k\geq 3$ and $1\leq l \leq k-2$. 
We prove the formula (\ref{scc2}) by induction on $l$. 

\par

By the Picard-Lefschetz formula, we obtain: 
{\allowdisplaybreaks
\begin{align*}
T_{2,1}(\gamma_{4k-7,1}) & =\gamma_{4k-7,1}-((4k-7)-2)\gamma_{2,1} \\
& =\gamma_{4k-7,1}-(4k-9)\gamma_{2,1} \\
& =\gamma_{-4k+11,-4k+10}. 
\end{align*}
}
So the formula (\ref{scc2}) holds for $l=1$. 

\par

For general $l$, by the induction hypothesis, we obtain: 
{\allowdisplaybreaks
\begin{align*}
& T_{2,2l-1}\circ T_{2,2l-3}\circ\cdots T_{2,1}(\gamma_{4k-7,1}) \\
= & T_{2,2l-1}(\gamma_{4k-3-4l,4lk-4k-4l^2+l+4}) \\
= & \gamma_{4k-3-4l,4lk-4k-4l^2+l+4}-((2l-1)(4k-3-4l)-2(4lk-4k-4l^2+l+4))\gamma_{2,2l-1} \\
= & \gamma_{4k-3-4l,4lk-4k-4l^2+l+4}-(-4l+4k-5)\gamma_{2,2l-1} \\
= & \gamma_{-4k+4l+7,-4lk+4l^2+7l-1}. 
\end{align*}
}
Thus the formula (\ref{scc2}) holds for general $l$ and we obtain: 
\[
\text{conj}(T_{2,2k-5}\circ T_{2,2k-7}\circ\cdots T_{2,1})(T_{4k-7,1})=T_{1,k-1}. 
\]
So we can transform $T_s$ into the following sequence: 
\[
(T_{2,1},T_{2,3},\ldots,T_{2,2s-5},T_{1,s-1},T_{1,-1}). 
\]
This completes the proof of Lemma \ref{transT_s}.  \hfill $\square$

\par

\begin{figure}[htbp]
\begin{center}
\includegraphics[width=140mm]{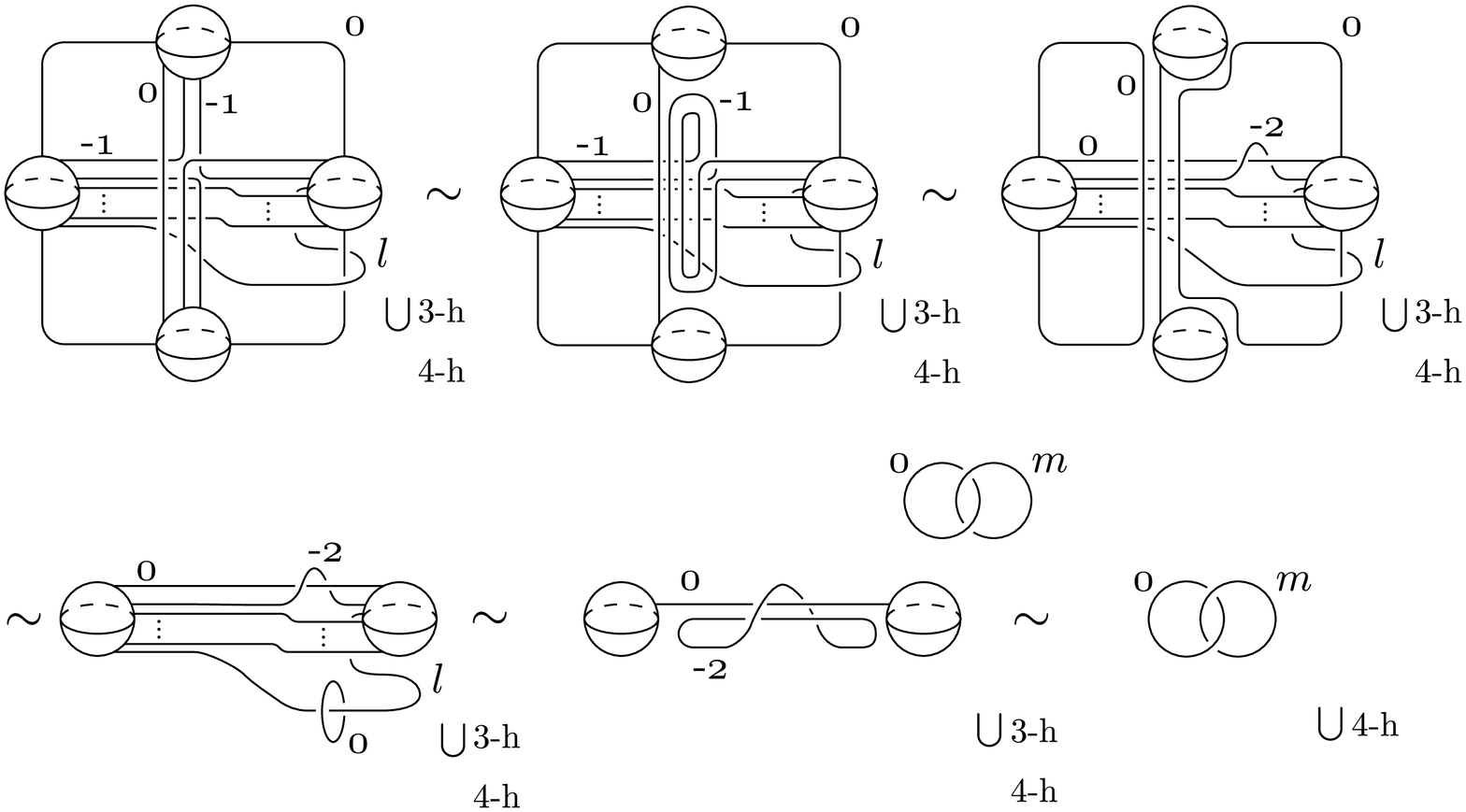}
\end{center}
\caption{A diagram of a SBLF whose Hurwitz system is $T_2$. }
\label{T_2mfd}
\end{figure}

{\it Proof of Theorem \ref{aboutT_s}}: Let $M_s$ be the total space of genus-$1$ SBLF with $W_f=T_s$. 
We first examine the case $s=2$. 
We can describe a diagram of $M_2$ as shown in Figure \ref{T_2mfd}. 
We slide the $2$-handles representing Lefschetz singularities and the $2$-handle of $D^2\times T^2$ to the $2$-handle of the round handle. 
Then we eliminate the obvious canceling pair and slide the $(-2)$-framed knot and the $l$-framed knot to the $0$-framed knot. 
We can change the $l$-framed knot and the $0$-framed meridian of this into a Hopf link by using the $0$-framed meridian. 
We can obtain the last diagram of Figure \ref{T_2mfd} by canceling two pairs of handles. 
Thus we obtain: 
\[
M_2=S, 
\]
where $S$ is $S^2\times S^2$ if $m$ is even and $S^2\tilde{\times}S^2$ if $m$ is odd. 

\begin{figure}[htbp]
\begin{center}
\includegraphics[width=145mm]{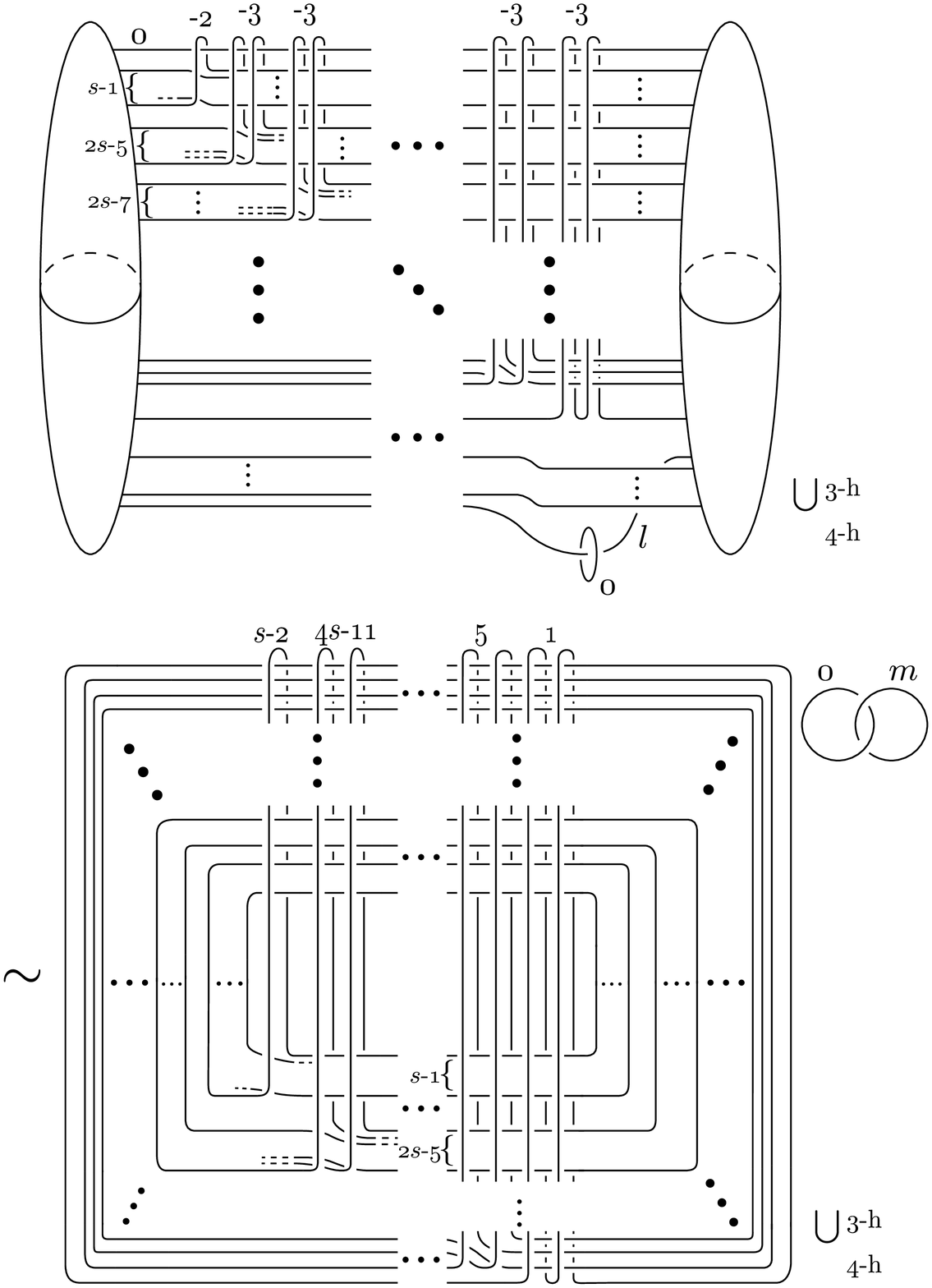}
\end{center}
\caption{}
\label{movesT_s}
\end{figure}

\begin{figure}[htbp]
\begin{center}
\includegraphics[width=145mm]{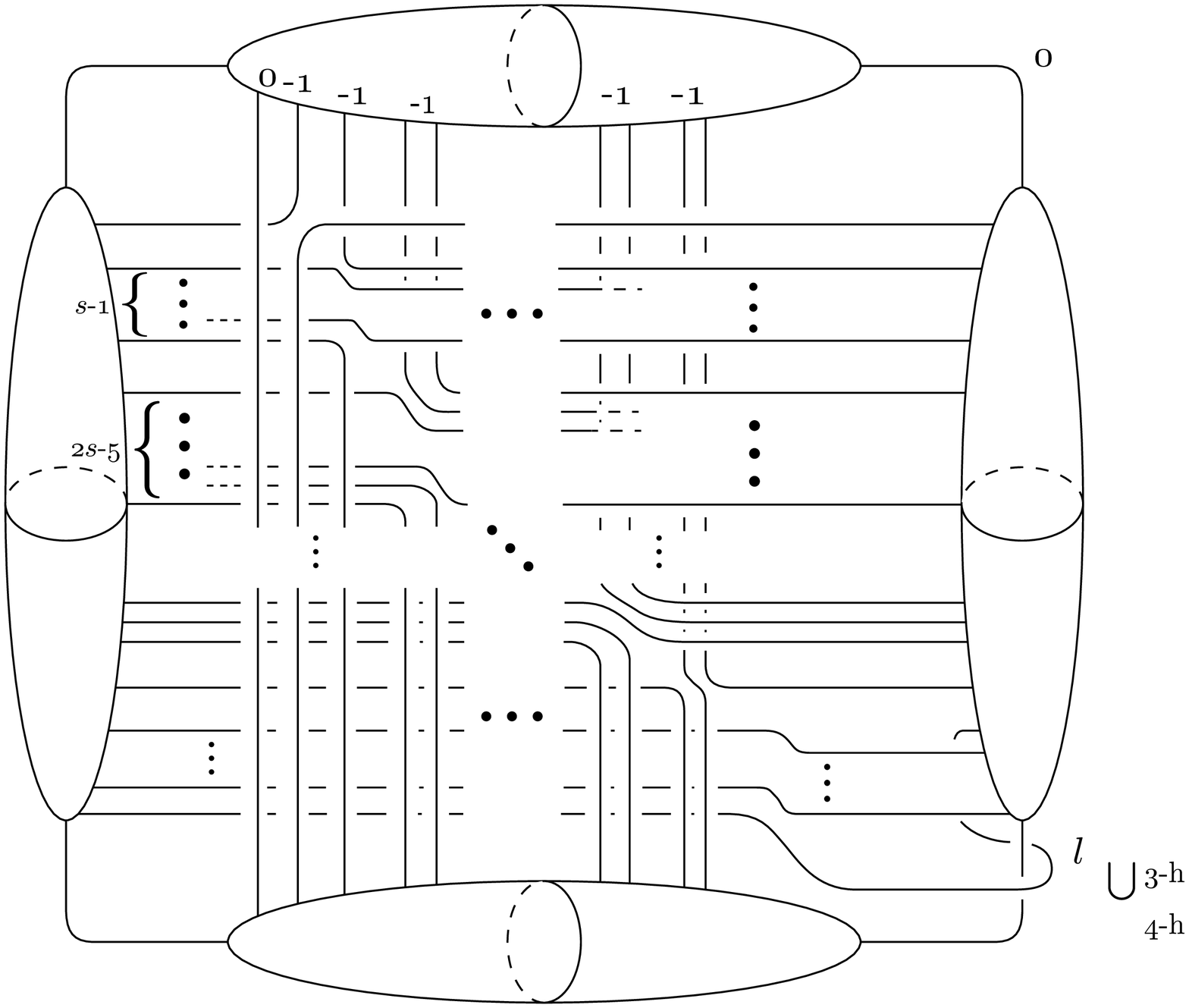}
\end{center}
\caption{A diagram of $M_s$. }
\label{T_smfd}
\end{figure}

For general $s$, we can describe a diagram of $M_s$ as shown in Figure \ref{T_smfd}. 
We can move this diagram to the upper diagram of Figure \ref{movesT_s} by sliding all the $2$-handles describing Lefschetz singularities and the $2$-handle of $D^2\times T^2$ to the $2$-handle of the round handle. 
To obtain the lower diagram of Figure \ref{movesT_s}, we slide the $2$-handles going through the $1$-handle to $0$-framed $2$-handle and eliminate the obvious canceling pair. 

\par

\begin{figure}[htbp]
\begin{center}
\includegraphics[width=145mm]{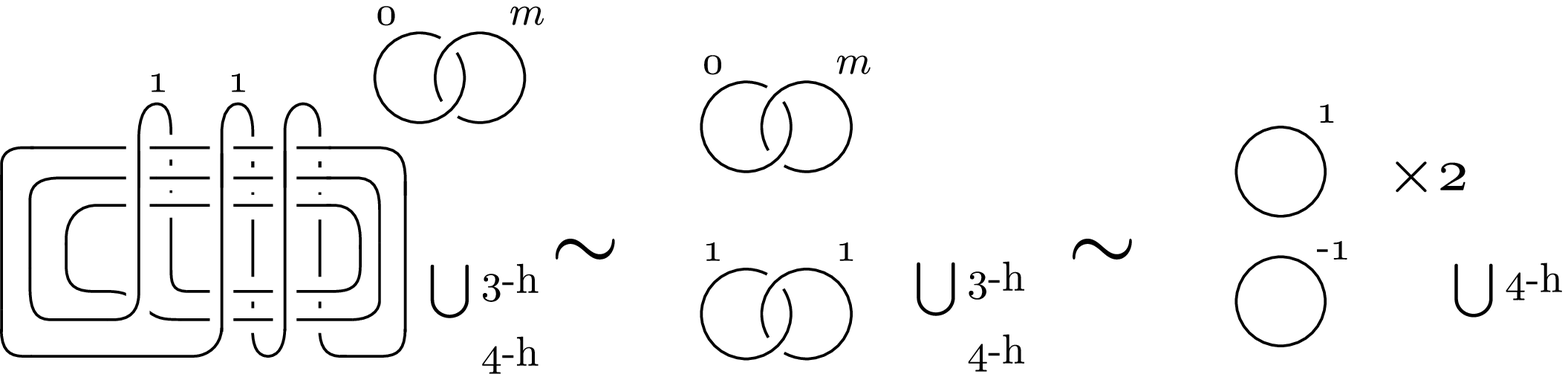}
\end{center}
\caption{}
\label{movesT_3}
\end{figure}

We will prove the conclusion by induction on $s$. 
We first examine the case $s=3$. 
The left diagram in Figure \ref{movesT_3} represents the manifold $M_3$ and we can move it to the center diagram of Figure \ref{movesT_3} by isotopy. 
It is easy to transform the center diagram of Figure \ref{movesT_3} to the right one and we obtain: 
\[
M_3=\sharp 2\mathbb{CP}^2\sharp\overline{\mathbb{CP}^2}. 
\]

\begin{figure}[htbp]
\begin{center}
\includegraphics[width=145mm]{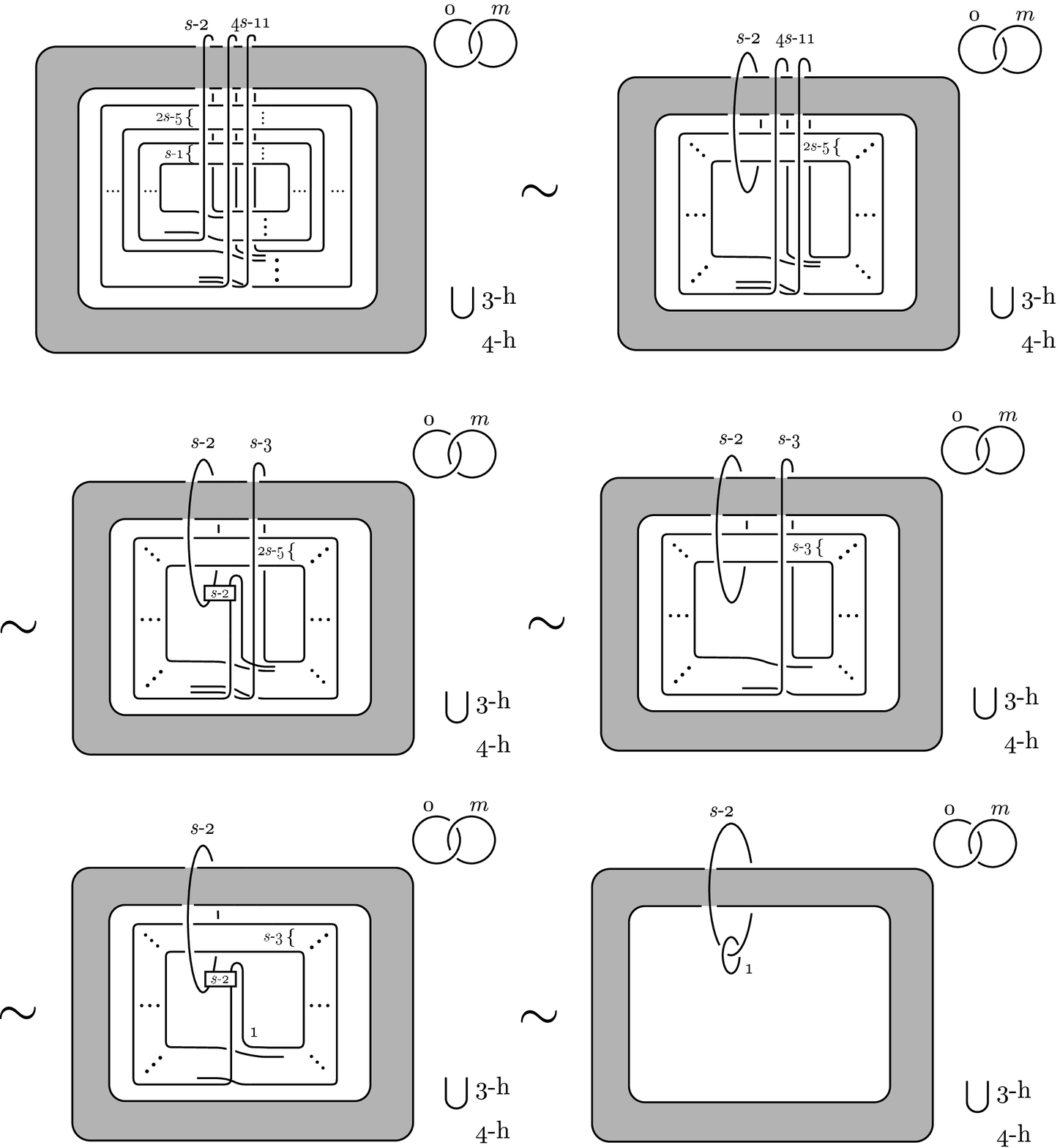}
\end{center}
\caption{Only inner two $2$-handles of the diagram and a Hopf link are described and the other handles are in the shaded part. 
}
\label{movesT_s2}
\end{figure}

To complete the proof, we will move inner two $2$-handles of the lower diagram in Figure \ref{movesT_s} as shown in Figure \ref{movesT_s2}. 
We first move the $(s-2)$-framed $2$-handle by isotopy. 
And we slide the $(4s-11)$-framed $2$-handle to the $(s-2)$-framed $2$-handle twice. 
We then obtain the last diagram of Figure \ref{movesT_s2}. 
Thus we obtain: 
{\allowdisplaybreaks
\begin{align*}
M_s & =M_{s-1}\sharp\mathbb{CP}^2 \\
& =\sharp (s-1)\mathbb{CP}^2\sharp\overline{\mathbb{CP}^2}, 
\end{align*}
}
where the second equality holds by induction hypothesis. 
This completes the proof of Theorem \ref{aboutT_s}.  \hfill $\square$

For sequences of elements of $SL(2,\mathbb{Z})$ $W$ and $W^\prime$, 
we denote by $WW^\prime$ a sequence obtained by standing elements of $W^\prime$ after standing elements of $W$ in a row. 

\begin{thm}\label{connectword}

Let $f:M\rightarrow S^2$ and $g:M^\prime\rightarrow S^2$ be genus-$1$ SBLFs. 

\begin{enumerate}[(1)]

\item If $W_f=S_rW_g$, then $M$ is diffeomorphic to $M^\prime\sharp r\overline{\mathbb{CP}^2}$. 

\item If $W_f=W_gT_s$, then $M$ is diffeomorphic to $M^\prime\sharp S\sharp(s-2)\mathbb{CP}^2$, 
where $S$ is either of the manifolds $S^2\times S^2$ and $S^2\tilde{\times}S^2$. 

\end{enumerate}

\end{thm}

We need the following lemma to prove Theorem \ref{connectword}. 

\begin{lem}\label{ex-lem}

Let $f:M\rightarrow D^2$ be a genus-$1$ LF. 
Suppose that $f$ has at least one Lefschetz singularity and the monodromy of $\partial M$ is equal to $\pm E\in SL(2,\mathbb{Z})$. 
Then every orientation and fiber preserving self-diffeomorphism $\varphi$ of $\partial M$ can be extended to an orientation and fiber preserving self-diffeomorphism $\Phi$ of $M$ up to isotopy. 
i.e. $\Phi|_{\partial M}$ is isotopic to $\varphi$. 

\end{lem}

We remark that the conclusion of Lemma \ref{ex-lem} holds by the following theorem if the monodromy of $\partial M$ is trivial. 

\begin{thm}[\cite{Moish}]\label{Moish-thm}

Let $f:M\rightarrow D^2$ be a genus-$1$ LF and $T\subset M$ a regular fiber of $f$. 
Suppose that $W_f=(X_1,X_2)$. 
Then for every orientation and fiber preserving self-diffeomorphism $\varphi$ of $\partial\nu T$, 
there exists an orientation and fiber preserving self-diffeomorphism $\Phi$ of $M\setminus\text{int}\nu T$ such that $\Phi$ is equal to $\varphi$ on $\partial\nu T$ 
and that $\Phi$ is the identity map on $\partial M$. 

\end{thm}

{\it Proof of Lemma \ref{ex-lem}}: We assume that the induced map on $\partial D^2$ by $\varphi$ is the identity map. 
It is sufficient to prove the conclusion under the above hypothesis. 
Both $M\cup_{id}\overline{M}$ and $M\cup_{\varphi}\overline{M}$ naturally have genus-$1$ ALF structures over $D^2\cup_{id}\overline{D^2}$. 
We denote these ALFs by $f_1$ and $f_2$, respectively. 
Then it is sufficient to prove that there exists an orientation and fiber preserving diffeomorphism $\tilde{\Phi}:M\cup_{\varphi}\overline{M}\rightarrow M\cup_{id}\overline{M}$ 
such that the image under $\Phi$ of $M$ (resp.$\overline{M}$) is in $M$ (resp.$\overline{M}$) and that $\tilde{\Phi}$ is identity map on $\overline{M}$. 

\par

$\partial M$ is isomorphic to $I\times T^2/((1,x)\sim(0,\pm Ex))$ as a torus bundle. 
We fix an identification $\xi_0:\partial M\rightarrow I\times T^2/\sim$ of the two bundles and we identify $\partial D^2$ with $[0,1]/\{0,1\}$ by using $\xi_0$. 

\par

{\it Step.1}: Suppose that an isotopy class of $\varphi:\{\frac{1}{2}\}\times T^2\rightarrow\{\frac{1}{2}\}\times T^2$ is equal to $E\in SL(2,\mathbb{Z})$. 
Then we can move $\varphi$ by isotopy to the map such that the restriction of the map to $[\frac{1}{2}-\varepsilon,\frac{1}{2}+\varepsilon]\times T^2$ is the identity map. 
Let $E\subset D^2$ be an embedded $2$-disk satisfying the following conditions: 

\begin{itemize}

\item $E\cap\partial D^2=[\frac{1}{2}-\varepsilon,\frac{1}{2}+\varepsilon]$; 

\item $E$ contains no critical values of $f$. 

\end{itemize}
\noindent
We denote by $\{p_1,\ldots,p_n\}\subset D^2$ a set of critical values of $f$. 
We put $p_0=\frac{1}{2}$ and $\alpha_1,\ldots,\alpha_n$ be embedded paths in $D^2$, beginning at $p_0$ and otherwise disjoint, intersecting with $\partial D^2$ at only $p_0$, 
connecting $p_0$ to the critical values $p_1,\ldots,p_n$, respectively. 
Let $W_\alpha$ be a regular neighborhood of $\alpha_1\cup\cdots\cup\alpha_n$ such that $W_\alpha\cap\partial D^2\subset E$. 
We denote by $\overline{W_\alpha}$, $\overline{E}$ and $\overline{p_i}$ the image under $id:D^2\rightarrow\overline{D^2}$ of $W_\alpha$, $E$ and $p_i$, respectively. 
Then we can define an orientation and fiber preserving diffeomorphism from $f_1^{-1}(E\cup W_\alpha\cup\overline{W_\alpha}\cup\overline{E})$ to $f_2^{-1}(E\cup W_\alpha\cup\overline{W_\alpha}\cup\overline{E})$
by using the identity maps of $M$ and $\overline{M}$. We denote this map by $\iota$ and $E\cup W_\alpha\cup\overline{W_\alpha}\cup\overline{E}$ by $D_0$. 

\par

\begin{figure}[t]
\begin{center}
\includegraphics[width=145mm]{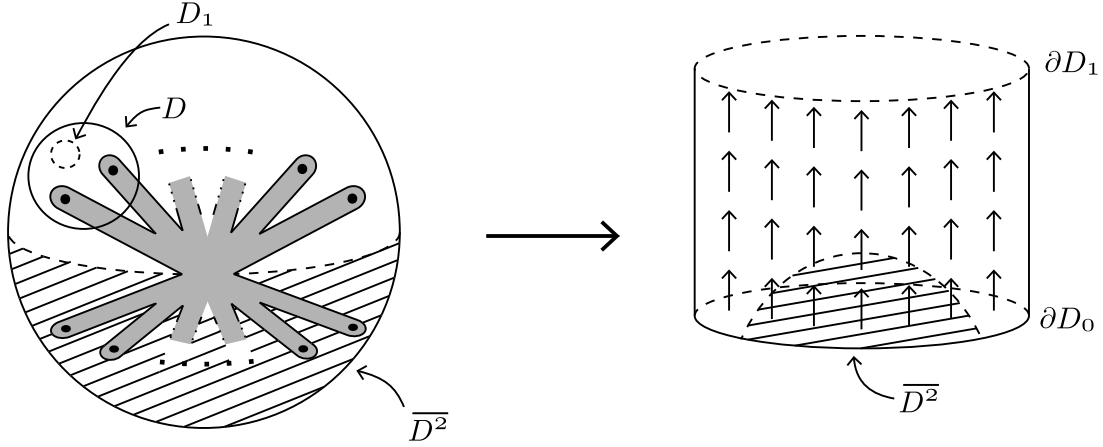}
\end{center}
\caption{Left: the figure of $D^2\cup\overline{D^2}$. 
The shaded part describes $D_0$ and the part with slanted lines describes $\overline{D^2}$. 
Right: the figure of the annulus $D^2\cup_{id}\overline{D^2}\setminus\text{int}(D_0\cup D_1)$. 
The arrows in the figure represent an example of $X_{A_0}$. 
}
\label{annulus}
\end{figure}

There exists a pair of disks $D$, $D_1$ satisfying the following conditions: 

\begin{itemize}

\item $D$ contains two critical values of $f$; 

\item a Hurwitz system of $f|_{f^{-1}(D)}$ is $(X_1,X_2)$; 

\item $D_1$ contains no critical values of $f$; 

\item $D_1\cap(E\cup W_\alpha)=\emptyset$. 

\end{itemize}
Indeed, we can see an example of such a pair in the left diagram of Figure \ref{annulus}. 

\par

$D^2\cup_{id}\overline{D^2}\setminus\text{int}(D_0\cup D_1)$ is an annulus. 
We denote the annulus by $A_0$. 
We can take a vector field $X_{A_0}$ on $A_0$ satisfying the following conditions: 

\begin{itemize}

\item $X_{A_0}$ is non-zero on every point of $A_0$; 

\item each integral curve of $X_{A_0}$ intersects $\partial D_0$ and $\partial D_1$ transversely; 

\item every integral curve of $X_{A_0}$ starting at a point on $\partial D_0$ ends at a point on $\partial D_1$; 

\item for an integral curve $c$ of $X_{A_0}$ and a real number $t$, if $c(t)$ is in $D^2$, $c(s)$ is not in $\overline{D^2}$ for every $s>t$. 

\end{itemize}

For example, we can take such an $X_{A_0}$ as shown in the right figure of Figure \ref{annulus}, 

Let $\mathcal{H}_i$ ($i=1,2$) be horizontal distributions of the restriction of $f_i$ to $M\cup\overline{M}\setminus f_i^{-1}(\{p_1,\ldots,p_n\})$ 
such that $\mathcal{H}_1$ and $\mathcal{H}_2$ are identical on $\overline{M}$. 
We can take vector fields $X_1$ and $X_2$ on $f_1^{-1}(A_0)$ and $f_2^{-1}(A_0)$, respectively, by lifting $X_{A_0}$ by using $\mathcal{H}_1$ and $\mathcal{H}_2$. 
Then we can define a map $\Theta:f_1^{-1}(D^2\cup\overline{D^2}\setminus\text{int}(D_0\cup D_1))\rightarrow f_2^{-1}(D^2\cup\overline{D^2}\setminus\text{int}(D_0\cup D_1))$ as follows: 
\[
\Theta(c_{1,x}(t))=c_{2,x}(t), 
\]
\noindent
where $c_{1,x}$ (resp.$c_{2,x}$) is an integral curve of $X_1$ (resp.$X_2$) starting at $x\in f_1^{-1}(\partial D_0)=f_2^{-1}(\partial D_0)$. 
$\Theta$ is an orientation and fiber preserving diffeomorphism 
and the restriction of $\Theta$ to $f_1^{-1}(\partial D_0)$ is the identity map. 
In particular, the following diagram commutes: 
\[
\begin{CD}
f_1^{-1}(D_0) @>\iota>> f_2^{-1}(D_0) \\
@V id|_{f_1^{-1}(\partial D_0)} VV @VV id|_{f_2^{-1}(\partial D_0)} V \\
f_1^{-1}(D^2\cup \overline{D^2}\setminus\text{int}(D_0\cup D_1)) @>\Theta>> f_2^{-1}(D^2\cup \overline{D^2}\setminus\text{int}(D_0\cup D_1)) \\
\end{CD}
\]
By the above commutative diagram, we obtain the orientation and fiber preserving diffeomorphism $\tilde{\Theta}:f_1^{-1}(D^2\cup \overline{D^2}\setminus\text{int}(D_1))\rightarrow f_2^{-1}(D^2\cup \overline{D^2}\setminus\text{int}(D_1))$. 
The restriction map of $\tilde{\Theta}$ to $\overline{M}$ is the identity map since the horizontal distributions $\mathcal{H}_1$ and $\mathcal{H}_2$ are identical on $\overline{M}$. 

\par

We can obtain by Theorem \ref{Moish-thm} an orientation and fiber preserving diffeomorphism $\tilde{\varphi}:f_1^{-1}(D\setminus\text{int}D_1)\rightarrow f_2^{-1}(D\setminus\text{int}D_1)$ satisfying the following conditions: 

\begin{itemize}

\item $\tilde{\varphi}|_{f_1^{-1}(\partial D_1)}={\tilde{\Theta}|_{f_1^{-1}(\partial D_1)}}^{-1}$; 

\item $\tilde{\varphi}|_{f_1^{-1}(\partial D)}$ is the identity map. 

\end{itemize}

Then we obtain the following commutative diagram: 

\[
\begin{CD}
f_1^{-1}(D^2\cup \overline{D^2}\setminus\text{int}D) @>\tilde{\Theta}>> f_2^{-1}(D^2\cup \overline{D^2}\setminus\text{int}D) \\
@V id|_{f_1^{-1}(\partial D)} VV @VV id|_{f_2^{-1}(\partial D)} V \\
f_1^{-1}(D\setminus\text{int}D_1) @>\tilde{\Theta}\circ\tilde{\varphi}>> f_2^{-1}(D\setminus\text{int}D_1) \\
@V id|_{f_1^{-1}(\partial D_1)} VV @VV id|_{f_2^{-1}(\partial D_1)} V \\
f_1^{-1}(D_1) @>id>> f_2^{-1}(D_1) ,
\end{CD}
\]
where vertical maps are defined only on boundaries. 
We can obtain by the above diagram the orientation and fiber preserving diffeomorphism $\tilde{\Phi}:M\cup_{id}\overline{M}\rightarrow M\cup_{\varphi}\overline{M}$ 
which is the identity map on $\overline{M}$. 
This proves the conclusion. 

\par

{\it Step.2}: We next prove the conclusion for general $\varphi$. 
By Step.1, it is sufficient to prove that, for any element $A\in SL(2,\mathbb{Z})$, there exists an orientation and fiber preserving diffeomorphism $\Psi_A:M\rightarrow M$ 
such that the image under $\Psi_A$ of $\{\frac{1}{2}\}\times T^2$ is in $\{\frac{1}{2}\}\times T^2$ and that the mapping class of the restriction of $\Psi_A$ to $\{\frac{1}{2}\}\times T^2$ is equal to $A$. 
To prove this claim, we use the following theorem. 

\begin{thm}[\cite{Matsumoto2}]\label{Matsumoto-thm}

Let $(B_1,\ldots,B_n)$ be a sequence of elements of $SL(2,\mathbb{Z})$. 
Suppose that each $B_i$ is conjugate to $X_1$ and that $B_1\cdot\cdots\cdot B_n=\pm E$. 
Then, by successive application of elementary transformations, we can change $(B_1,\ldots,B_n)$ into the following sequence: 
\[
(X_1,X_2,X_1,X_2,\ldots,X_1,X_2). 
\]

\end{thm}

We put $\tilde{p_0}\in\text{int}D^2\setminus\{p_1,\ldots,p_n\}$ and let $\gamma$ be an embedded path in $D^2$, connecting $\tilde{p_0}$ to $p_0$, 
intersecting $\partial D^2$ at only $p_0$ transversely. 
We can take an isotopy class of an orientation preserving diffeomorphism between $f^{-1}(\tilde{p_0})$ and $T^2$ by using $\gamma$. 
We fix a representative element $\theta_0:f^{-1}(\tilde{p_0})\rightarrow T^2$ of this class. 
By Theorem \ref{Matsumoto-thm}, we can take embedded paths $\alpha_1,\ldots,\alpha_n$ in $\text{int}D^2$ satisfying the following conditions: 

\begin{enumerate}[(a)]

\setlength{\itemindent}{10pt}

\item $\alpha_i$ starts at $\tilde{p_0}$ and connects $\tilde{p_0}$ to $p_i$; 

\item $\alpha_i\cap\alpha_j=\{\tilde{p_0}\}$ if $i\neq j$; 

\item $\alpha_i\cap\gamma=\{\tilde{p_0}\}$; 

\item $\gamma,\alpha_1,\ldots,\alpha_n$ appear in this order when we travel counterclockwise around $\tilde{p_0}$; 

\item a Hurwitz system of $f$ determined by $\theta_0,\alpha_1,\ldots,\alpha_n$ is 
\[
(X_1,X_2,\ldots,X_1,X_2). 
\]

\end{enumerate}

Similarly, we can take embedded paths $\beta_1,\ldots,\beta_n$ in $\text{int}D^2$ so that they satisfy the following conditions: 

\begin{enumerate}[(a)$^\prime$]

\setlength{\itemindent}{10pt}

\item $\beta_i$ starts at $\tilde{p_0}$ and connects $\tilde{p_0}$ to $p_{k_i}$, where $\{p_{k_1},\ldots,p_{k_n}\}=\{p_1,\ldots,p_n\}$; 

\item $\beta_i\cap\beta_j=\{\tilde{p_0}\}$ if $i\neq j$; 

\item $\beta_i\cap\gamma=\{\tilde{p_0}\}$; 

\item $\gamma,\beta_1,\ldots,\beta_n$ appear in this order when we travel counterclockwise around $\tilde{p_0}$; 

\item a Hurwitz system of $f$ determined by $\theta_0,\beta_1,\ldots,\beta_n$ is 
\[
(A^{-1}X_1A,A^{-1}X_2A,\ldots,A^{-1}X_1A,A^{-1}X_2A). 
\]

\end{enumerate}

We denote by $E_0$ a disk neighborhood of $\tilde{p_0}$ in $\text{int}D^2\setminus\{p_1,\ldots,p_n\}$. 
Let $\tilde{\theta_0}:f^{-1}(E_0)\rightarrow E_0\times T^2$ be an extension of an identification $\theta_0$. 
We define a diffeomorphism $\Psi_1:f^{-1}(E_0)\rightarrow f^{-1}(E_0)$ as follows: 
\[
\Psi_1(p,x)=(p,Ax), 
\]
where $(p,x)\in D^2\times T^2$. 
Let $W_\alpha$ and $W_\beta$ be regular neighborhoods of $\alpha_1\cup\cdots\cup\alpha_n$ and $\beta_1\cup\cdots\cup\beta_n$, respectively. 
We denote by $V_1,\ldots,V_n$ sufficiently small neighborhoods in $\text{int}D^2$ of $p_1,\ldots,p_n$, respectively. 
We can easily extend $\Psi_1$ to: 
\[
\Psi_2:f^{-1}(E_0\cup W_\alpha\setminus(V_0\cup\cdots\cup V_n))\rightarrow f^{-1}(E_0\cup W_\beta\setminus(V_0\cup\cdots\cup V_n)). 
\]

Moreover, we can take $\Psi_2$ so that the image under $\Psi_2$ of the inverse image of $\alpha_i$ is in the inverse image of $\beta_i$. 
Then the image under $\Psi_2$ of a vanishing cycle in a regular neighborhood of $f$ near $f^{-1}(p_i)$ is isotopic to a vanishing cycle in a regular neighborhood of $f$ near $f^{-1}(p_{k_i})$. 
Thus $\Psi_2$ can be extended to: 
\[
\Psi_3:f^{-1}(E_0\cup W_\alpha)\rightarrow f^{-1}(E_0\cup W_\beta). 
\]
We can easily extend this map to a self-diffeomorphism of $M$ and this extended self-diffeomorphism satisfies the desired condition. 
This completes the proof of Lemma \ref{ex-lem}.  \hfill $\square$

{\it Proof of Theorem \ref{connectword}}: By Lemma \ref{ex-lem}, without loss of generality, we can assume that 
the attaching circle of the $2$-handle of the round handle in $M$ is isotopic to $\gamma_{1,0}$ in a regular fiber. 

\par

\begin{figure}[htbp]
\begin{center}
\includegraphics[width=65mm]{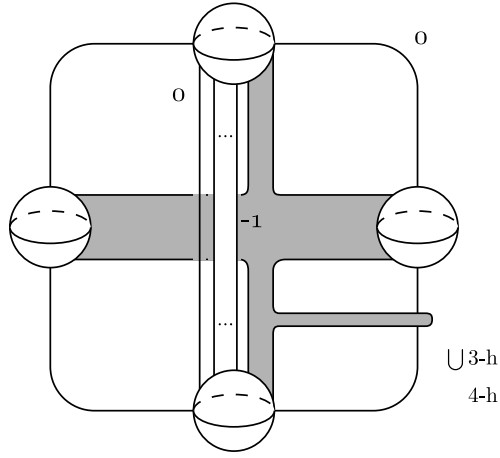}
\end{center}
\caption{A diagram of a total space of a SBLF with associated sequence $S_rW_g$. 
Only handles of $D^2\times T^2$, round handle and $2$-handles corresponding to elements in $S_r$ are described 
and the other handles are in a shaded part. }
\label{S_rWmfd}
\end{figure}

We first prove the statement (1). 
We can describe a diagram of $M$ as shown in Figure \ref{S_rWmfd}. 
We can divide the diagram into $(-1)$-framed unknots and a diagram of $M^\prime$ by sliding $2$-handles corresponding to elements in $S_r$ to the $2$-handle of the round handle. 
This proves the statement (1). 

\begin{figure}[htbp]
\begin{center}
\includegraphics[width=140mm]{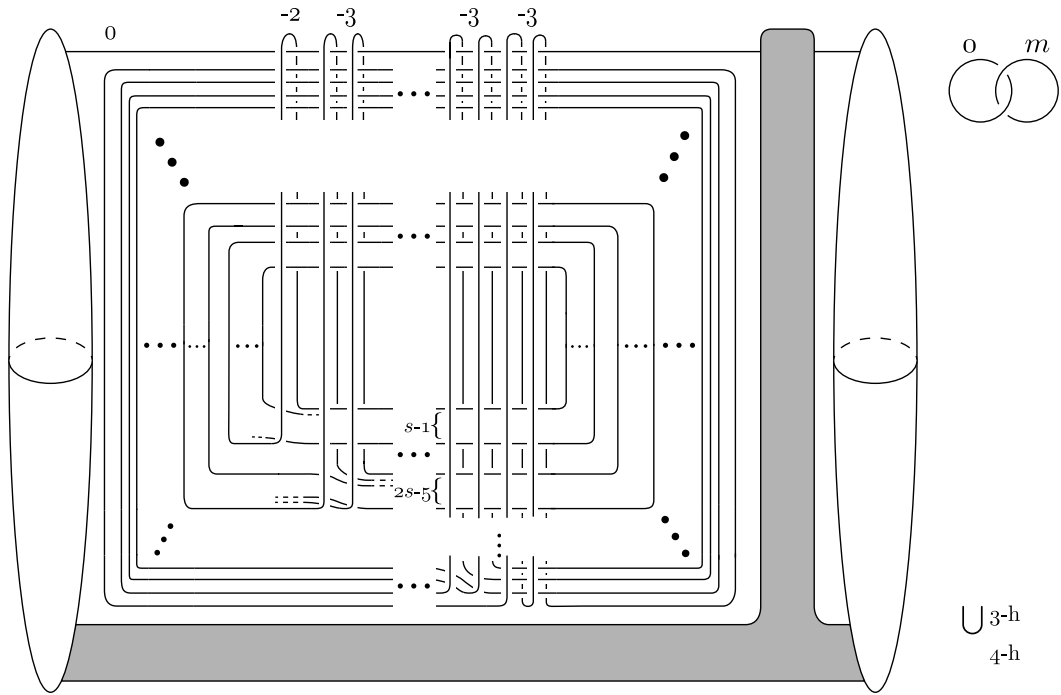}
\end{center}
\caption{A diagram of a total space of a SBLF whose Hurwitz system is $W_gT_s$. 
Only handles of $D^2\times T^2$, round handle and $2$-handles corresponding to elements in $T_s$ are described 
and the $2$-handle of the round handle is eliminated. }
\label{WT_smfd}
\end{figure}

\begin{figure}[htbp]
\begin{center}
\includegraphics[width=140mm]{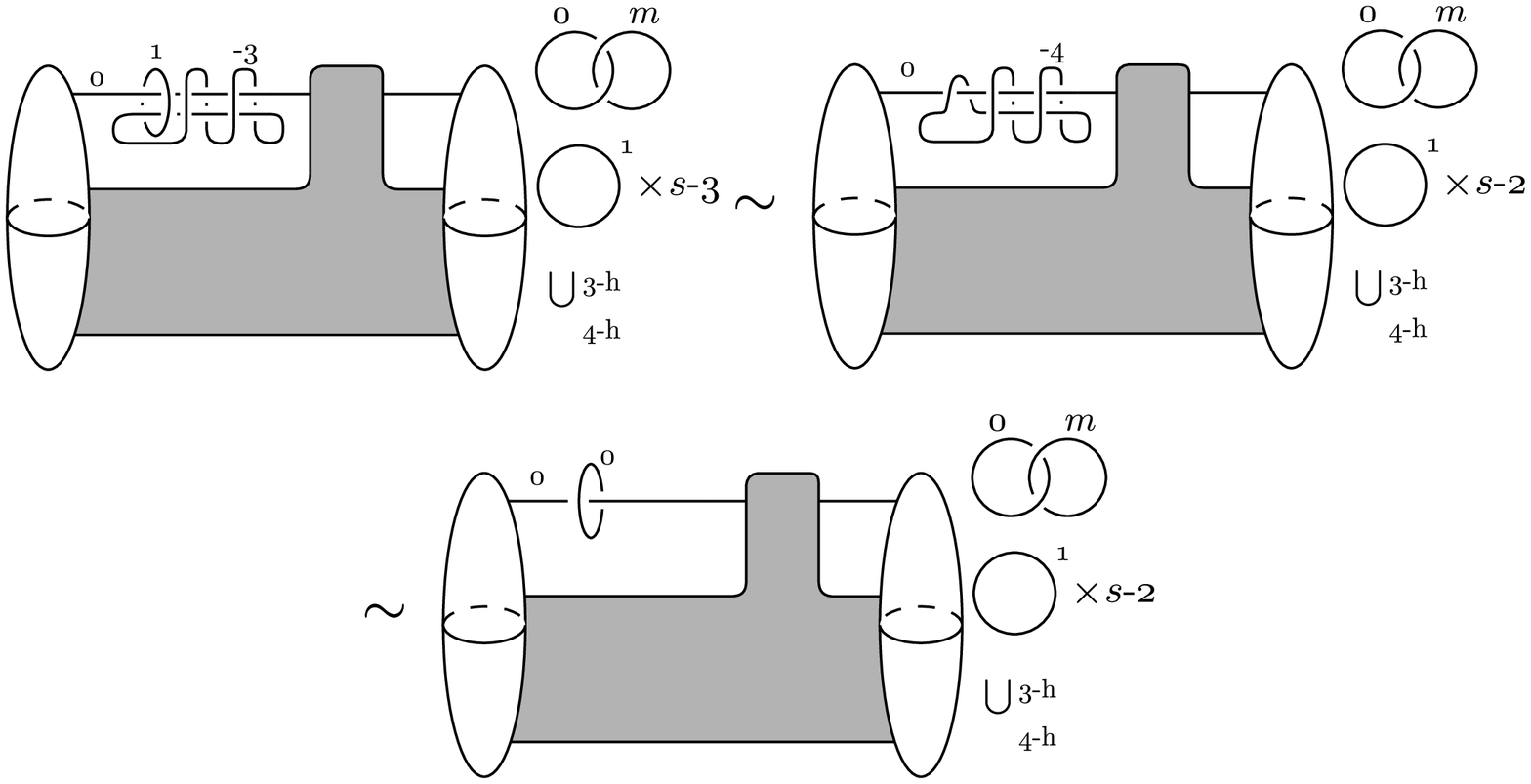}
\end{center}
\caption{}
\label{movesWT_s}
\end{figure}

We next prove the statement (2) by induction on the number of elements in $W_g$. 
The statement is true if $W_g=\emptyset$ by Theorem \ref{aboutT_s}. 
If $W_g$ contains an element $X_1$, the statement holds by the statement (1) and the induction hypothesis. 
Suppose that $W_g\neq\emptyset$ and $W_g$ does not contain $X_1$. 
By Theorem \ref{main in chart sec}, $W_g$ is equal to the following normal form: 
\[
(X_1^{n_1}X_2X_1^{-n_1},\ldots,X_1^{n_u}X_2X_1^{-n_u}). 
\]
Then a diagram of $M$ can be described as shown in Figure \ref{WT_smfd}. 
By applying moves as we performed to prove Theorem \ref{aboutT_s}, we obtain the first diagram in Figure \ref{movesWT_s}. 
After handle slidings, we obtain the last diagram in Figure \ref{movesWT_s}. 
There exists a $2$-handle whose attaching circle intersects the belt sphere of the $1$-handle geometrically once in the shaded part. 
We denote the $2$-handle by $H$. 
By sliding the $0$-framed $2$-handle intersecting the belt sphere of the $1$-handle geometrically once to $H$ and using the $0$-framed meridian, 
we can divide the diagram into a diagram of $S$, $s-2$ $(-1)$-framed unknots and a diagram of $M^\prime$. 
This proves the statement (2). 

\par

Combining the conclusions obtained above, we complete the proof of Theorem \ref{connectword}. \hfill $\square$

We can prove the following theorem by an argument similar to that in the proof of Theorem \ref{connectword}. 

\begin{thm}\label{indef}

Let $f:M\rightarrow S^2$ be a genus-$1$ SBLF. 
Suppose that $W_f=S_rT(n_1,\ldots,n_s)$ and $s>0$. 
Then there exists a $4$-manifold $M^\prime$ such that $M=M^\prime\sharp S$, where $S$ is either of the manifolds $S^2\times S^2$ and $S^2\tilde{\times}S^2$. 

\end{thm}

\begin{figure}[htbp]
\begin{center}
\includegraphics[width=60mm]{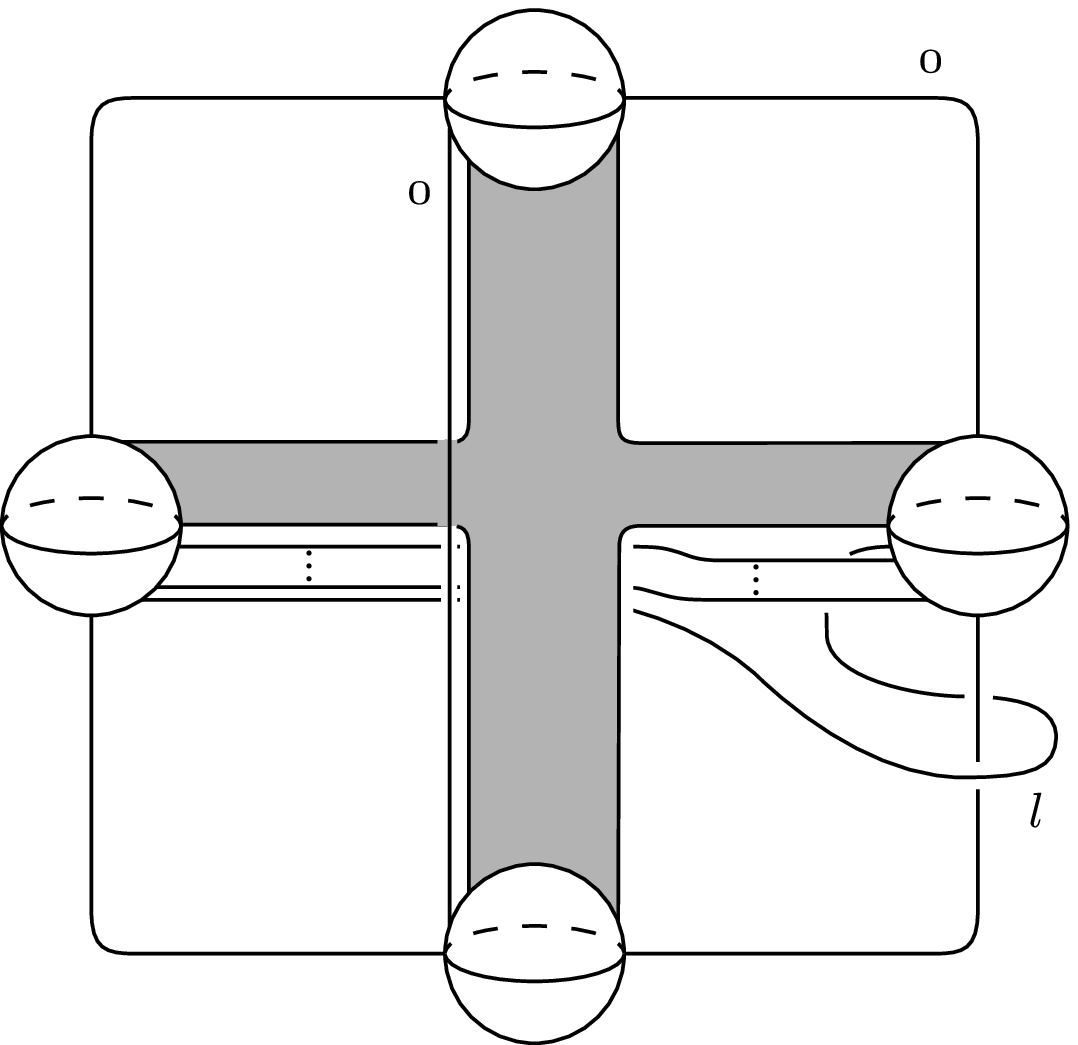}
\end{center}
\caption{}
\label{indefinite}
\end{figure}

{\it Proof}: By assumption about $W_f$, we can draw a diagram of $M$ as shown in Figure \ref{indefinite}. 
We first slide the $2$-handle of $D^2\times T^2$ to the $2$-handle of the round $2$-handle twice. 
Then the $2$-handle of $D^2\times T^2$ becomes a $0$-framed meridian of $l$-framed $2$-handle in the diagram. 
Since we assume that $s>0$, there exists a $2$-handle $H$ in a shaded part such that an attaching circle of $H$ goes through the $1$-handle that the $l$-framed $2$-handle goes through. 
By sliding $l$-framed $2$-handle to $H$ and using the $0$-framed meridian, we can change the $l$-framed $2$-handle and the $0$-framed meridian into the Hopf link. 
Thus we complete the proof of Theorem \ref{indef}.  \hfill $\square$

\begin{cor}\label{definite}

Simply connected $4$-manifolds with positive definite intersection form cannot admit any genus-$1$ SBLF structures except $S^4$. 
Especially, $k\mathbb{CP}^2$ cannot admits any such fibrations for any positive integers $k$. 

\end{cor}

{\it Proof}: If such a $4$-manifold admitted a genus-$1$ SBLF structure, the $4$-manifold would contain at least one $S^2\times S^2$ or $S^2\tilde{\times}S^2$ connected summand by Theorem \ref{indef}. 
This contradicts the hypothesis about the intersection form. \hfill $\square$

We end this section with the proof of Main Theorem A. 
\\[5pt]
\par
\begin{figure}[htbp]
\begin{center}
\includegraphics[width=80mm]{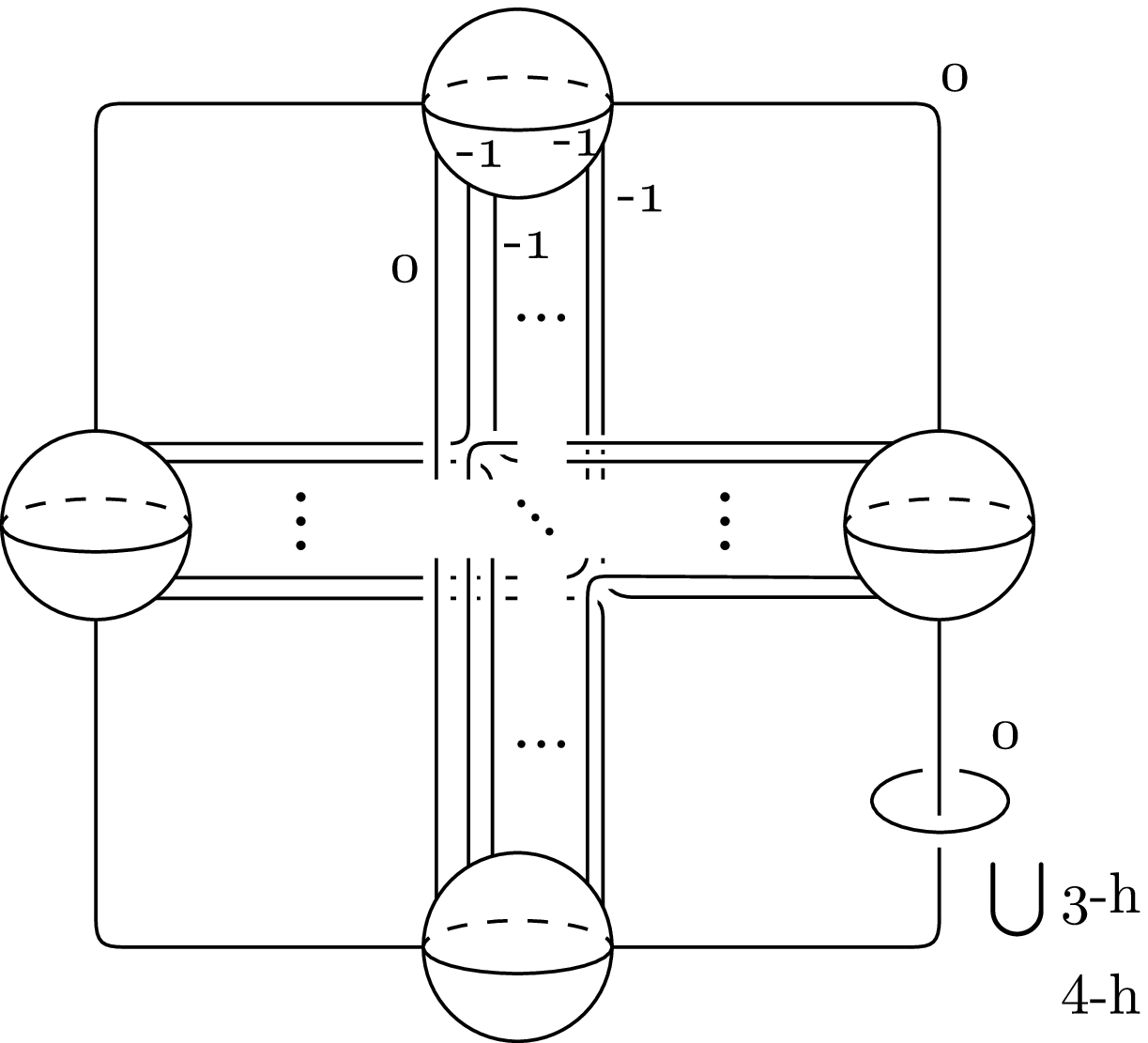}
\end{center}
\caption{}
\label{kS^2timesS^2}
\end{figure}

{\it Proof of Main Theorem A}: By Theorem \ref{aboutS_r}, \ref{aboutT_s} and \ref{connectword}, it is sufficient to prove only that $\sharp kS^2\times S^2$ admits a genus-$1$ SBLF structure. 
A diagram in Figure \ref{kS^2timesS^2} represents the total space of a genus-$1$ SBLF whose monodromy representation is represented by the sequence $kT_2$. 
We denote this manifold by $M_k$. 

We can prove that $M_k$ is diffeomorphic to $S^2\times S^2$ by applying same moves as we apply in the proof of the statement (2) of Theorem \ref{connectword}. 
This completes the proof of Main Theorem A.  \hfill $\square$

\section{Proof of Main Theorem B}

We will use some properties of $PSL(2,\mathbb{Z})=SL(2,\mathbb{Z})/\{\pm E\}$ for proving Main Theorem B. 
So we first review these properties before proving the theorem. 

\par

We define matrices $A$ and $B$ in $SL(2,\mathbb{Z})$ as follows: 

\[
A=
\begin{pmatrix}
0 & -1 \\
1 & 1 
\end{pmatrix}
\hspace{.2em},\hspace{.2em}
B=
\begin{pmatrix}
1 & 2 \\
-1 & -1 
\end{pmatrix}.
\]
Then $A$ and $B$ generate $SL(2,\mathbb{Z})$. 
In fact, both $X_1$ and $X_2$ are represented by $A$ and $B$ as follows: 
\[
X_1=ABA\hspace{.2em},\hspace{.2em}X_2=BA^2. 
\]
Let $a$, $b$, $x_1$ and $x_2$ be elements of $PSL(2,\mathbb{Z})$ represented by $A$, $B$, $X_1$ and $X_2$, respectively. 
Then $PSL(2,\mathbb{Z})$ has the following finite presentation \cite{MKS}: 
\[
PSL(2,\mathbb{Z})=<a,b|a^3,b^2>. 
\]
Especially, $PSL(2,\mathbb{Z})$ is isomorphic to a free product $\mathbb{Z}/3\ast\mathbb{Z}/2$. 

\par

The sequence $(w_1,\ldots,w_n)$ of elements of $PSL(2,\mathbb{Z})$ is called {\it reduced} 
if the set $\{w_i,w_{i+1}\}$ is equal to either of the sets $\{a,b\}$ and $\{a^2,b\}$ for each $i\in\{1,\ldots,n-1\}$. 
The proof of Main Theorem B is based on the following theorem. 

\begin{thm}[\cite{MKS}]\label{key thm of B}

For every element $g$ of $PSL(2,\mathbb{Z})$, there exists a reduced sequence $(w_1,\ldots,w_n)$ of $PSL(2,\mathbb{Z})$ such that $g=w_1\cdot\cdots\cdot w_n$. 
Moreover, such a sequence is unique. 

\end{thm}

For an element $g\in PSL(2,\mathbb{Z})$, we take the reduced sequence $(w_1,\ldots,w_n)$ such that $g=w_1\cdot\cdots\cdot w_n$. 
We define an element $t(g)$ as follows: 
\[
t(g)=w_n\cdot\cdots\cdot w_1. 
\]
Then a map $t:PSL(2,\mathbb{Z})\rightarrow PSL(2,\mathbb{Z})$ is well-defined by Theorem \ref{key thm of B} and it is an antihomomorphism. 
Since $x_1=aba$ and $x_2=ba^2$, we obtain $t(x_1)=x_1$ and $t(x_2)=x_1^{-1}x_2x_1$. 

\par

Now we are ready to prove Main Theorem B. 
\\[5pt]
{\it Proof of Main Theorem B}: Let $f:M\rightarrow S^2$ be a genus-$1$ SBLF with non-empty round singular locus and with $r$ Lefschetz singularities. 
Let $W_f$ be a Hurwitz system of $f$ which satisfies $w(W_f)=\pm X_1^m$, where $m\in\mathbb{Z}$. 
We can assume by Theorem \ref{main in chart sec} that $W_f=S_kT(n_1,\ldots,n_{r-k})$. 
\par

For the case $r=0$, the conclusion holds directly by Theorem \ref{aboutS_r}. 

\par

For the case $r=1$, $W_f$ is either of the sequence $S_1$ and $T(n_1)$. 
If $W_f$ were $T(n_1)$, $w(W_f)$ would not be equal to $\pm X_1^m$. 
So $W_f$ is equal to $S_1$ and the conclusion holds by Theorem \ref{aboutS_r}. 

\par

For the case $r=2$, $W_f$ is equal to $S_2$, $S_1T(n_1)$ or $T(n_1,n_2)$. 
$W_f$ is not equal to $S_1T(n_1)$ since $w(W_f)=\pm X_1^m$. 
If $W_f$ is equal to $S_2$, then $M$ is diffeomorphic to $\sharp 2\overline{\mathbb{CP}^2}$, $S^1\times S^3\sharp S\sharp 2\overline{\mathbb{CP}^2}$ or $L\sharp 2\overline{\mathbb{CP}^2}$. 
If $W_f$ is equal to $T(n_1,n_2)$ and $n_1-n_2=1$, we can change $W_f$ by elementary transformations as follows: 
{\allowdisplaybreaks
\begin{align*}
W_f = & (X_1^{-n_1}X_2X_1^{n_1},X_1^{-n_2}X_2X_1^{n_2}) \\
\mapsto & ((X_1^{-n_1}X_2X_1^{n_1})(X_1^{-n_2}X_2X_1^{n_2})(X_1^{-n_1}X_2^{-1}X_1^{n_1}),X_1^{-n_1}X_2X_1^{n_1}) \\
= & (X_1^{-n_1}(X_2X_1X_2X_1^{-1}X_2^{-1})X_1^{n_1},X_1^{-n_1}X_2X_1^{n_1}) \\
= & (X_1^{-n_1}X_1X_1^{n_1},X_1^{-n_1}X_2X_1^{n_1}) \\
= & (X_1,X_1^{-n_1}X_2X_1^{n_1}) \\
= & S_1T(n_1). 
\end{align*}
}
This contradiction says that $n_1-n_2\neq 1$. 
Since $w(W_f)=\pm X_1^m$, the following equation holds: 
{\allowdisplaybreaks
\begin{align}
\label{r=2}
& x_1^{-n_1}x_2x_1^{n_1}x_1^{-n_2}x_2x_1^{n_2}=x_1^m \nonumber \\
\Rightarrow & x_2x_1^{n_1-n_2}x_2=x_1^n, 
\end{align}
}
where $n=m+n_1-n_2$. 
The left side of Equation (\ref{r=2}) is represented by $a$ and $b$ as follows: 
{\allowdisplaybreaks
\begin{align*}
x_2x_1^{n_1-n_2}x_2 & = \begin{cases}
ba^2\cdot a^2(ba)^{-n_1+n_2-1}ba^2ba^2 & (n_1-n_2\leq 0) \\
ba^2\cdot a(ba^2)^{n_1-n_2-1}baba^2 & (n_1-n_2\geq 2) 
\end{cases} \\
& = \begin{cases}
(ba)^{-n_1+n_2}ba^2ba^2 & (n_1-n_2\leq 0), \\
a^2(ba^2)^{n_1-n_2-2}baba^2 & (n_1-n_2\geq 2). 
\end{cases}
\end{align*}
}
$x_1^n$ is represented by $a$ and $b$ as follows: 
{\allowdisplaybreaks
\begin{align}
\label{x_1^n}
x_1^n & = \begin{cases}
a(ba^2)^{n-1}ba & (n\geq 1), \\
a^2(ba)^{-n-1}ba^2 & (n\leq -1). 
\end{cases}
\end{align}
}
By Theorem \ref{key thm of B}, we obtain $n_1-n_2=2$. 
Then we can change $W_f$ into $T_2$ by applying simultaneous conjugation. 
Thus $M$ is diffeomorphic to either of the manifolds $S^2\times S^2$ and $S^2\tilde{\times}S^2$. 

\par

For the case $r=3$, $W_f$ is equal to $S_3$, $S_2T(n_1)$, $S_1T(n_1,n_2)$ or $T(n_1,n_2,n_3)$. 
The case $W_f=S_2T(n_1)$ never occurs since $w(W_f)=\pm X_1^m$. 
If $W_f=S_3$ or $S_1T(n_1,n_2)$, $4$-manifolds that $M$ can be is well-known by Theorem \ref{aboutS_r}, \ref{aboutT_s}, \ref{connectword} and the argument for the case $r=2$ 
and these are contained in a family of $4$-manifolds in the statement of Main Theorem B. 
So all we need to examine is the case $W_f=T(n_1,n_2,n_3)$. 
We can assume that $n_i-n_{i+1}\neq 1$ ($i=1,2$) by the same argument as in the case $r=2$. 
If $n_i-n_{i+1}=2$ for $i=1$ or $2$, then $X_1^{-n_i}X_2X_1^{n_i}X_1^{-n_{i+1}}X_2X_1^{n_{i+1}}=-X_1^{-4}$ and we can change $W_f$ into $T(\tilde{n})T_2$ by elementary transformations and simultaneous conjugations, where $\tilde{n}$ is an integer. 
This contradicts $w(W_f)=\pm X_1^m$. 
So $n_i-n_{i+1}\neq 2$ for $i=1,2$. 
Since $W_f=\pm X_1^m$, the following equation holds: 
{\allowdisplaybreaks
\begin{align}
\label{r=3}
& x_1^{-n_1}x_2x_1^{n_1}x_1^{-n_2}x_2x_1^{-n_3}x_2x_1^{n_3}=x_1^m \nonumber \\
\Rightarrow & x_2x_1^{n_1-n_2}x_2x_1^{n_2-n_3}x_2=x_1^n, 
\end{align}
}
where $n=m+n_1-n_3$. 
If $n_1-n_2\leq 0$, Then the left side of Equation (\ref{r=3}) is represented by $a$ and $b$ as follows: 
{\allowdisplaybreaks
\begin{align*}
x_2x_1^{n_1-n_2}x_2x_1^{n_2-n_3}x_2 & =\begin{cases}
ba^2\cdot a^2(ba)^{-n_1+n_2-1}ba^2\cdot ba^2\cdot a(ba^2)^{n_2-n_3-1}\cdot ba^2 & (n_2-n_3\geq 3) \\
ba^2\cdot a^2(ba)^{-n_1+n_2-1}ba^2\cdot ba^2\cdot a^2(ba)^{-n_2+n_3-1}\cdot ba^2 & (n_2-n_3\leq 0) 
\end{cases} \\
& = \begin{cases}
(ba)^{-n_1+n_2}ba\cdot (ba^2)^{n_2-n_3-3}ba^2 & (n_2-n_3\geq 3), \\
(ba)^{-n_1+n_2}ba^2(ba)^{-n_2+n_3}ba^2 & (n_2-n_3\leq 0).
\end{cases}
\end{align*}
}
This would not be equal to $x_1^n$ by (\ref{x_1^n}) and Theorem \ref{key thm of B}. 
Thus we obtain $n_1-n_2\geq 3$. 
By operating a map $t$ to both sides of Equation (\ref{r=3}), we obtain: 
{\allowdisplaybreaks
\begin{align*}
& (x_1^{-1}x_2x_1)x_1^{n_2-n_3}(x_1^{-1}x_2x_1)x_1^{n_1-n_2}(x_1^{-1}x_2x_1)=x_1^n \\
\Rightarrow & x_2x_1^{n_2-n_3}x_2x_1^{n_1-n_2}x_2=x_1^n. 
\end{align*}
}
So we obtain $n_2-n_3\geq 3$ by the same argument as above. 
Thus the left side of Equation (\ref{r=3}) is represented by $a$ and $b$ as follows: 
{\allowdisplaybreaks
\begin{align*}
& x_2x_1^{n_1-n_2}x_2x_1^{n_2-n_3}x_2 \\
= & ba^2\cdot a(ba^2)^{n_1-n_2-1}ba\cdot ba^2\cdot a(ba^2)^{n_2-n_3-1}ba\cdot ba^2 \\
= & a^2(ba^2)^{n_1-n_2-2}ba\cdot a^2(ba^2)^{n_2-n_3-2}baba^2 \\
= & a^2(ba^2)^{n_1-n_2-3}ba(ba^2)^{n_2-n_3-3}baba^2. 
\end{align*}
}
By Theorem \ref{key thm of B}, we obtain $n_1-n_2=n_2-n_3=3$. 
Then we can change $W_f$ into $T_3$ by applying simultaneous conjugation and $M$ is diffeomorphic to $\sharp 2\mathbb{CP}^2\sharp\overline{\mathbb{CP}^2}$. 

\par

For the case $r=4$, $W_f$ is equal to $S_4$, $S_3T(n_1)$, $S_2T(n_1,n_2)$, $S_1T(n_1,n_2,n_3)$ or $T(n_1,n_2,n_3,n_4)$. 
But all we need to examine is the case $W_f=T(n_1,n_2,n_3,n_4)$ as the case $r=3$. 
We can assume that $n_i-n_{i+1}\neq 1,2$ for $i=1,2,3$ by the same argument as in the case $r=3$. 
If $n_i-n_{i+1}=n_{i+1}-n_{i+2}=3$ for $i=1$ or $2$, then we can change $W_f$ into $T(\tilde{n})T_3$ by elementary transformations and simultaneous conjugations, where $\tilde{n}$ is an integer. 
This contradicts $w(W_f)=\pm X_1^m$. 
So we obtain $(n_i-n_{i+1},n_{i+1}-n_{i+2})\neq (3,3)$ for $i=1,2$. 
Since $W_f=\pm X_1^m$, we obtain the following equation: 
{\allowdisplaybreaks
\begin{align}
\label{r=4}
& x_1^{-n_1}x_2x_1^{n_1}x_1^{-n_2}x_2x_1^{n_2}x_1^{-n_3}x_2x_1^{n_3}x_1^{-n_4}x_2x_1^{n_4}=x_1^m \nonumber \\
\Rightarrow & x_2x_1^{n_1-n_2}x_2x_1^{n_2-n_3}x_2x_1^{n_3-n_4}x_2=x_1^n, 
\end{align}
}
where $n=m+n_1-n_4$. 
If $n_1-n_2\leq 0$, the left side of Equation (\ref{r=4}) is represented by $a$ and $b$ as follows: 
{\allowdisplaybreaks
\begin{align*}
& x_2x_1^{n_1-n_2}x_2x_1^{n_2-n_3}x_2x_1^{n_3-n_4}x_2 \\
= & \begin{cases}
ba^2\cdot a^2(ba)^{-n_1+n_2-1}ba^2\cdot ba^2\cdot a(ba^2)^{n_2-n_3-1}ba\cdot ba^2 \\
\hspace{2em}\cdot a(ba^2)^{n_3-n_4-1}ba\cdot ba^2 & (n_2-n_3\geq 3, n_3-n_4\geq 3) \\
ba^2\cdot a^2(ba)^{-n_1+n_2-1}ba^2\cdot ba^2\cdot a^2(ba)^{-n_2+n_3-1}ba^2 \cdot ba^2 \\
\hspace{2em}\cdot a(ba^2)^{n_3-n_4-1}ba\cdot ba^2 & (n_2-n_3\leq 0, n_3-n_4\geq 3) \\
ba^2\cdot a^2(ba)^{-n_1+n_2-1}ba^2\cdot ba^2\cdot a(ba^2)^{n_2-n_3-1}ba\cdot ba^2 \\
\hspace{2em}\cdot a^2(ba)^{-n_3+n_4-1}ba^2\cdot ba^2 & (n_2-n_3\geq 3, n_3-n_4\leq 0) \\
ba^2\cdot a^2(ba)^{-n_1+n_2-1}ba^2\cdot ba^2\cdot a^2(ba)^{-n_2+n_3-1}ba^2\cdot ba^2 \\
\hspace{2em}\cdot a^2(ba)^{-n_3+n_4-1}ba^2\cdot ba^2 & (n_2-n_3\leq 0, n_3-n_4\leq 0)
\end{cases} \\
= & \begin{cases}
(ba)^{-n_1+n_2}ba(ba^2)^{n_2-n_3-3}ba(ba^2)^{n_3-n_4-3}baba^2 & (n_2-n_3\geq 3, n_3-n_4\geq 3), \\
(ba)^{-n_1+n_2}ba^2(ba)^{-n_2+n_3}ba(ba^2)^{n_3-n_4-2}baba^2 & (n_2-n_3\leq 0, n_3-n_4\geq 3), \\
(ba)^{-n_1+n_2}ba(ba^2)^{n_2-n_3-2}(ba)^{-n_3+n_4+1}ba^2ba^2 & (n_2-n_3\geq 3, n_3-n_4\leq 0), \\
(ba)^{-n_1+n_2}ba^2(ba)^{-n_2+n_3}ba^2(ba)^{-n_3+n_4}ba^2ba^2 & (n_2-n_3\leq 0, n_3-n_4\leq 0).
\end{cases} 
\end{align*}
}
This would not be equal to $x_1^n$ by (\ref{x_1^n}) and Theorem \ref{key thm of B}. 
Thus we obtain $n_1-n_3\geq 3$. 
By operating a map $t$ to both sides of Equation (\ref{r=4}), we obtain: 
\[
x_2x_1^{n_3-n_4}x_2x_1^{n_2-n_3}x_2x_1^{n_1-n_2}x_2=x_1^n. 
\]
So we obtain $n_3-n_4\geq 3$ by the same argument above. 
If $n_2-n_3\leq 0$, the left side of Equation (\ref{r=4}) is represented by $a$ and $b$ as follows: 
{\allowdisplaybreaks
\begin{align*}
& x_2x_1^{n_1-n_2}x_2x_1^{n_2-n_3}x_2x_1^{n_3-n_4}x_2 \\
= & ba^2\cdot a(ba^2)^{n_1-n_2-1}ba\cdot ba^2\cdot a^2(ba)^{-n_2+n_3-1}ba^2\cdot ba^2 \cdot a(ba^2)^{n_3-n_4-1}ba\cdot ba^2 \\
= & a^2(ba^2)^{n_1-n_2-1}(ba)^{-n_2+n_3+1}ba(ba^2)^{n_3-n_4-2}ba\cdot ba^2. 
\end{align*}
}
This would not be equal to $x_1^n$ and we obtain $n_2-n_3\geq 3$. 
Thus the left side of Equation (\ref{r=4}) is represented by $a$ and $b$ as follows: 
{\allowdisplaybreaks
\begin{align*}
& x_2x_1^{n_1-n_2}x_2x_1^{n_2-n_3}x_2x_1^{n_3-n_4}x_2 \\
= & ba^2\cdot a(ba^2)^{n_1-n_2-1}ba\cdot ba^2\cdot a(ba^2)^{n_2-n_3-1}ba\cdot ba^2 \cdot a(ba^2)^{n_3-n_4-1}ba\cdot ba^2 \\
= & a^2(ba^2)^{n_1-n_2-3}ba(ba^2)^{n_2-n_3-3}a^2(ba^2)^{n_3-n_4-3}baba^2 \\
= & \begin{cases}
a^2(ba^2)^{n_1-n_2-4}ba(ba^2)^{n_3-n_4-4}baba^2 & (n_2-n_3=3), \\
a^2(ba^2)^{n_1-n_2-3}ba(ba^2)^{n_2-n_3-4}ba(ba^2)^{n_3-n_4-3}baba^2 & (n_2-n_3\geq 4). \\
\end{cases}
\end{align*}
}
By Theorem \ref{key thm of B}, we obtain: 
\[
(n_1-n_2,n_2-n_3,n_3-n_4)=(3,4,3)\text{ or }(4,3,4). 
\]
Then we can easily change $W_f$ into $T_4$ by elementary transformations and simultaneous conjugations. 
Thus $M$ is diffeomorphic to $\sharp 3\mathbb{CP}^2\sharp\overline{\mathbb{CP}^2}$. 

\par

For the case $r=5$, $W_f$ is equal to $S_5$, $S_4T(n_1)$, $S_3T(n_1,n_2)$, $S_2T(n_1,n_2,n_3)$, $S_1T(n_1,n_2,n_3,n_4)$ or $T(n_1,n_2,n_3,n_4,n_5)$. 
But all we need to examine is the case $W_f=T(n_1,n_2,n_3,n_4,n_5)$ as the case $r=4$. 
We can assume that $n_i-n_{i+1}\neq 1,2$ for $i=1,2,3,4$ and $(n_j-n_{j+1},n_{j+1}-n_{j+2})\neq(3,3)$ for $j=1,2,3$ by the same argument as in the case $r=4$. 
If $(n_i-n_{i+1},n_{i+1}-n_{i+2},n_{i+2}-n_{i+3})=(3,4,3)$ or $(4,3,4)$ for $i=1$ or $2$, then we can change $W_f$ into $T(\tilde{n})T_4$ by elementary transformations and simultaneous conjugations, where $\tilde{n}$ is an integer. 
This contradicts $w(W_f)=\pm X_1^m$. 
So we obtain $(n_i-n_{i+1},n_{i+1}-n_{i+2},n_{i+2}-n_{i+3})=(3,4,3), (4,3,4)$ for $i=1,2$. 
Since $w(W_f)=\pm X_1^m$, we obtain the following equation: 
{\allowdisplaybreaks
\begin{align}
\label{r=5}
& x_1^{-n_1}x_2x_1^{n_1}x_1^{-n_2}x_2x_1^{n_2}x_1^{-n_3}x_2x_1^{n_3}x_1^{-n_4}x_2x_1^{n_4}x_1^{-n_5}x_2x_1^{n_5}=x_1^m \nonumber \\
\Rightarrow & x_2x_1^{n_1-n_2}x_2x_1^{n_2-n_3}x_2x_1^{n_3-n_4}x_2x_1^{n_4-n_5}x_2=x_1^n, 
\end{align}
}
where $n=m+n_1-n_5$. 
If $n_1-n_2\leq 0$ and $n_2-n_3\leq 0$, the left side of Equation (\ref{r=5}) is represented by $a$ and $b$ as follows: 
{\allowdisplaybreaks
\begin{align*}
& x_2x_1^{n_1-n_2}x_2x_1^{n_2-n_3}x_2x_1^{n_3-n_4}x_2x_1^{n_4-n_5}x_2 \\
= & \begin{cases}
ba^2\cdot a^2(ba)^{-n_1+n_2-1}ba^2\cdot ba^2\cdot a^2(ba)^{-n_2+n_3-1}ba^2\cdot ba^2 \\
\cdot a(ba^2)^{n_3-n_4-1}ba\cdot ba^2\cdot a(ba^2)^{n_4-n_5-1}ba\cdot ba^2 & (n_3-n_4\geq 3, n_4-n_5\geq 3) \\
ba^2\cdot a^2(ba)^{-n_1+n_2-1}ba^2\cdot ba^2\cdot a^2(ba)^{-n_2+n_3-1}ba^2\cdot ba^2 \\
\cdot a^2(ba)^{-n_3+n_4-1}ba^2\cdot ba^2\cdot a(ba^2)^{n_4-n_5-1}ba\cdot ba^2 & (n_3-n_4\leq 0, n_4-n_5\geq 3) \\
ba^2\cdot a^2(ba)^{-n_1+n_2-1}ba^2\cdot ba^2\cdot a^2(ba)^{-n_2+n_3-1}ba^2\cdot ba^2 \\
\cdot a(ba^2)^{n_3-n_4-1}ba\cdot ba^2\cdot a^2(ba)^{-n_4+n_5-1}ba^2\cdot ba^2 & (n_3-n_4\geq 3, n_4-n_5\leq 0) \\
ba^2\cdot a^2(ba)^{-n_1+n_2-1}ba^2\cdot ba^2\cdot a^2(ba)^{-n_2+n_3-1}ba^2\cdot ba^2 \\
\cdot a^2(ba)^{-n_3+n_4-1}ba^2\cdot ba^2\cdot a^2(ba)^{-n_4+n_5-1}ba^2\cdot ba^2 & (n_3-n_4\leq 0, n_4-n_5\leq 0) \\
\end{cases} \\
= & \begin{cases}
(ba)^{-n_1+n_2}ba^2(ba)^{-n_2+n_3}ba(ba^2)^{n_3-n_4-3}ba(ba^2)^{n_4-n_5-3}baba^2 & (n_3-n_4\geq 3, n_4-n_5\geq 3), \\
(ba)^{-n_1+n_2}ba^2(ba)^{-n_2+n_3}ba^2(ba)^{-n_3+n_4}ba(ba^2)^{n_4-n_5-2}baba^2 & (n_3-n_4\leq 0, n_4-n_5\geq 3), \\
(ba)^{-n_1+n_2}ba^2(ba)^{-n_2+n_3}ba(ba^2)^{n_3-n_4-2}ba(ba)^{-n_4+n_5}ba^2ba^2 & (n_3-n_4\geq 3, n_4-n_5\leq 0), \\
(ba)^{-n_1+n_2}ba^2(ba)^{-n_2+n_3}ba^2(ba)^{-n_3+n_4}ba^2(ba)^{-n_4+n_5}ba^2ba^2 & (n_3-n_4\leq 0, n_4-n_5\leq 0). 
\end{cases}
\end{align*}
}
This would not be equal to $x_1^n$ by Theorem \ref{key thm of B}. 
So either of the cases $n_1-n_2\geq 3$ and $n_2-n_3\geq 3$ holds. 
If $n_1-n_2\leq 0$, then $n_2-n_3\geq 3$ and the left side of Equation (\ref{r=5}) is represented by $a$ and $b$ as follows: 
{\allowdisplaybreaks
\begin{align*}
& x_2x_1^{n_1-n_2}x_2x_1^{n_2-n_3}x_2x_1^{n_3-n_4}x_2x_1^{n_4-n_5}x_2 \\
= & \begin{cases}
ba^2\cdot a^2(ba)^{-n_1+n_2-1}ba^2\cdot ba^2\cdot a(ba^2)^{n_2-n_3-1}ba\cdot ba^2 \\
\cdot a(ba^2)^{n_3-n_4-1}ba\cdot ba^2\cdot a(ba^2)^{n_4-n_5-1}ba\cdot ba^2 & (n_3-n_4\geq 3, n_4-n_5\geq 3) \\
ba^2\cdot a^2(ba)^{-n_1+n_2-1}ba^2\cdot ba^2\cdot a(ba^2)^{n_2-n_3-1}ba\cdot ba^2 \\
\cdot a^2(ba)^{-n_3+n_4-1}ba^2\cdot ba^2\cdot a(ba^2)^{n_4-n_5-1}ba\cdot ba^2 & (n_3-n_4\leq 0, n_4-n_5\geq 3) \\
ba^2\cdot a^2(ba)^{-n_1+n_2-1}ba^2\cdot ba^2\cdot a(ba^2)^{n_2-n_3-1}ba\cdot ba^2 \\
\cdot a(ba^2)^{n_3-n_4-1}ba\cdot ba^2\cdot a^2(ba)^{-n_4+n_5-1}ba^2\cdot ba^2 & (n_3-n_4\geq 3, n_4-n_5\leq 0) \\
ba^2\cdot a^2(ba)^{-n_1+n_2-1}ba^2\cdot ba^2\cdot a(ba^2)^{n_2-n_3-1}ba\cdot ba^2 \\
\cdot a^2(ba)^{-n_3+n_4-1}ba^2\cdot ba^2\cdot a^2(ba)^{-n_4+n_5-1}ba^2\cdot ba^2 & (n_3-n_4\leq 0, n_4-n_5\leq 0) \\
\end{cases} \\
= & \begin{cases}
(ba)^{-n_1+n_2}ba(ba^2)^{n_2-n_3-3}ba(ba^2)^{n_3-n_4-3}a^2(ba^2)^{n_4-n_5-3}baba^2 & (n_3-n_4\geq 3, n_4-n_5\geq 3) \\
(ba)^{-n_1+n_2}ba(ba^2)^{n_2-n_3-2}(ba)^{-n_3+n_4+1}ba(ba^2)^{n_4-n_5-2}baba^2 & (n_3-n_4\leq 0, n_4-n_5\geq 3) \\
(ba)^{-n_1+n_2}ba(ba^2)^{n_2-n_3-3}ba(ba^2)^{n_3-n_4-3}ba(ba)^{-n_4+n_5}ba^2ba^2 & (n_3-n_4\geq 3, n_4-n_5\leq 0) \\
(ba)^{-n_1+n_2}ba(ba^2)^{n_2-n_3-2}(ba)^{-n_3+n_4+1}ba^2(ba)^{-n_4+n_5}ba^2ba^2 & (n_3-n_4\leq 0, n_4-n_5\leq 0) \\
\end{cases} \\
= & \begin{cases}
(ba)^{-n_1+n_2}ba(ba^2)^{n_2-n_3-3}ba(ba^2)^{n_3-n_4-4}ba(ba^2)^{n_4-n_5-3}baba^2 & (n_3-n_4\geq 4, n_4-n_5\geq 3), \\
(ba)^{-n_1+n_2}ba(ba^2)^{n_2-n_3-4}ba(ba^2)^{n_4-n_5-4}baba^2 & (n_3-n_4=3, n_4-n_5\geq 3), \\
(ba)^{-n_1+n_2}ba(ba^2)^{n_2-n_3-2}(ba)^{-n_3+n_4+1}ba(ba^2)^{n_4-n_5-2}baba^2 & (n_3-n_4\leq 0, n_4-n_5\geq 3), \\
(ba)^{-n_1+n_2}ba(ba^2)^{n_2-n_3-3}ba(ba^2)^{n_3-n_4-3}ba(ba)^{-n_4+n_5}ba^2ba^2 & (n_3-n_4\geq 3, n_4-n_5\leq 0), \\
(ba)^{-n_1+n_2}ba(ba^2)^{n_2-n_3-2}(ba)^{-n_3+n_4+1}ba^2(ba)^{-n_4+n_5}ba^2ba^2 & (n_3-n_4\leq 0, n_4-n_5\leq 0). \\
\end{cases}
\end{align*}
}
This would not be equal to $x_1^n$ by Theorem \ref{key thm of B}. 
So we obtain $n_1-n_2\geq 3$. 
By operating a map $t$ to both sides of Equation (\ref{r=5}), we obtain: 
\begin{align}
\label{t:r=5}
x_2x_1^{n_4-n_5}x_2x_1^{n_3-n_4}x_2x_1^{n_2-n_3}x_2x_1^{n_1-n_2}x_2=x_1^n
\end{align}
So we obtain $n_4-n_5\geq 3$. 
If $n_2-n_3\leq 0$, the left side of Equation (\ref{r=5}) is represented by $a$ and $b$ as follows: 
{\allowdisplaybreaks
\begin{align*}
& x_2x_1^{n_1-n_2}x_2x_1^{n_2-n_3}x_2x_1^{n_3-n_4}x_2x_1^{n_4-n_5}x_2 \\
= & \begin{cases}
ba^2\cdot a(ba^2)^{n_1-n_2-1}ba\cdot ba^2\cdot a^2(ba)^{-n_2+n_3-1}ba^2\cdot ba^2 \\
\cdot a(ba^2)^{n_3-n_4-1}ba\cdot ba^2\cdot a(ba^2)^{n_4-n_5-1}ba\cdot ba^2 & (n_3-n_4\geq 3) \\
ba^2\cdot a(ba^2)^{n_1-n_2-1}ba\cdot ba^2\cdot a^2(ba)^{-n_2+n_3-1}ba^2\cdot ba^2 \\
\cdot a^2(ba)^{-n_3+n_4-1}ba^2\cdot ba^2\cdot a(ba^2)^{n_4-n_5-1}ba\cdot ba^2 & (n_3-n_4\leq 0) \\
\end{cases} \\
= & \begin{cases}
a^2ba^2(ba^2)^{n_1-n_2-3}(ba)^{-n_2+n_3+1}ba(ba^2)^{n_3-n_4-3}ba(ba^2)^{n_4-n_5-3}baba^2 & (n_3-n_4\geq 3), \\
a^2ba^2(ba^2)^{n_1-n_2-3}(ba)^{-n_2+n_3+1}ba^2(ba)^{-n_3+n_4}ba(ba^2)^{n_4-n_5-2}baba^2 & (n_3-n_4\leq 0). \\
\end{cases}
\end{align*}
}
This would not be equal to $x_1^n$ by Theorem \ref{key thm of B} and we obtain $n_2-n_3\geq 3$. 
We also obtain $n_3-n_4\geq 3$ by Equation (\ref{t:r=5}). 
Thus the left side of Equation (\ref{r=5}) is represented by $a$ and $b$ as follows: 
{\allowdisplaybreaks
\begin{align*}
& x_2x_1^{n_1-n_2}x_2x_1^{n_2-n_3}x_2x_1^{n_3-n_4}x_2x_1^{n_4-n_5}x_2 \\
= & ba^2\cdot a(ba^2)^{n_1-n_2-1}ba\cdot ba^2\cdot a(ba^2)^{n_2-n_3-1}ba\cdot ba^2\cdot a(ba^2)^{n_3-n_4-1}ba\cdot ba^2\cdot a(ba^2)^{n_4-n_5-1}ba\cdot ba^2 \\
= & a^2(ba^2)^{n_1-n_2-2}a^2(ba^2)^{n_2-n_3-3}a^2(ba^2)^{n_3-n_4-3}a^2(ba^2)^{n_4-n_5-3}baba^2 \\
= & a^2(ba^2)^{n_1-n_2-3}ba(ba^2)^{n_2-n_3-3}a^2(ba^2)^{n_3-n_4-3}a^2(ba^2)^{n_4-n_5-3}baba^2 \\
= & \begin{cases}
a^2(ba^2)^{n_1-n_2-3}ba(ba^2)^{n_2-n_3-4}ba(ba^2)^{n_3-n_4-3}a^2(ba^2)^{n_4-n_5-3}baba^2 & (n_2-n_3\geq 4) \\
a^2(ba^2)^{n_1-n_2-4}ba(ba^2)^{n_3-n_4-5}ba(ba^2)^{n_4-n_5-3}baba^2  & (n_2-n_3=3) \\
\end{cases} \\
= & \begin{cases}
a^2(ba^2)^{n_1-n_2-3}ba(ba^2)^{n_2-n_3-4}ba(ba^2)^{n_3-n_4-4}ba(ba^2)^{n_4-n_5-3}baba^2 & (n_2-n_3\geq 4, n_3-n_4\geq 4), \\
a^2(ba^2)^{n_1-n_2-3}ba(ba^2)^{n_2-n_3-5}ba(ba^2)^{n_4-n_5-4}baba^2 & (n_2-n_3\geq 4, n_3-n_4=3), \\
a^2(ba^2)^{n_1-n_2-4}ba(ba^2)^{n_3-n_4-5}ba(ba^2)^{n_4-n_5-3}baba^2  & (n_2-n_3=3). \\
\end{cases}
\end{align*}
}
By Theorem \ref{key thm of B}, we obtain: 
\[
(n_1-n_2,n_2-n_3,n_3-n_4,n_4-n_5)=(3,4,4,3)\text{, }(3,5,3,4)\text{ or }(4,3,5,3). 
\]
Then we can easily change $W_f$ into $T_5$ by elementary transformations and simultaneous conjugations. 
Thus $M$ is diffeomorphic to $\sharp 4\mathbb{CP}^2\sharp\overline{\mathbb{CP}^2}$. 

\par

Combining the above arguments, we complete the proof of Main Theorem B.  \hfill $\square$

\begin{rem}

The case $r\geq 6$ is algebraically too complicated to solve the classification problem. 
For example, in the case $r=6$, the following sequence appears as a candidate of a Hurwitz system of a SBLF: 
\[
\tilde{T_6}=(X_1^{-10}X_2X_1^{10},X_1^{-6}X_2X_1^{6},X_1^{-1}X_2X_1,X_1^{2}X_2X_1^{-2},X_1^{6}X_2X_1^{-6},X_1^{11}X_2X_1^{-11}). 
\]
Indeed, it is easy to see that $w(\tilde{T_6})=-X_1^{-24}$ and so it satisfies the condition in Theorem \ref{main in chart sec}. 
We can change $\tilde{T_6}$ into the following sequence by successive application of elementary transformations and simultaneous conjugations: 
\[
(T_{3,1},T_{3,2},T_{2,3},T_{2,5},T_{1,4},T_{-1,1}). 
\]
We can draw a Kirby diagram of the total space of a SBLF whose Hurwitz systems $\tilde{T_6}$. 
We can prove by Kirby calculus that this $4$-manifold is diffeomorphic to $\sharp 5\mathbb{CP}^2\sharp\overline{\mathbb{CP}^2}$, 
which is the total space of a SBLF whose Hurwitz system is $T_6$. 
However, the author does not know whether $\tilde{T_6}$ can be changed into $T_6$ by elementary transformations and simultaneous conjugations. 

\end{rem}

\begin{conj}

If $M$ admits a genus-$1$ SBLF structure with non-empty round singular locus, then it is diffeomorphic to one of the following $4$-manifolds: 

\begin{itemize}

\item $\sharp k\mathbb{CP}^2\sharp l\overline{\mathbb{CP}^2}$, where $l>0$ and $k\geq 0$. 

\item $\sharp kS^2\times S^2$, where $k\geq 0$ . 

\item $S^1\times S^3\sharp S\sharp k\overline{\mathbb{CP}^2}$, where $k\geq 0$ and $S$ is either $S^2\times S^2$ or $S^2\tilde{\times}S^2$. 

\item $L\sharp k\overline{\mathbb{CP}^2}$, where $k\geq 0$ and $L$ is either $L_n$ or $L_n^{\prime}$. 

\end{itemize}

More strongly, the families of genus-$1$ SBLF obtained in Section \ref{examples} contain all genus-$1$ SBLFs with non-empty round singular locus. 

\end{conj}

\end{document}